\documentclass[12pt]{article}
\usepackage{amsmath,amssymb}
\usepackage{graphicx}

\newtheorem{theorem}{Theorem}[section]
\newtheorem{proposition}{Proposition}[section]
\newtheorem{corollary}{Corollary}[section]

\newcommand{\CQFD}{\nolinebreak\hfill\rule{2mm}{2mm}\medbreak\par}  
\newtheorem{Rem}{Remark}[section]  
  
\newtheorem{Example}{Example}[section]

\numberwithin{equation}{section}

\def\bbe{\mathbb E}

\def\bbp{\mathbb P}
\def\bbr{\mathbb R}

\begin{document}

\title{On the Limiting Shape of Markovian Random Young Tableaux}
\author{Christian Houdr\'e \thanks{Georgia Institute of Technology, 
School of Mathematics, Atlanta, Georgia, 30332-0160,
houdre@math.gatech.edu} 
\and Trevis J. Litherland \thanks{Georgia Institute of Technology, 
School of Mathematics, Atlanta, Georgia, 30332-0160,
trevisl@math.gatech.edu} }

\maketitle

\vspace{0.5cm}

\begin{abstract}
\noindent Let $(X_n)_{n \ge 0}$ 
be an irreducible, aperiodic, homogeneous
Markov chain, with state space 
an ordered finite alphabet of size $m$.
Using combinatorial constructions and 
weak invariance principles, we obtain
the limiting shape of the associated Young tableau
as a multidimensional Brownian functional.  Since the 
length of the top row of
the Young tableau is also
the length of the longest (weakly) increasing
subsequence of $(X_k)_{1\le k \le n}$, 
the corresponding limiting law follows.
We relate our results to a conjecture of
Kuperberg by showing that, under a cyclic condition,
a spectral characterization of the Markov
transition matrix delineates precisely when
the limiting shape is the spectrum of the
traceless GUE.  For $m=3$, all cyclic Markov chains 
have such a limiting shape,
a fact previously known for $m=2$. 
However, this is no longer true for $m \ge 4$.

\end{abstract}





\section{Introduction}  

The identification of the limiting distribution of $LI_n$,
the length of the longest increasing
subsequence  of a random word of length $n$,
whose letters are iid and chosen uniformly from an ordered, $m$-letter alphabet,
was first made by Tracy and Widom \cite{TW}.
They showed that the limiting distribution of $LI_n$,
properly centered and normalized,
is that of the largest eigenvalue
of the traceless $m \times m$ 
Gaussian unitary ensemble (GUE). 
In the non-uniform iid case,
Its, Tracy, and Widom \cite{ITW1,ITW2}
described the corresponding
limiting distribution as 
that of the largest eigenvalue of one
of the diagonal blocks 
(corresponding to the highest probability)
in a direct sum of certain independent GUE matrices.  
The number and respective dimensions of these matrices
are determined by the multiplicities 
of the probabilities of choosing the
letters, and the direct sum is 
subject again to an overall zero-trace type of condition.

  The well-known Robinson-Schensted-Knuth (RSK)
correspondence between sequences and pairs of Young
tableaux led Tracy and Widom \cite{TW} to conjecture that
the (necessarily $m$-row) Young tableau of a
random word generated by an $m$-letter, uniform iid sequence
has a limiting shape given by the joint distribution 
of the eigenvalues of a
$m \times m$ traceless element of the GUE.
Since the length of the longest row of the 
Young tableau is precisely $LI_n$,
this appears to be a natural generalization.
Johansson \cite{Jo} proved this conjecture 
using orthogonal polynomial methods.
Further, Okounkov \cite{Ok}, and
Borodin, Okounkov, and Olshankii \cite{BOO},
as well as Johansson \cite{Jo}, also answered a conjecture
of Baik, Deift, and Johansson \cite{BDJ2,BDJ2Add} regarding the limiting
shape of the Young tableau associated to a random
permutation of $\{1,2,\dots,n\}$.
In particular, as $n$ grows without bound,
the lengths $\lambda_1, \lambda_2, \dots ,\lambda_k$
of the first $k$ rows of the Young tableau,
appropriately centered and scaled, have the same limiting law 
as the $k$ largest eigenvalues of a $n \times n$ 
element of the GUE,
a result first proved, for $k=2$, in \cite{BDJ2,BDJ2Add}.

The extension to the non-uniform iid case
was addressed to some degree in 
Its, Tracy, and Widom \cite{ITW1,ITW2},
who focused primarily on the top row
of the Young tableau.
Here the obvious conjecture is that
the limiting shape has rows 
whose suitably centered and normalized lengths
have a joint distribution
which is that of the whole spectrum of 
the direct sum of GUE matrices
described above.  Below, we prove this result as a special
case of the Markovian framework.

  Kuperberg \cite{Ku} conjectured that
if the word is generated by an irreducible,
doubly-stochastic, cyclic Markov chain,
then the limiting distribution of the shape is still
that of the joint distribution of the eigenvalues
of a traceless $m \times m$ element of the GUE.  
For $m=2$, this was shown to be true 
by Chistyakov and G\"otze \cite{ChG},
who, in view of further simulations,
expressed doubts concerning the validity for $m \ge 4$.
For $m=3$, we will show that the conjecture
holds as well. 
However, for $m\ge4$, this is no longer the case.
Indeed, some, but not all, cyclic Markov chains 
lead to a limiting law as in the iid uniform
case already obtained by Johansson \cite{Jo}.

  The precise class of homogeneous Markov chains 
with which Kuperberg's conjecture is concerned is more specific than
the ones we shall study.  The irreducibility of
the chain is a basic property we certainly must demand:
each letter has to occur at some point following the
occurrence of any given letter.
Moreover, the doubly-stochastic hypothesis
ensures that we have a uniform stationary distribution.
However, the cyclic criterion, {\it i.e.}, the Markov transition
matrix $P$ has entries satisfying $p_{i,j} = p_{i+1,j+1}$,
for $1 \le i,j \le m$ (where $m+1 = 1$), 
is more restrictive: cyclicity implies but
is not equivalent to $P$ being doubly stochastic.  Kuperberg was led
to introduce this latter restriction via simulations \cite{Ku}
inspired by mathematical physics considerations, 
which appear to show that at least some irreducible, 
doubly-stochastic, {\it non-}cyclic Markov chains do
{\it not} produce such limiting behavior.

Let us also note that Kuperberg implicitly assumes 
the Markov chain to also be aperiodic.  Indeed, the simple
$2$-state Markov chain for the letters $\alpha_1$ and $\alpha_2$
described by $\bbp(X_{n+1} = \alpha_i| X_n = \alpha_j) = 1$ for $i \ne j$,
produces a sequence of alternating letters, so that
$LI_n$ is always either $n/2$ or $n/2 + 1$, for $n$ even,
and $(n+1)/2$, for $n$ odd, and so has a degenerate
limiting distribution.  Even though this Markov chain
is irreducible, doubly-stochastic, and cyclic,
it is periodic.

The paper is organized in the following manner.
In Section $2$, we present the simple combinatorial
formulation of the $LI_n$
problem the authors developed in \cite{HL}.  
Next, in Section $3$, we use this formulation
to rederive the two-letter Markov case first studied by 
Chistyakov and G\"otze \cite{ChG}.
Then, in order to extend these results to alphabets
of size $m \ge 3$, we introduce, in Section $4$,  a slight 
modification of our original combinatorial
development, and so obtain a functional of
combinatorial quantities which describes
the shape of the entire Young tableau with $n$ cells,
along with a concise expression for the
associated asymptotic covariance structure.  
In Section $5$, we apply
Markovian Invariance Principles
to express the limiting shape of the Young 
tableau as a Brownian functional
for all irreducible, aperiodic, homogeneous
Markov chains (without the cyclic or even
the doubly-stochastic constraint.)
Using this functional we are then able to
answer Kuperberg's conjecture. 
In Section $6$, we investigate, in
further detail, various symmetries exhibited
by the Brownian functional.  
In particular, we clarify the asymptotic covariance
structure in the cyclic case, and
obtain, for $m$ arbitrary, a precise description of
the class of cyclic Markov chains having
the same limiting law as in the uniform iid case.
In Section $7$, we further explore connections between the
various Brownian functionals obtained as limiting laws
and eigenvalues of random matrices.
Finally, in Section $8$, we conclude with a brief discussion
of natural extensions and complements to some of the ideas
and results presented in the paper.

\section{Combinatorics}

As in \cite{HL}, one can express 
$LI_n$ in a fundamentally combinatorial manner.
For convenience, this section recapitulates that development.\\

Let $(X_n)_{n\ge1}$ consist of a sequence of values taken from an
$m$-letter ordered alphabet, 
${\cal A}_m = \{\alpha_1 < \alpha_2 < \cdots < \alpha_m\}$. 
Let $a^r_k$ be the number of occurrences
of $\alpha_r$ among $(X_i)_{1\le i \le k}$.  
Each increasing subsequence of $(X_i)_{1\le i \le k}$ consists simply of 
consecutive identical values, 
with these values forming an increasing subsequence of $\alpha_r$.
Moreover, the number of occurrences of $\alpha_r\in \{\alpha_1,
\dots,\alpha_m\}$ among
$(X_i)_{k+1 \le i \le \ell}$,
where $1 \le k < \ell \le n$, is simply
$a^r_{\ell}-a^r_k$. The length of the longest increasing subsequence 
of $X_1,X_2,\dots, X_n$ is thus given by

\begin{equation}\label{item1}
LI_n=\max_{\stackrel{\scriptstyle 0\le k_1\le\cdots}{\le k_{m-1}\le n}}
     [(a^1_{k_1}-a^1_0)+(a^2_{k_2}-a^2_{k_1})+\cdots +
     (a^m_n-a^m_{k_{m-1}})],
\end{equation}

\noindent {\it i.e.},

\begin{equation}\label{item2}
LI_n=\max_{\stackrel{\scriptstyle 0\le k_1\le\cdots}{\le k_{m-1}\le n}}
     [(a^1_{k_1}-a^2_{k_1})+(a^2_{k_2}-a^3_{k_2})+\cdots +
     (a^{m-1}_{k_{m-1}}-a^m_{k_{m-1}})+a^m_n],
\end{equation}

\noindent where $a^r_0=0$.

For $i = 1, \dots ,n$ and $r = 1, \dots ,{m-1}$, let 

\begin{equation}\label{item3}
Z^r_i=	\begin{cases} 	1, &\text{if $X_i=\alpha_r,$}\\
			-1, & \text{if $X_i=\alpha_{r+1},$}\\
			0, & \text{otherwise,}
	\end{cases}
\end{equation}

\noindent and let $S^r_k=\sum^k_{i=1}Z^r_i$, $k = 1, \dots ,n$, with also $S^r_0=0$. 
Then clearly $S^r_k = a^r_k - a^{r+1}_k$.  Hence,

\begin{equation}\label{item4}
LI_n=\max_{\stackrel{\scriptstyle 0\le k_1\le\cdots}{\le k_{m-1}\le n}}
     \{S^1_{k_1}+S^2_{k_2}+ \cdots + S^{m-1}_{k_{m-1}}+a^m_n\}.
\end{equation}

By the telescoping nature of the sum $\sum_{k=r}^{m-1}S_n^k = \sum_{k=r}^{m-1}(a_n^k - a_n^{k+1})$,
we find that, for each $1 \le r \le m-1$,
$a_n^r = a_n^m + \sum_{k=r}^{m-1}S_n^k$.
Since $a^1_k, \dots ,a^m_k$ must evidently sum to $k$, we have

\begin{align*}
n &= \sum^m_{r=1}a^r_n \\
&= \sum^{m-1}_{r=1} \left(a_n^m + \sum_{k=r}^{m-1}S_n^k\right) + a^m_n\\
&= \sum^{m-1}_{r=1}rS^r_n + ma^m_n.\\
\end{align*}

Solving for $a^m_n$ gives us
$$a^m_n=\frac nm- \frac1m \sum^{m-1}_{r=1} rS^r_n.$$

Substituting into \eqref{item4}, we finally obtain

\begin{equation}\label{item5}
 LI_n=\frac nm-\frac1m\sum^{m-1}_{r=1} rS^r_n+
      \max_{\stackrel{\scriptstyle 0\le k_1\le\cdots}{\le k_{m-1}\le n}}
      \{S^1_{k_1}+S^2_{k_2}+ \cdots + S^{m-1}_{k_{m-1}}\}.
\end{equation}

As was emphasized in \cite{HL} \eqref{item5} is of a {\it purely
combinatorial nature or, in more probabilistic terms, is of a pathwise nature}.  
We now proceed to analyze \eqref{item5} in the case
of a Markovian sequence.

\section{Markovian Alphabet: $2$-Letter Case}

We begin our study of Markovian alphabets by
concentrating on the $2$-letter case.  Here
$(X_n)_{n\ge0}$ is described by the following
transition probabilities between the two states
(which we identify with the two letters 
$\alpha_1$ and $\alpha_2$):
$\bbp(X_{n+1} = \alpha_2 | X_n = \alpha_1) = a$ and
$\bbp(X_{n+1} = \alpha_1 | X_n = \alpha_2) = b$,
where $0 < a + b < 2$.  We later examine the degenerate
cases $a = b = 0$ and $a = b = 1$.
In keeping with the common usage within the Markov chain
literature, we begin our sequence at $n=0$, although
our focus will be on $n \ge 1$.
Denoting by $(p_n^1,p_n^2)$
the vector describing the probability distribution
on $\{\alpha_1,\alpha_2\}$ at time $n$, 
we have

\begin{equation}\label{item5a}
  \begin{pmatrix} p_{n+1}^1, p_{n+1}^2 \end{pmatrix}
= \begin{pmatrix} p_{n}^1, p_{n}^2   \end{pmatrix}
  \begin{pmatrix} 1-a & a\\ b & 1-b   \end{pmatrix}.
\end{equation}

The eigenvalues of the
matrix in \eqref{item5a} are $\lambda_1 = 1$ and 
$-1 < \lambda_2 = 1 - a - b < 1$, 
with respective left eigenvectors
$(\pi_1, \pi_2) = (b/(a+b),a/(a+b))$ and
$(1,-1)$.  Moreover, $(\pi_1, \pi_2)$ is also the
stationary distribution.  Given any initial
distribution $(p_0^1,p_0^2)$, we find that

\begin{equation}\label{item5b}
\begin{pmatrix} p_{n}^1, p_{n}^2 \end{pmatrix}
=  \begin{pmatrix} \pi_1, \pi_2  \end{pmatrix}
+ \lambda_2^{n}\frac{a p_0^1-b p_0^2}{a+b} 
  \begin{pmatrix} 1,-1\end{pmatrix}
  \rightarrow \begin{pmatrix} \pi_1, \pi_2  \end{pmatrix},
\end{equation}

\noindent as $n \rightarrow \infty$, since $\lambda_2 < 1$.\\

Our goal is now to use these probabilistic
expressions to describe the random variables
$Z_k^1$ and $S_k^1$ defined in the previous section.  
(We retain the redundant
superscript ``$1$'' in $Z_k^1$ and $S_k^1$ in the 
interest of uniformity.)

Setting $\beta = a p_0^1-b p_0^2$,
we easily find that 

\begin{align}\label{item5c}
\bbe Z_k^1 &= (+1)\left(\pi_1 + \frac{\beta}{a+b}\lambda_2^{k} \right)
            + (-1)\left(\pi_2 - \frac{\beta}{a+b}\lambda_2^{k} \right)\nonumber\\
	   &= \frac{b-a}{a+b}  + 2 \frac{\beta}{a+b}\lambda_2^{k},
\end{align}

\noindent for each $1 \le k \le n$.  Thus,

\begin{equation}\label{item5d}
\bbe S_k^1 = \frac{b-a}{a+b}k  
           + 2 \left(\frac{\beta \lambda_2}{a+b}\right) \left( \frac{1-\lambda_2^k}{1-\lambda_2} \right),
\end{equation}

\noindent and so $\bbe S_k^1/k \rightarrow (b-a)/(a+b)$,
as $k \rightarrow \infty$.

Turning to the second moments of 
$Z_k^1$ and $S_k^1$, first note
that $\bbe (Z_k^1)^2 = 1$, since
$(Z_k^1)^2 = 1$ a.s.
Next, we consider $\bbe Z_k^1 Z_{\ell}^1$,
for $k < \ell$.  Using
the Markovian structure of $(X_n)_{n\ge0}$,
it quickly follows that

\begin{align}\label{item5e}
&\bbp((X_k,X_{\ell}) = (x_k,x_{\ell}))\nonumber\\
 &\qquad =
\begin{cases}
  \left(\pi_1 + \lambda_2^{\ell-k}\frac{a}{a+b}\right)\left(\pi_1 + \lambda_2^{k}\frac{\beta}{a+b}\right),
     &\text{if $(x_k,x_{\ell})=(\alpha_1,\alpha_1)$},\\
  \left(\pi_1 - \lambda_2^{\ell-k}\frac{b}{a+b}\right)\left(\pi_2 - \lambda_2^{k}\frac{\beta}{a+b}\right),
     &\text{if $(x_k,x_{\ell})=(\alpha_1,\alpha_2)$},\\
  \left(\pi_2 - \lambda_2^{\ell-k}\frac{a}{a+b}\right)\left(\pi_1 + \lambda_2^{k}\frac{\beta}{a+b}\right),
     &\text{if $(x_k,x_{\ell})=(\alpha_2,\alpha_1)$},\\
  \left(\pi_2 + \lambda_2^{\ell-k}\frac{b}{a+b}\right)\left(\pi_2 - \lambda_2^{k}\frac{\beta}{a+b}\right),
     &\text{if $(x_k,x_{\ell})=(\alpha_2,\alpha_2)$}.
\end{cases}
\end{align}

For simplicity, we will henceforth
assume that our initial
distribution is the stationary one, {\it i.e.},
$(p_0^1,p_0^2) = (\pi_1, \pi_2)$.
Later, (see Concluding Remarks) we drop this assumption 
and deal with initial distributions 
concentrated on an arbitrary state.
Under this assumption, $\beta = 0$, 
$\bbe S_k^1 = k\mu$, where 
$\mu = \bbe Z_k^1 = (b-a)/(a+b)$, and 
\eqref{item5e} simplifies to

\begin{align}\label{item5f}
&\bbp((X_k,X_{\ell}) = (x_k,x_{\ell}))\nonumber\\
 &\qquad =
\begin{cases}
  \left(\pi_1 + \lambda_2^{\ell-k}\frac{a}{a+b}\right)\pi_1,
     &\text{if $(x_k,x_{\ell})=(\alpha_1,\alpha_1)$},\\
  \left(\pi_1 - \lambda_2^{\ell-k}\frac{b}{a+b}\right)\pi_2,
     &\text{if $(x_k,x_{\ell})=(\alpha_1,\alpha_2)$},\\
  \left(\pi_2 - \lambda_2^{\ell-k}\frac{a}{a+b}\right)\pi_1,
     &\text{if $(x_k,x_{\ell})=(\alpha_2,\alpha_1)$},\\
  \left(\pi_2 + \lambda_2^{\ell-k}\frac{b}{a+b}\right)\pi_2,
     &\text{if $(x_k,x_{\ell})=(\alpha_2,\alpha_2)$}.
\end{cases}
\end{align}

We can now compute $\bbe Z_k^1 Z_{\ell}^1$:

{\allowdisplaybreaks
\begin{align}\label{item5ga}
\bbe Z_k^1 Z_{\ell}^1 
&= \bbp(Z_k^1 Z_{\ell}^1 = +1) - \bbp(Z_k^1 Z_{\ell}^1 = -1)\nonumber\\
&= \bbp( (X_k,X_{\ell}) \in \{(\alpha_1,\alpha_1),(\alpha_2,\alpha_2) \})\nonumber\\
  &\qquad\qquad - \bbp( (X_k,X_{\ell}) \in \{(\alpha_1,\alpha_2),(\alpha_2,\alpha_1) \})\nonumber\\
&= \left( \pi_1^2 + \lambda_2^{\ell-k}\frac{a}{a+b}\pi_1 + \pi_2^2 + \lambda_2^{\ell-k}\frac{b}{a+b}\pi_2\right)\nonumber\\
  &\qquad\qquad -\left( \pi_1 \pi_2 - \lambda_2^{\ell-k}\frac{b}{a+b}\pi_2 + \pi_1 \pi_2  - \lambda_2^{\ell-k}\frac{a}{a+b}\pi_1\right)\nonumber\\
&= \left(\pi_1^2 + \pi_2^2 + \frac{2ab}{(a+b)^2}\lambda_2^{\ell-k}\right)
  -\left(2\pi_1 \pi_2 - \frac{2ab}{(a+b)^2}\lambda_2^{\ell-k}\right)\nonumber\\
&= \frac{(b-a)^2}{(a+b)^2} + \frac{4ab}{(a+b)^2}\lambda_2^{\ell-k}.
\end{align}
}

Hence, recalling that $\beta = 0$,

\begin{align}\label{item5gb}
\sigma^2 := \mbox{Var} Z_k^1 
&= 1 - \left(\frac{b-a}{a+b}\right)^2\nonumber\\
&= \frac{4ab}{(a+b)^2},
\end{align}

\noindent for all $k\ge 1$,
and, for $k < \ell$, the covariance
of $Z_k^1$ and $Z_{\ell}^1$ is

\begin{align}\label{item5h}
\mbox{Cov} (Z_k^1,Z_{\ell}^1) &= \frac{(b-a)^2}{(a+b)^2} + \sigma^2\lambda_2^{\ell-k} - \left(\frac{b-a}{a+b}\right)^2 = \sigma^2\lambda_2^{\ell-k}.
\end{align}

Proceeding to the covariance structure
of $S_k^1$, we first find that

\begin{align}\label{item5i}
\mbox{Var}S_k^1 &= \sum_{j=1}^k \mbox{Var} Z_j^1 + 2 \sum_{j<\ell}\mbox{Cov}(Z_j^1,Z_l^1)\nonumber\\
&= \sigma^2k + 2\sigma^2\sum_{j<\ell}\lambda_2^{\ell-j}\nonumber\\
&= \sigma^2k 
  + 2\sigma^2\left(\frac{\lambda_2^{k+1} - k\lambda_2^2 + (k-1)\lambda_2}{(1-\lambda_2)^2}\right)\nonumber\\
&= \sigma^2\left(\frac{1+\lambda_2}{1-\lambda_2}\right)k
 + 2\sigma^2 \left(\frac{\lambda_2(\lambda_2^k-1)}{(1-\lambda_2)^2}\right).
\end{align}

Next, for $k < \ell$, and using \eqref{item5h} and \eqref{item5i},
the covariance of $S_k^1$ and $S_{\ell}^1$ is
given by

{\allowdisplaybreaks
\begin{align}
\mbox{Cov} (S_k^1,S_{\ell}^1) 
  &= \sum_{i=1}^k \sum_{j=1}^{\ell} \mbox{Cov} (Z_i^1,Z_j^1)\nonumber\\
  &= \sum_{i=1}^k \mbox{Var} Z_i^1 + 2\sum_{i<j<k} \mbox{Cov} (Z_i^1,Z_j^1) 
      + \sum_{i=1}^k \sum_{j=k+1}^{\ell} \mbox{Cov} (Z_i^1,Z_j^1)\nonumber\\ 
  &= \mbox{Var} S_k^1 + \sum_{i=1}^k \sum_{j=k+1}^{\ell} \mbox{Cov} (Z_i^1,Z_j^1)\nonumber\\ 
  &= \mbox{Var} S_k^1 
    +  \sigma^2\left(\frac{\lambda_2(1-\lambda_2^k)(1-\lambda_2^{\ell-k})}{(1-\lambda_2)^2}\right)\nonumber\\ 
  &= \sigma^2\left( \left(\frac{1+\lambda_2}{1-\lambda_2}\right)k
    - \frac{\lambda_2(1-\lambda_2^k)(1+\lambda_2^{\ell-k})}{(1-\lambda_2)^2}\right).\label{item5jb}
\end{align}
}

From \eqref{item5i} and \eqref{item5jb} we see that,
as $k \rightarrow \infty$,

\begin{equation}\label{item5k}
\frac{\mbox{Var}S_k^1}{k} \rightarrow \sigma^2\left(\frac{1+\lambda_2}{1-\lambda_2}\right),
\end{equation}

\noindent and, moreover,
as $k \wedge \ell \rightarrow \infty$,

\begin{equation}\label{item5l}
\frac{\mbox{Cov} (S_k^1,S_{\ell}^1)}{(k \wedge \ell)} \rightarrow \sigma^2\left(\frac{1+\lambda_2}{1-\lambda_2}\right).
\end{equation}

\noindent When $a = b$, $\bbe S_k^1 = 0$, and in \eqref{item5k}
the asymptotic variance becomes

\begin{align*}
 \frac{\mbox{Var}S_k^1}{k} 
 &\rightarrow  \frac{4a^2}{(2a)^2}\left(\frac{1+(1-2a)}{1-(1-2a)}\right)\nonumber\\
 &= \frac{1}{a} - 1.
\end{align*}

For $a$ small, we have a ''lazy'' Markov chain,
that is, a Markov chain which tends to remain
in a given state for long periods of time.
In this regime, the random variable $S_k^1$ has long
periods of increase followed by long periods
of decrease.  In this way, linear asymptotics
of the variance with large constants
occur.  If, on the other hand, $a$ is close to $1$,
the Markov chain rapidly shifts back and forth
between $\alpha_1$ and $\alpha_2$, and so
the constant associated with the linearly 
increasing variance of  $S_k^1$ is small.

As in \cite{HL}, Brownian functionals play a central r\^ole in
describing the limiting distribution of $LI_n$.  
By a {\it Brownian motion} on $[0,1]$
we shall mean an a.s.~continuous, centered 
Gaussian process 
having stationary, independent increments,
and which is zero at the origin.
By a {\it standard Brownian motion} 
$B(t), 0 \le t \le 1$,
we shall further require 
that Var$B(t) = t$, $0 \le t \le 1$, {\it i.e.}, 
we endow  $C[0,1]$ with the Wiener measure. 
A {\it standard $m$-dimensional Brownian 
motion} will be defined to be 
a multivariate process consisting of  
$m$ independent Brownian motions.
More generally, an {\it $m$-dimensional Brownian motion} 
shall refer to a linear transformation of
a standard $m$-dimensional Brownian motion.
Throughout the paper, we assume that our underlying
probability space is rich enough so that all 
the Brownian motions and sequences we study
can be defined on it.

To move towards a Brownian functional
expression for the limiting law of $LI_n$, 
define the polygonal function

\begin{equation}\label{item5n}
\hat B_n(t)=\frac{S^1_{[nt]} - [nt]\mu}{\sigma\sqrt{n(1+\lambda_2)/(1-\lambda_2)}} +\frac{(nt-[nt])(Z^1_{[nt]+1} - \mu)}{\sigma\sqrt{n(1+\lambda_2)/(1-\lambda_2)}}, 
\end{equation}

\noindent for $0 \le t \le 1$.  In our finite-state,
irreducible, aperiodic, stationary Markov chain setting, 
we may conclude that $\hat B_n \Rightarrow B$, as desired.
(See, for example, Gordin's martingale approach to dependent
invariance principles \cite{Go}, and the stationary ergodic
invariance principle found in Theorem 19.1 of Billingsley
\cite{Bill}.)

Turning now to $LI_n$, we see that
for the present $2$-letter situation,
\eqref{item5} simply becomes

$$LI_n = \frac{n}{2} - \frac{1}{2}S^1_n + \max_{1 \le k \le n} S^1_k.$$

To find the limiting distribution of $LI_n$ from this
expression, recall that $\pi_1 = b/(a+b)$,
$\pi_2 = a/(a+b)$, $\mu = \pi_1 - \pi_2 = (b-a)/(a+b)$,
$\sigma^2 = 4ab/(a+b)^2$, 
and that $\lambda_2 = 1-a-b$.
Define $\pi_{max} = \max\{\pi_1,\pi_2\}$ and
$\tilde{\sigma}^2 = \sigma^2(1+\lambda_2)/(1-\lambda_2)$.
Rewriting \eqref{item5n} as


$$ \hat B_n(t)=\frac{S^1_{[nt]} - [nt]\mu}{\tilde{\sigma}\sqrt{n}} +\frac{(nt-[nt])(Z^1_{[nt]+1} - \mu)}{\tilde{\sigma}\sqrt{n}}, $$
\
\noindent $LI_n$ becomes

\begin{align}\label{item5q}
LI_n &= \frac{n}{2}  -\frac12 \left(\tilde{\sigma}\sqrt{n} \hat B_n(1) + \mu n\right) + \max_{0\le t \le 1} \left( \tilde{\sigma}\sqrt{n} \hat B_n(t) + \mu nt\right) \nonumber\\
&= n\pi_2  -\frac12 \left(\tilde{\sigma}\sqrt{n} \hat B_n(1)\right) + \max_{0\le t \le 1} \left( \tilde{\sigma}\sqrt{n} \hat B_n(t) + (\pi_1-\pi_2) nt\right) \nonumber\\
&= n\pi_{max}  -\frac12 \left(\tilde{\sigma}\sqrt{n} \hat B_n(1)\right) \nonumber\\
&\qquad + \max_{0\le t \le 1} \left( \tilde{\sigma}\sqrt{n} \hat B_n(t) + (\pi_1-\pi_2) nt - (\pi_{max} - \pi_2)n\right).
\end{align}

\noindent This immediately gives

\begin{align}\label{item5qa}
\frac{LI_n  - \pi_{max}n}{\tilde{\sigma}\sqrt{n}} &= -\frac12 \hat B_n(1) \nonumber\\
&+ \max_{0\le t \le 1} \left(  \hat B_n(t) + \frac{\sqrt{n}}{\tilde{\sigma}}((\pi_1-\pi_2)t - (\pi_{max} - \pi_2)) \right).
\end{align}

Let us examine \eqref{item5qa} on a case-by-case basis.
First, if $\pi_{max} = \pi_1 = \pi_2 = 1/2$, {\it i.e.},
if $a = b$, then
$\sigma = 1$ and $\tilde{\sigma} = (1-a)/a$, and
so \eqref{item5qa} becomes 

\begin{align}\label{item5qb}
\frac{LI_n  - n/2}{\sqrt{(1-a)n/a}} &= -\frac12 \hat B_n(1) + \max_{0\le t \le 1}  \hat B_n(t).
\end{align}

\noindent Then, by the Invariance Principle and the Continuous Mapping Theorem,

\begin{equation}\label{item5qba}
\frac{LI_n  - n/2}{\sqrt{(1-a)n/a}} \Rightarrow -\frac12 B(1) + \max_{0\le t \le 1}  B(t).
\end{equation}

Next, if $\pi_{max} = \pi_2 > \pi_1$, \eqref{item5qa} becomes

\begin{align}\label{item5qd}
\frac{LI_n  - \pi_{max}n}{\tilde{\sigma}\sqrt{n}} &= -\frac12 \hat B_n(1) \nonumber\\
&\qquad + \max_{0\le t \le 1} \left(  \hat B_n(t) - \frac{\sqrt{n}}{\tilde{\sigma}}(\pi_{max} - \pi_1)t \right).
\end{align}

On the other hand, if $\pi_{max} = \pi_1 > \pi_2$, \eqref{item5qa} becomes

\begin{align}\label{item5qc}
\frac{LI_n  - \pi_{max}n}{\tilde{\sigma}\sqrt{n}} 
&= -\frac12 \hat B_n(1) \nonumber\\
&\qquad + \max_{0\le t \le 1} \left(  \hat B_n(t) - \frac{\sqrt{n}}{\tilde{\sigma}}(\pi_{max} - \pi_2)(1-t) \right)\nonumber\\
&= \frac12 \hat B_n(1) \nonumber\\
&\qquad + \max_{0\le t \le 1} \left(  \hat B_n(t) - \hat B_n(1)- \frac{\sqrt{n}}{\tilde{\sigma}}(\pi_{max} - \pi_2)(1-t) \right).
\end{align}

In both \eqref{item5qd} and \eqref{item5qc} we have a term in
our maximal functional which is linear in $t$ or $1-t$,
with a negative slope.  We now show,
in an elementary fashion, that in both cases,
as $n \rightarrow \infty$, the maximal functional goes
to zero in probability.  

Consider first \eqref{item5qd}.  Let
$c_n = \sqrt{n}(\pi_{max} - \pi_1)/\tilde{\sigma} > 0,$
and for any $c > 0$, let 
$M_c = \max_{0\le t \le 1} (B(t) - ct)$,
where $(B(t))$ is a standard Brownian motion.
Now for $n$ large enough,

$$\hat B_n(t) - c t \ge \hat B_n(t) - c_nt$$

\noindent a.s., for all $0 \le t \le 1$.
Then for any $z > 0$, and $n$ large enough, 

\begin{align}\label{item5qe}
\bbp( \max_{0 \le t \le 1} (\hat B_n(t) - c_n t) > z ) 
&\le \bbp( \max_{0 \le t \le 1} (\hat B_n(t) - c t) > z ),
\end{align}

\noindent and so by the Invariance Principle and the
Continuous Mapping Theorem,

\begin{align}\label{item5qf}
\limsup_{n \rightarrow \infty} \bbp( \max_{0 \le t \le 1} (\hat B_n(t) - c_n t) > z ) 
&\le \lim_{n \rightarrow \infty} \bbp( \max_{0 \le t \le 1} (\hat B_n(t) - c t) > z )\nonumber\\
&= \bbp( M_c > z ).
\end{align}

Now, as is well-known,
$\bbp( M_c > z ) \rightarrow 0$ as $c \rightarrow \infty$.
One can confirm this intuitive fact with the following simple argument.
For $z > 0$, $c > 0$, and $0 < \varepsilon < 1$, we have that

{\allowdisplaybreaks
\begin{align}\label{item5qfa}
\bbp ( M_c > z) &\le \bbp (\max_{0 \le t \le \varepsilon} (B(t) - c t) > z ) + \bbp (\max_{\varepsilon < t \le 1} (B(t) - c t) > z )\nonumber\\
&\le  \bbp (\max_{0 \le t \le \varepsilon} B(t) > z ) + \bbp (\max_{\varepsilon < t \le 1} (B(t) - c \varepsilon) > z )\nonumber\\
&\le  \bbp (\max_{0 \le t \le \varepsilon} B(t) > z ) + \bbp (\max_{0 < t \le 1} B(t) >  c \varepsilon + z )\nonumber\\
&= 2\left(1 - \Phi\left(\frac{z}{\sqrt{\varepsilon}}\right) \right) + 2\left(1 - \Phi(c\varepsilon + z) \right).
\end{align}
}

\noindent But, as $c$ and $\varepsilon$ are arbitrary, we can
first take the limsup of \eqref{item5qfa} as $c \rightarrow \infty$,
and then let $\varepsilon \rightarrow 0$, 
proving the claim.

We have thus shown that 

$$\limsup_{n \rightarrow \infty} \bbp( \max_{0 \le t \le 1} (\hat B_n(t) - c_n t) > z ) \le 0,$$

\noindent and since the functional clearly is equal to zero when
$t = 0$, we have

\begin{equation}\label{item5qg}
\max_{0 \le t \le 1} (\hat B_n(t) - c_n t) \stackrel{\bbp}{\rightarrow} 0,
\end{equation}

\noindent as $n \rightarrow \infty$.
Thus, by the Continuous Mapping Theorem, and the Converging Together Lemma,
we obtain the weak convergence result

\begin{equation}\label{item5qh}
\frac{LI_n  - \pi_{max}n}{\tilde{\sigma}\sqrt{n}} \Rightarrow -\frac12 B(1).
\end{equation}

Lastly, consider \eqref{item5qc}.  
Here we need simply note the 
following equality in law, which
follows from the stationary and Markovian
nature of the underlying sequence $(X_n)_{n\ge0}$:

{\allowdisplaybreaks
\begin{align}\label{item5qi}
\hat B_n(t) - \hat B_n(1)&- \frac{\sqrt{n}}{\tilde{\sigma}}(\pi_{max} - \pi_2))(1-t)\nonumber\\
&\stackrel{\cal{L}}{=} -\hat B_n(1-t)- \frac{\sqrt{n}}{\tilde{\sigma}}(\pi_{max} - \pi_2))(1-t),
\end{align}
}

\noindent 
for $t = 0,1/n,\dots,(n-1)/n,1$.
With a change of variables
$(u = 1 - t)$, and noting that $B(t)$
and $-B(t)$ are equal in law,
our previous convergence result
\eqref{item5qg} implies that

\begin{align}\label{item5qj}
\max_{0 \le t \le 1} (\hat B_n(t) - \hat B_n(1)- c_n(1-t))
\stackrel{\cal{L}}{=}  \max_{0 \le u \le 1} (-\hat B_n(u) - c_n u)
\stackrel{\bbp}{\rightarrow} 0,
\end{align}

\noindent as $n \rightarrow \infty$.
Our limiting functional is thus of the form

\begin{equation}\label{item5ql}
\frac{LI_n  - \pi_{max}n}{\tilde{\sigma}\sqrt{n}} \Rightarrow \frac12 B(1).
\end{equation}

\noindent Since $B(1)$ is simply a standard normal random variable,
the different signs in \eqref{item5qh} and \eqref{item5ql}
are inconsequential.

Finally, consider the degenerate cases.
If either $a=0$ or $b = 0$, 
then the sequence $(X_n)_{n \ge 0}$ will be a.s.~constant,
regardless of the starting state, and so
$LI_n \sim n$.
On the other hand, if $a=b=1$, then the sequence
oscillates back and forth between $\alpha_1$
and $\alpha_2$, so that $LI_n \sim n/2$.
Combining these trivial cases with the previous development,
we have proved the following theorem:

\begin{theorem}\label{thm1}
Let $(X_n)_{n \ge 0}$ be a $2$-state Markov chain,
with
$\bbp(X_{n+1} = \alpha_2|$\\
$X_n = \alpha_1) = a$ and
$\bbp(X_{n+1} = \alpha_1 | X_n = \alpha_2) = b$.
Let the law of $X_0$ be the invariant distribution
$(\pi_1,\pi_2) = (b/(a+b),a/(a+b))$, for $0 < a + b \le 2$, and
$(\pi_1,\pi_2) = (1,0)$, for  $a = b = 0$.
Then, for $a=b > 0$,

\begin{equation}\label{item5r}
\frac{LI_n-n/2}{\sqrt {n}} 
\Rightarrow  \sqrt{\frac{1-a}{a}}\left(-\frac12 B(1)+ \max_{0\le t \le 1}B(t)\right),
\end{equation}

\noindent 
where $(B(t))_{t \in [0,1]}$ is a standard Brownian motion,
and for $a \ne b$ or $a=b=0$, 
and $\pi_{max} = \max\{\pi_1,\pi_2\}$,

\begin{equation}\label{item5s}
\frac{LI_n  - \pi_{max}n}{\sqrt{n}} \Rightarrow N(0,\tilde{\sigma}^2/4),
\end{equation}

\noindent where $N(0,\tilde{\sigma}^2/4)$ 
is a centered normal
random variable with 
variance $\tilde{\sigma}^2/4 = ab(2-a-b)/(a+b)^3$, for $a \ne b$,
and $\tilde{\sigma}^2=0$, for $a=b=0$.
(If $a=b=1$, or $\tilde{\sigma}^2 = 0$, then the distributions in \eqref{item5r} and \eqref{item5s},
respectively, are understood to be degenerate at the origin.)
\end{theorem}

To extend this result to the entire Young tableau, let us introduce 
the following notation.  By 

\begin{equation}\label{item5sa}
(Y^{(1)}_n,Y^{(2)}_n,\dots,Y^{(k)}_n) \Rightarrow (Y^{(1)}_{\infty},Y^{(2)}_{\infty},\dots,Y^{(k)}_{\infty})
\end{equation}

\noindent we shall mean the weak convergence of the 
{\it joint} law of the $k$-vector
$(Y^{(1)}_n,Y^{(2)}_n,$\\
$\dots,Y^{(k)}_n)$ to that of
$(Y^{(1)}_{\infty},Y^{(2)}_{\infty},\dots,Y^{(k)}_{\infty})$,
as $n \rightarrow \infty$.
Since $LI_n$ is the length of the top row of the associated
Young tableau, the length of the second row is simply $n - LI_n$. 
Denoting the length of the $i^{th}$ row by $LY^i_n$, \eqref{item5sa},
together with an application of the Cram\'er-Wold Theorem,
recovers the result of Chistyakov and G\"otze \cite{ChG}
as part of the following easy corollary,
which is in fact equivalent to Theorem \ref{thm1}:\\

\begin{corollary} \label{cor1}
For the sequence in Theorem \ref{thm1},
if $a=b>0$, then

\begin{equation}\label{item5t}
\left(\frac{LY^1_n-n/2}{\sqrt {n}}, \frac{LY^2_n-n/2}{\sqrt{n}}\right)  \Rightarrow  Y_{\infty} := (Y^{(1)}_{\infty},Y^{(2)}_{\infty}),
\end{equation}

\noindent where the law of $Y_{\infty}$ is 
supported on the $2^{nd}$ main diagonal of $\bbr^2$,
and with

$$Y^{(1)}_{\infty} \stackrel{\cal{L}}{=}  \sqrt{\frac{1-a}{a}}\left( - \frac12 B(1)+ \max_{0\le t \le 1}B(t)\right).$$

\noindent If $a \ne b$ or $a=b=0$, then setting
$\pi_{min} = \min\{\pi_1,\pi_2\}$, we have

\begin{equation}\label{item5u}
\left(\frac{LY^1_n  - \pi_{max}n}{\sqrt{n}} ,\frac{LY^2_n  - \pi_{min}n}{\sqrt{n}}\right)
\Rightarrow   N((0,0),\tilde{\Sigma}),
\end{equation}

\noindent where $\tilde{\Sigma}$ is the covariance matrix

$$(\tilde{\sigma}^2/4) \begin{pmatrix} 1 &-1  \\ -1 & 1   \end{pmatrix},$$

\noindent where $\tilde{\sigma}^2 = 4ab(2-a-b)/(a+b)^3$, for $a \ne b$,
and $\tilde{\sigma}^2=0$, for $a=b=0$.

\end{corollary}

\begin{Rem}
The joint distributions in \eqref{item5t} and \eqref{item5u}
are of course degenerate, in  that the sum of the two
components is a.s.~identically zero in each case.  In \eqref{item5t},
the density of the first component of $Y_{\infty}$ is
easy to find, and is given by
(e.g., see \cite{HLM})

\begin{equation}\label{item5ua}
f(y) = \frac{16}{\sqrt{2\pi}} \left(\frac{a}{1-a}\right)^{3/2} y^2 e^{-2ay^2/(1-a)}, \qquad y \ge 0.
\end{equation}

\noindent As in Chistyakov and G\"otze \cite{ChG}, 
\eqref{item5t} can then be stated as:
For any bounded, continuous function 
$g:\bbr^2 \rightarrow \bbr$,

\begin{align*}
&\lim_{n \rightarrow \infty} \left( g \left(\frac{LY^1_n-n/2}{\sqrt {(1-a)n/a}}, \frac{LY^2_n-n/2}{\sqrt{(1-a)n/a}}\right) \right)\nonumber\\
&\qquad = 2\sqrt{2\pi}\int_0^{\infty}  g(x,-x)\phi_{GUE,2}(x,-x) dx,
\end{align*}

\noindent where $\phi_{GUE,2}$ is the density of the
eigenvalues of the $2 \times 2$ GUE, and
is given by

$$\phi_{GUE,2}(x_1,x_2) = \frac{1}{\pi}(x_1-x_2)^2 e^{-(x_1^2 + x_2^2)}.$$

To see the GUE connection more explicitly, consider the  
$2 \times 2$ traceless GUE matrix

\begin{equation*}
M_0 =\begin{pmatrix}
X_1 & Y + iZ\\
Y - iZ  &X_2
\end{pmatrix},
\end{equation*}

\noindent where $X_1, X_2, Y$, and $Z$
are centered, normal random variables.
Since $\mbox{Corr }(X_1,X_2) = -1$, 
the largest eigenvalue of $M_0$ is

$$\lambda_{1,0} = \sqrt{X_1^2 + Y^2 + Z^2},$$

\noindent almost surely, 
so that $\lambda_{1,0}^2 \sim \chi_3^2$
if $\mbox{Var } X_1 =  \mbox{Var } Y = \mbox{Var } Z = 1$.
Hence, up to a scaling factor, the density of $\lambda_{1,0}$ is
given by \eqref{item5ua}.  
Next, let us perturb $M_0$ to

$$M = \alpha GI + \beta M_0,$$

\noindent where $\alpha$ and $\beta$ are constants, $G$
is a standard normal random variable independent of $M_0$,
and $I$ is the identity matrix.  
The covariance of the diagonal elements of $M$ is then computed to be
$\rho := \alpha^2 - \beta^2$.  Hence, to obtain a given value
of $\rho$, we may take $\alpha = \sqrt{(1+\rho)/2}$ and
$\beta = \sqrt{(1-\rho)/2}.$  Clearly, the largest eigenvalue
of $M$ can then be expressed as

\begin{equation}\label{item5ub}
\lambda_1 = \sqrt{\frac{1+\rho}{2}}G + \sqrt{\frac{1-\rho}{2}}\lambda_{1,0}.
\end{equation}

\noindent At one extreme, $\rho = -1$, we recover
$\lambda_1 = \lambda_{1,0}$.  At the other extreme, $\rho = 1$,
we obtain $\lambda_1 = Z$.  Midway between these two extremes,
at $\rho=0$, we have a standard GUE matrix, so that

$$\lambda_1 =  \sqrt{\frac{1}{2}}\left(G + \lambda_{1,0}\right).$$
\end{Rem}

\section{Combinatorics Revisited}

The original combinatorial development for the $m$-letter
alphabet resulted in $m-1$ quantities $S_n^r$, $1 \le r \le m-1$.  
In the $2$-letter case we were then able to proceed with a 
probabilistic development which involved a {\it single}
Brownian motion.  Using an even more straightforward
development which involves $m$ quantities instead,
we can obtain more symmetric expressions for $LI_n$.
This is done next, and will prove useful when studying
the shape of the whole Young tableau.\\

Recall that $a^r_k$ counts the number of occurrences
of $\alpha_r$ among $(X_i)_{1\le i \le k}$.
Moving beyond the purely combinatorial setting, 
assume that $(X_k)_{k \ge 0}$ is a doubly-infinite
sequence generated by an irreducible homogeneous Markov chain
having a stationary distribution $(\pi_1,\pi_2,\dots,\pi_m)$.
(For no $k \ge 0$ is the law of $X_k$ necessarily
assumed to be the stationary distribution.)
For each $1 \le r \le m$, 
set $T^r_k = a^r_k - \pi_r k$,
for $k \ge 1$, and $T^r_0 = 0$.
Beginning again with \eqref{item1},
we find that

\begin{align}\label{item6a}
LI_n &=\max_{\stackrel{\scriptstyle 0\le k_1\le\cdots}{\le k_{m-1}\le n}}
      \Bigl[(a^1_{k_1}-a^1_0)+(a^2_{k_2}-a^2_{k_1})+\cdots +
      (a^m_n-a^m_{k_{m-1}})\Bigr]\nonumber\\
     &= \max_{\stackrel{\scriptstyle 0\le k_1\le\cdots}{\le k_{m-1}\le n}}
      \Bigl[( (T^1_{k_1}+ \pi_1 k_1) - (T^1_{k_0}+ \pi_1 k_0) )+ ( (T^2_{k_2}+ \pi_2 k_2) - (T^2_{k_1}+ \pi_2 k_1) )\nonumber\\
     &\qquad +\cdots + ( (T^m_{k_m}+ \pi_m k_m) - (T^m_{k_{m-1}}+ \pi_m k_{m-1}) )\Bigr]\nonumber\\
     &= \max_{\stackrel{\scriptstyle 0\le k_1\le\cdots}{\le k_{m-1}\le n}}
      \Bigl[( T^1_{k_1} - T^1_{k_0} )+ ( T^2_{k_2} - T^2_{k_1} ) +\cdots + ( T^m_{k_m} - T^m_{k_{m-1}} ) \nonumber\\
     & \qquad + \pi_1( k_1 - k_0 ) + \pi_2( k_2 - k_1 ) + \cdots  + \pi_m( k_m - k_{m-1} )    \Bigr].
\end{align}

\noindent Setting $\pi_{max} = \max\{\pi_1, \pi_2, \dots, \pi_m\}$, 
\eqref{item6a} becomes

\begin{align}\label{item6b}
LI_n - \pi_{max}n 
&= \max_{\stackrel{\scriptstyle 0=k_0 \le k_1\le\cdots}{\le k_{m-1}\le k_m=n}}
      \sum_{r=1}^{m} \bigl[ (T^r_{k_r} - T^r_{k_{r-1}})
     + (\pi_r -  \pi_{max}) ( k_r - k_{r-1} ) \bigr].
\end{align}

\noindent For a uniform alphabet, 
$\pi_{max} = \pi_r = 1/m$, for all $r$, and \eqref{item6b} 
simplifies to

\begin{equation}\label{item6c}
LI_n - \frac{n}{m} = \max_{\stackrel{\scriptstyle 0=k_0\le k_1\le\cdots}{\le k_{m-1}\le k_m=n}}
      \sum_{r=1}^{m}  (T^r_{k_r} - T^r_{k_{r-1}}).
\end{equation}

To introduce a random walk formalism into the picture, 
we next set, for $i = 1, \dots ,n$ and $r = 1,2, \dots ,m$,  

\begin{equation}\label{item6d}
W^r_i=	\begin{cases} 	1, &\text{if $X_i=\alpha_r,$}\\
			0, & \text{otherwise.}
	\end{cases}
\end{equation}

\noindent Clearly, $a^r_k = \sum^k_{i=1}W^r_i$,
and so $T^r_k = \sum^k_{i=1}(W^r_i - \pi_r)$,
for $1 \le r \le m$.

To understand the limiting law of \eqref{item6b} or
\eqref{item6c}, we must have a more precise description
of the underlying Markovian structure.  To that end, let 
$p_{r,s} = \bbp(X_{k+1} = \alpha_s|X_k = \alpha_r)$
be the transition probability  from state 
$\alpha_r$ to state $\alpha_s$, and let $P = (p_{r,s})$
be the associated Markov transition matrix.  In this setting,

$$ (p^{n+1}_1, p^{n+1}_2,\dots,p^{n+1}_m) = (p^{n}_1, p^{n}_2,\dots,p^{n}_m)P.$$

\noindent Moreover, as usual, let $p^{(k)}_{r,s}$ denote the $k$-step transition
probability from $\alpha_r$ to $\alpha_s$; its associated transition
matrix is simply $P^k$.

Assume now that the law of $X_0$ is the stationary distribution.
Thus, by construction, $\bbe T^r_k = 0$ for all $1\le r \le m$ and $1 \le k \le n$, and
our primary task is to describe the covariance structure of
these random variables $T_k^r$.

Since $W^r_i$ is, simply, a Bernoulli random variable
with parameter $\pi_r$, $\mbox{Var}W^r_i = \pi_r(1-\pi_r)$.
We then find that, for $k \ge 1$,

\begin{align}\label{item6e}
\mbox{Var}T^r_k 
&= \mbox{Var}\left( \sum^k_{i=1}(W^r_i - \pi_r) \right)\nonumber\\
&= \sum^k_{i=1}\mbox{Var}W^r_i + \sum^{k-1}_{i_1=1} \sum^{k}_{i_2=i_1+1} \mbox{Cov}(W^r_{i_1},W^r_{i_2})\nonumber\\
&\qquad \qquad + \sum^k_{i_1=2} \sum^{i_1 - 1}_{i_2=1} \mbox{Cov}(W^r_{i_1},W^r_{i_2}).
\end{align}

\noindent By stationarity, \eqref{item6e} becomes

{\allowdisplaybreaks
\begin{align}\label{item6f}
\mbox{Var}T^r_k 
&= \sum^k_{i=1}\mbox{Var}W^r_i + \sum^{k-1}_{i_1=1} \sum^{k}_{i_2=i_1+1} \mbox{Cov}(W^r_{0},W^r_{i_2-i_1})\nonumber\\
&\qquad \qquad + \sum^k_{i_1=2} \sum^{i_1 - 1}_{i_2=1} \mbox{Cov}(W^r_{0},W^r_{i_1-i_2})\nonumber\\
&= k\pi_r(1-\pi_r) + \sum^{k-1}_{i_1=1} \sum^{k}_{i_2=i_1+1} (\pi_r p^{(i_2-i_1)}_{r,r} - \pi_r^2)\nonumber\\
&\qquad \qquad + \sum^k_{i_1=2} \sum^{i_1 - 1}_{i_2=1} (\pi_r p^{(i_1-i_2)}_{r,r} - \pi_r^2)\nonumber\\
&= k\pi_r -k^2\pi_r^2 +  \pi_r\sum^{k-1}_{i_1=1} \sum^{k}_{i_2=i_1+1}  e_r P^{i_2-i_1} e_r^T\nonumber\\
&\qquad \qquad + \pi_r\sum^k_{i_1=2} \sum^{i_1 - 1}_{i_2=1}  e_r P^{i_1-i_2} e_r^T,
\end{align}
}

\noindent where $e_r = (0,0,\dots, 0,1,0, \dots 0)$ 
is the $r^{th}$ standard basis vector of $\bbr^m$.  
Setting

\begin{equation}\label{item6fa}
Q_k = \sum^{k-1}_{i_1=1} \sum^{k}_{i_2=i_1+1}  P^{i_2-i_1} = \sum^{k}_{i=1}  (k-i) P^{i},
\end{equation}

\noindent we can rewrite \eqref{item6f} in the simple form

\begin{align}\label{item6g}
\mbox{Var}T^r_k &= k\pi_r - k^2\pi_r^2 + 2\pi_r e_r Q_k e_r^T .
\end{align}

Our description of the covariance structure 
can now be completed using the above results.  
For $r_1 \ne r_2$ and $k \ge 1$,

{\allowdisplaybreaks
\begin{align}\label{item6h}
\mbox{Cov}(T^{r_1}_k, T^{r_2}_k) 
&= \sum^k_{i=1}\mbox{Cov}(W^{r_1}_i,W^{r_2}_i) +\sum^{k-1}_{i_1=1} \sum^{k}_{i_2=i_1+1} \mbox{Cov}(W^{r_1}_{i_1},W^{r_2}_{i_2})\nonumber\\
&\qquad \qquad + \sum^k_{i_1=2} \sum^{i_1 - 1}_{i_2=1} \mbox{Cov}(W^{r_1}_{i_1},W^{r_2}_{i_2})\nonumber\\
&= \sum^k_{i=1}\mbox{Cov}(W^{r_1}_i,W^{r_2}_i)  + \sum^{k-1}_{i_1=1} \sum^{k}_{i_2=i_1+1} \mbox{Cov}(W^{r_1}_{0},W^{r_2}_{i_2-i_1})\nonumber\\
&\qquad \qquad + \sum^k_{i_1=2} \sum^{i_1 - 1}_{i_2=1} \mbox{Cov}(W^{r_2}_{0},W^{r_1}_{i_1-i_2})\nonumber\\
&= -k\pi_{r_1}\pi_{r_2} + \sum^{k-1}_{i_1=1} \sum^{k}_{i_2=i_1+1} (\pi_{r_1} p^{(i_2-i_1)}_{r_1,r_2} - \pi_{r_1}\pi_{r_2})\nonumber\\
&\qquad \qquad + \sum^k_{i_1=2} \sum^{i_1 - 1}_{i_2=1} (\pi_{r_2} p^{(i_1-i_2)}_{r_2,r_1} - \pi_{r_1}\pi_{r_2})\nonumber\\
&=  -k^2\pi_{r_1}\pi_{r_2} +  \pi_{r_1}\sum^{k-1}_{i_1=1} \sum^{k}_{i_2=i_1+1}  e_{r_1} P^{i_2-i_1} e_{r_2}^T\nonumber\\
&\qquad \qquad + \pi_{r_2}\sum^k_{i_1=2} \sum^{i_1 - 1}_{i_2=1}  e_{r_2} P^{i_1-i_2} e_{r_1}^T\nonumber\\
&= -k^2\pi_{r_1}\pi_{r_2} + \pi_{r_1} e_{r_1} Q_k e_{r_2}^T + \pi_{r_2} e_{r_2} Q_k e_{r_1}^T.
\end{align}
}

\begin{Rem} Both \eqref{item6g} and \eqref{item6h} 
appear to be asymptotically quadratic in $k$.  However, 
since $Q_k  = \sum^{k}_{i=i}  (k-i) P^{i}$,
cancellations will show that
when the Markov chain is irreducible and aperiodic, 
the order is, in fact, linear in $k$. 
\end{Rem}

In order to further analyze the asymptotics of $Q_k$, 
we first examine the diagonalization of $P$ for a very
general class of transition matrices.  

\begin{proposition}\label{prop1}
Let $P$ be the $m \times m$ transition matrix of
an irreducible, aperiodic, homogeneous Markov chain
with eigenvalues $1 > |\lambda_2| \ge \cdots \ge |\lambda_m|$,
and let $\Lambda = diag(1,\lambda_2,\dots,\lambda_m)$.
Let $P = S^{-1}\Lambda S$ be the diagonalization of $P$, where the rows
of $S$ consist of the {\it left}-eigenvectors of $P$,
with, moreover, the first row of $S$ being the
stationary distribution $(\pi_1,\pi_2,\dots,\pi_m)$.
Then the first column of $S^{-1}$
is $(1,1,\dots,1)^T$.
\end{proposition}

\noindent \begin{Proof}
Since $P = S^{-1}\Lambda S$, then 
$PS^{-1} = S^{-1}\Lambda$.  Denoting the first column of
$S^{-1}$ by $c_1$, we have
$Pc_1 = c_1$.  But since the rows of $P$ sum to $1$,
we see that $c_1 = (1,1,\dots,1)^T$
satisfies $Pc_1 = c_1$.  Moreover, $c_1$ must be
unique, up to normalization, since the irreducibility
of $P$ implies that $\lambda_1 = 1$ has multiplicity $1$.
Finally, since the inner product of
the first row of $S$ and the first column of $S^{-1}$
is $1$, the correct normalization is indeed
$(1,1,\dots,1)^T$.\CQFD
\end{Proof}

Returning to $Q_k$, as given in \eqref{item6fa},
and using Proposition \ref{prop1}, we then
obtain:

\begin{theorem}\label{thm2}
Let $(X_n)_{n \ge 0}$ be a sequence generated
by an $m$-letter, aperiodic, irreducible, homogeneous
Markov chain with state space
${\cal A}_m = \{\alpha_1 < \cdots < \alpha_m\}$,
transition matrix $P$,
and stationary distribution $(\pi_1,\pi_2,\dots,\pi_m)$.  
Let also the law of $X_0$ be the
stationary distribution.
Moreover,
for $1 \le r \le m$,
let $T^r_k = a^r_k - \pi_r k$,
for $k \ge 1$, and $T^r_0 = 0$, where 
$a^r_k$ is the number of occurrences
of $\alpha_r$ among $(X_i)_{1\le i \le k}$.
Then, for $1 \le r \le m$,

\begin{equation}\label{item6k}
\lim_{k \rightarrow \infty} \frac{\text{Var } T^r_k}{k} = \pi_r\left(1 + 2e_r S^{-1} D S e_r^T  \right),
\end{equation}

\noindent and for $r_1 \ne r_2$,

\begin{equation}\label{item6l}
\lim_{k \rightarrow \infty} \frac{\mbox{Cov}(T^{r_1}_k, T^{r_2}_k)}{k} = \pi_{r_1}e_{r_1} S^{-1} D S e_{r_2}^T  + \pi_{r_2}e_{r_2} S^{-1} D S e_{r_1}^T ,
\end{equation}

\noindent where $P = S^{-1}\Lambda S$ is the standard diagonalization
of $P$ in Proposition \ref{prop1}, and 
$D = \text{diag}(-1/2, \lambda_2/(1-\lambda_2),\dots ,\lambda_m/(1-\lambda_m))$.
That is, the asymptotic covariance matrix of
$(T^1_k,T^2_k,\dots,T^m_k)$ is given by

\begin{equation}\label{item6la}
 \Sigma = \Pi + \Pi(S^{-1} D S) + (S^{-1} D S)^T \Pi,
\end{equation}

\noindent where $ \Pi = diag(\pi_1,\pi_2,\dots,\pi_m)$.
\end{theorem}

\noindent \begin{Proof} Beginning with \eqref{item6fa},
we diagonalize $P$ and find that

{\allowdisplaybreaks
\begin{align}\label{itm6m}
Q_k &= \sum^{k-1}_{i=1} (k-i) (S^{-1}\Lambda S)^{i}\nonumber\\
    &= S^{-1} \left( \sum^{k-1}_{i=1} (k-i) \Lambda ^{i} \right) S\nonumber\\
    &= S^{-1} \text{diag}( h(1), h(\lambda_2), \dots, h(\lambda_m) ) S,
\end{align}
}

\noindent where 
$h(\lambda) := \sum_{k=1}^{n-1} (n-k)\lambda^k$.
Now $h(1) = k(k-1)/2$ is quadratic in $k$, while for $\lambda \ne 1$,

$$h(\lambda) = k \frac{\lambda}{(1-\lambda)} + \frac{\lambda(\lambda^{k} - 1)}{(1-\lambda)^2},$$

\noindent so that $h(\lambda)$ is linear in $k$.  
We thus can write $Q_k$ as the sum
of terms which are, respectively,
quadratic and linear in $k$.  Recalling, moreover,
that the first row of $S$ contains the
stationary distribution, and that the first
column of $S^{-1}$ is $(1,1,\dots,1)^T$, we have

{\allowdisplaybreaks
\begin{align}\label{itm6n}
Q_k &= S^{-1} \text{diag}( h(1), h(\lambda_2), \dots, h(\lambda_m) ) S,\nonumber\\
    &= \frac{k^2}{2} S^{-1} \text{diag}(  1, 0, \dots, 0) S\nonumber\\ 
    &\qquad + k S^{-1} \text{diag}\left( -\frac12, \frac{\lambda_2}{1-\lambda_2}, \dots,  \frac{\lambda_m}{1-\lambda_m}\right) S + o(k)\nonumber\\
    &= \frac{k^2}{2} 
       \begin{pmatrix} 
          \pi_1 & \pi_2 & \cdots & \pi_m \\
	  \pi_1 & \pi_2 & \cdots & \pi_m \\
	  \vdots & \vdots & \cdots & \vdots\\
	  \pi_1 & \pi_2 & \cdots & \pi_m 
       \end{pmatrix} + k S^{-1} D S + o(k).
\end{align}
}

Starting with the variance in \eqref{item6g},
we now find that, for each $1 \le r \le m$,

\begin{align}\label{item6o}
\mbox{Var }T^r_k &= k\pi_r - k^2\pi_r^2 + 2\pi_r e_r Q_k e_r^T\nonumber\\
&= k\pi_r - k^2\pi_r^2 + 2\pi_r\left( \frac{k^2}{2} \pi_r  + k e_r S^{-1} D S e_r^T\right) + o(k)\nonumber\\
&= k\pi_r\left(1 + 2e_r S^{-1} D S e_r^T  \right) + o(k),
\end{align}

\noindent from which the asymptotic result \eqref{item6k}
follows immediately.

An identical development shows that, for $r_1 \ne r_2$,
\eqref{item6h} simplifies to

\begin{align}\label{item6q}
\mbox{Cov}(T^{r_1}_k, T^{r_2}_k) &= -k^2\pi_{r_1}\pi_{r_2} + \pi_{r_1} e_{r_1} Q_k e_{r_2}^T + \pi_{r_2} e_{r_2} Q_k e_{r_1}^T\nonumber\\
&= -k^2\pi_{r_1}\pi_{r_2} + \pi_{r_1}\left(\frac{k^2}{2} \pi_{r_2}  + k e_{r_1} S^{-1} D S e_{r_2}^T\right) \nonumber\\
&\qquad +  \pi_{r_2}\left(\frac{k^2}{2} \pi_{r_1}  + k e_{r_2} S^{-1} D S e_{r_1}^T\right) + o(k)\nonumber\\
&= k\left(\pi_{r_1}e_{r_1} S^{-1} D S e_{r_2}^T  + \pi_{r_2}e_{r_2} S^{-1} D S e_{r_1}^T\right) + o(k),\nonumber\\
\end{align}

\noindent from which the asymptotic result \eqref{item6l}
follows, and so does \eqref{item6la}.\CQFD
\end{Proof}

\begin{Rem}\label{Remiidnonunif}  
To see that \eqref{item6k} and \eqref{item6l} both
recover the covariance results for the iid case investigated by the
authors in \cite{HL}, let $P$ be the transition matrix whose
rows each consist of the stationary distribution
$(\pi_1,\pi_2,\dots,\pi_m)$.
In this case $\lambda_2 = \cdots = \lambda_m = 0$,
and so $D = \text{diag}( -1/2, 0,\dots,0)$. Hence,

\begin{align*}
e_{r_1} S^{-1} D S e_{r_2}^T
&= \left(1,*,\dots,*\right)D\left(\pi_{r_2},*,\dots,*\right)^T\\
&= -\frac{\pi_{r_2}}{2},
\end{align*}

\noindent for all $r_1$ and $r_2$, and so, for each $r$,

$$\lim_{k \rightarrow \infty} \frac{\mbox{Var }T^r_k}{k} = \pi_r\left(1 + 2\left(-\frac{\pi_r}{2}\right) \right) = \pi_r(1-\pi_r),$$

\noindent while, for $r_1 \ne r_2$,

$$\lim_{k \rightarrow \infty} \frac{\mbox{Cov}(T^{r_1}_k, T^{r_2}_k)}{k} 
= \pi_{r_1}\left( -\frac{\pi_{r_2}}{2} \right) + \pi_{r_2}\left( -\frac{\pi_{r_1}}{2} \right)     = - \pi_{r_1}\pi_{r_2}.$$

\noindent  Note that, in the uniform iid case, we have
$\pi_r = 1/m$, for all $1 \le r \le m$.  
Hence, for $r_1 \ne r_2$, the 
asymptotic correlation between $T^{r_1}_k$ and $T^{r_2}_k$
is given by $(-1/(m^2))/((1/m)(1-1/m)) = -1/(m-1)$, so that the covariance matrix
is indeed the permutation-symmetric one obtained in the iid uniform
case in \cite{HL}.  There is, moreover, another Brownian functional
representation for the iid uniform case in \cite{HL} in which the Brownian
motions have a tridiagonal covariance matrix.
\end{Rem}

\section{The Limiting Shape of the Young Tableau}
Thus far, our results have centered on $LI_n$ alone, 
essentially ignoring the larger
question of the structure of the entire Young tableau.
The present section extends
the combinatorial development of the previous section to answer the
question of the limiting shape of the Young tableau.

Our first result in this direction is a purely combinatorial
expression generalizing \eqref{item1}.  
It is standard in the Young tableau literature
to have entries chosen from the set $\{1,2,\dots,m\}$.  
Below, without loss of generality, we allow our entries to be chosen from 
the $m$-letter ordered alphabet 
${\cal A}_m = \{\alpha_1 < \cdots < \alpha_m\}$.

\begin{theorem}\label{thm3}
Let $R^1_n, R^2_n, \dots, R^r_n$ be the lengths of the first $1 \le r \le m$ rows
of the Young tableau generated by the sequence $(X_k)_{1\le k \le n}$
whose elements belong to an ordered alphabet 
${\cal A}_m = \{\alpha_1 < \cdots < \alpha_m\}$.
Then, for each $1 \le r \le m$, the sum of the lengths of the first
$r$ rows of the Young tableau is given by

\begin{equation}\label{item7e}
\sum_{j=1}^r R^j_n =
  \max_{k_{j,\ell} \in J_{r,m}} \sum_{j=1}^r \sum_{\ell=j}^{m-r+j} \left( a^{\ell}_{k_{j,\ell}} - a^{\ell}_{k_{j,\ell-1}} \right),
\end{equation}

\noindent where  
$J_{r,m} = \{(k_{j,\ell}, 1 \le j \le r, 0 \le \ell \le m): k_{j,j-1} = 0, k_{j,m-r+j} = n, 1 \le j \le r;
k_{j,\ell-1} \le k_{j,\ell}, 1 \le j \le r, 1 \le \ell \le m;
k_{j,\ell} \le k_{j-1,\ell}, 2 \le j \le r, 1 \le \ell \le m\}$,
and where $a^{\ell}_k$ is the number of occurrences of
$\alpha_{\ell}$ among $\{X_1, X_2, \dots ,X_k\}$.
\end{theorem}

\noindent \begin{Proof}
Recall that the sum of the lengths of the first $r$ rows of
the Young tableau generated by a sequence $(X_k)_{1\le k \le n}$,
whose letters arise from an $m$-letter alphabet,
has an interpretation in terms of the length of certain increasing sequences.
Indeed, the sum $R^1_n + R^2_n + \cdots + R^r_n$ is equal to the maximum
sum of the lengths of $r$ disjoint, increasing subsequences
of $(X_k)_{1\le k \le n}$, where by {\it disjoint} it is meant
that each element of $(X_k)_{1\le k \le n}$
occurs in at most one of the $r$ subsequences.
(See Lemma 1 of Section 3.2 in \cite{Fu}).
More general results of this sort, involving partial orderings of the
alphabet and associated antichains, are known
as Greene's Theorem \cite{Greene}.  However, such results are not
enough for our purpose.  Below we need a different
way of reconstructing disjoint subsequences.

We begin by examining
an arbitrary collection of $r$ disjoint, increasing subsequences
of $(X_k)_{1\le k \le n}$, and show that we can always map
these $r$ subsequences onto another collection of 
$r$ disjoint, increasing subsequences
whose properties will be amenable to our combinatorial analysis.

Specifically, with the number of rows $r$ fixed, suppose that,
for each $1 \le j \le r$, we have an increasing subsequence 
$(X_{k_\ell^j}^j)_{1\le \ell \le n_j}$ of length $n_j \le n$,
and that the $r$ subsequences are disjoint.

We first construct the new subsequence 
$(\tilde{X}_{\tilde{k}_{\ell}^{1}}^1)_{1 \le \ell \le \tilde{n}_1}$
as follows.  First, place all $\alpha_1$s occurring among the $r$
original subsequences into 
$(\tilde{X}_{\tilde{k}_{\ell}^{1}}^1)_{1 \le \ell \le \tilde{n}_1}$,
if there are any.  If the last $\alpha_1$ occurs at the $n^{th}$
index, then 
$(\tilde{X}_{\tilde{k}_{\ell}^{1}}^1)_{1 \le \ell \le \tilde{n}_1}$
is complete.
Otherwise, place all $\alpha_2$s which occur after the
final $\alpha_1$ into 
$(\tilde{X}_{\tilde{k}_{\ell}^{1}}^1)_{1 \le \ell \le \tilde{n}_1}$,
if there are any.  If the last $\alpha_2$ occurs at the $n^{th}$
index, then 
$(\tilde{X}_{\tilde{k}_{\ell}^{1}}^1)_{1 \le \ell \le \tilde{n}_1}$
is complete.  Otherwise, continue adding, successively,
$\alpha_3, \dots, \alpha_{m-r+1}$ in the same manner.  Thus,
$(\tilde{X}_{\tilde{k}_{\ell}^{1}}^1)_{1 \le \ell \le \tilde{n}_1}$
consists of a weakly increasing sequence 
of length $\tilde{n}_1$
having values in
$\{ \alpha_1, \dots, \alpha_{m-r+1} \}$.

Next, we construct the new subsequence
$(\tilde{X}_{\tilde{k}_{\ell}^{2}}^2)_{1 \le \ell \le \tilde{n}_2}$
similarly.
By considering only those letters among the $r$ original
subsequences which have not already been moved to the first new subsequence,
start with the smallest available letter, $\alpha_2$,
and continue adding, successively, $\alpha_3, \dots, \alpha_{m+r-2}$.
Note that, crucially, all $\alpha_2$s added to 
$(\tilde{X}_{\tilde{k}_{\ell}^{2}}^2)_{1 \le \ell \le \tilde{n}_2}$
occur before the last index at which $\alpha_1$ was added to
the first subsequence.  More generally, each $\alpha_j$, 
$2 \le j \le m-r+2$, added to 
$(\tilde{X}_{\tilde{k}_{\ell}^{2}}^2)_{1 \le \ell \le \tilde{n}_2}$
occurs before the last $\alpha_{j-1}$ was added to the first subsequence.
Thus, $(\tilde{X}_{\tilde{k}_{\ell}^{2}}^2)_{1 \le \ell \le \tilde{n}_2}$
consists of a weakly increasing subsequence of length $\tilde{n}_2$
having values in $\{ \alpha_2, \dots, \alpha_{m-r+2} \}$.

The construction of 
$(\tilde{X}_{\tilde{k}_{\ell}^{j}}^j)_{1 \le \ell \le \tilde{n}_j}$,
for $3 \le j \le r$, continues in the same manner, 
with 
$(\tilde{X}_{\tilde{k}_{\ell}^{j}}^j)_{1 \le \ell \le \tilde{n}_j}$,
constructed from among the entries of the $r$ original subsequences
which were not moved into any of the first $j-1$ new subsequences,
so that
$(\tilde{X}_{\tilde{k}_{\ell}^{j}}^j)_{1 \le \ell \le \tilde{n}_j}$,
consists of a weakly increasing sequence 
of length $\tilde{n}_j$
having values in
$\{ \alpha_j, \dots, \alpha_{m-r+j} \}$.  It is possible that beyond some
$j \ge 2$ the new subsequences may be empty.

We claim that, indeed, the construction of the $r^{th}$ new subsequence
exhausts the set of available entries.  Indeed, without loss of generality,
assume that after we have created the $(r-1)^{th}$ new subsequence,
the set of available entries is non-empty,
and designate the location of the final $\alpha_{\ell}$ to be included
in the $j^{th}$ new subsequence by $k_{j,\ell}$, for
$1 \le j \le r$ and $1 \le \ell \le m$.  (If no $\alpha_{\ell}$ was
available for inclusion, set $k_{j,\ell} = k_{j,\ell-1}$,
where $k_{j,0} = 0,$ for all $1 \le j \le r$.)
Clearly, all $\alpha_1,\alpha_2,\dots,\alpha_{r-1}$ have been
included in the first $r-1$ new subsequences.
If $r=m$, we are done: simply put the remaining $\alpha_r$s
into the $r^{th}$ new subsequence. 
If $r < m$, we may still ask whether there was,
for some $r+1 \le \ell \le m$, 
an $\alpha_{\ell}$
from among the available entries
which occurred before $k_{r,\ell-1}$.
Assume that there is such an $\alpha_{\ell}$.
Now by construction, $k_{j+1,\ell-r+j} \le k_{j,\ell-r+j-1},$
for $1 \le j \le r-1$.  Hence, there exist letters 
$\alpha_{j_1} < \alpha_{j_2} < \dots < \alpha_{j_r} \le \alpha_{\ell-1}$ 
among the original subsequences which occurred after $k_{r,\ell-1}$, and,
moreover, each letter must come from a different subsequence.
But since each original subsequence was increasing, none
of them could have contained an $\alpha_{\ell}$ before
$k_{r,\ell-1}$, and we have a contradiction.

To better understand this construction, consider
the first row of Figure~\ref{figa},
which shows an initial sequence of
length $n=12$, with $m=4$ letters, 
broken into $r=3$ disjoint, increasing subsequences 
of lengths $n_1=3, n_2=4$, and $n_3=3$, and so
with total length $10$.
The final three rows of the diagram show the results
of the operations described above, producing $3$ new
increasing subsequences of length
$\tilde{n}_1=4, \tilde{n}_2=3$, and $\tilde{n}_3=3$.

\begin{figure}
  \begin{center}
    \includegraphics[width=0.8\textwidth]{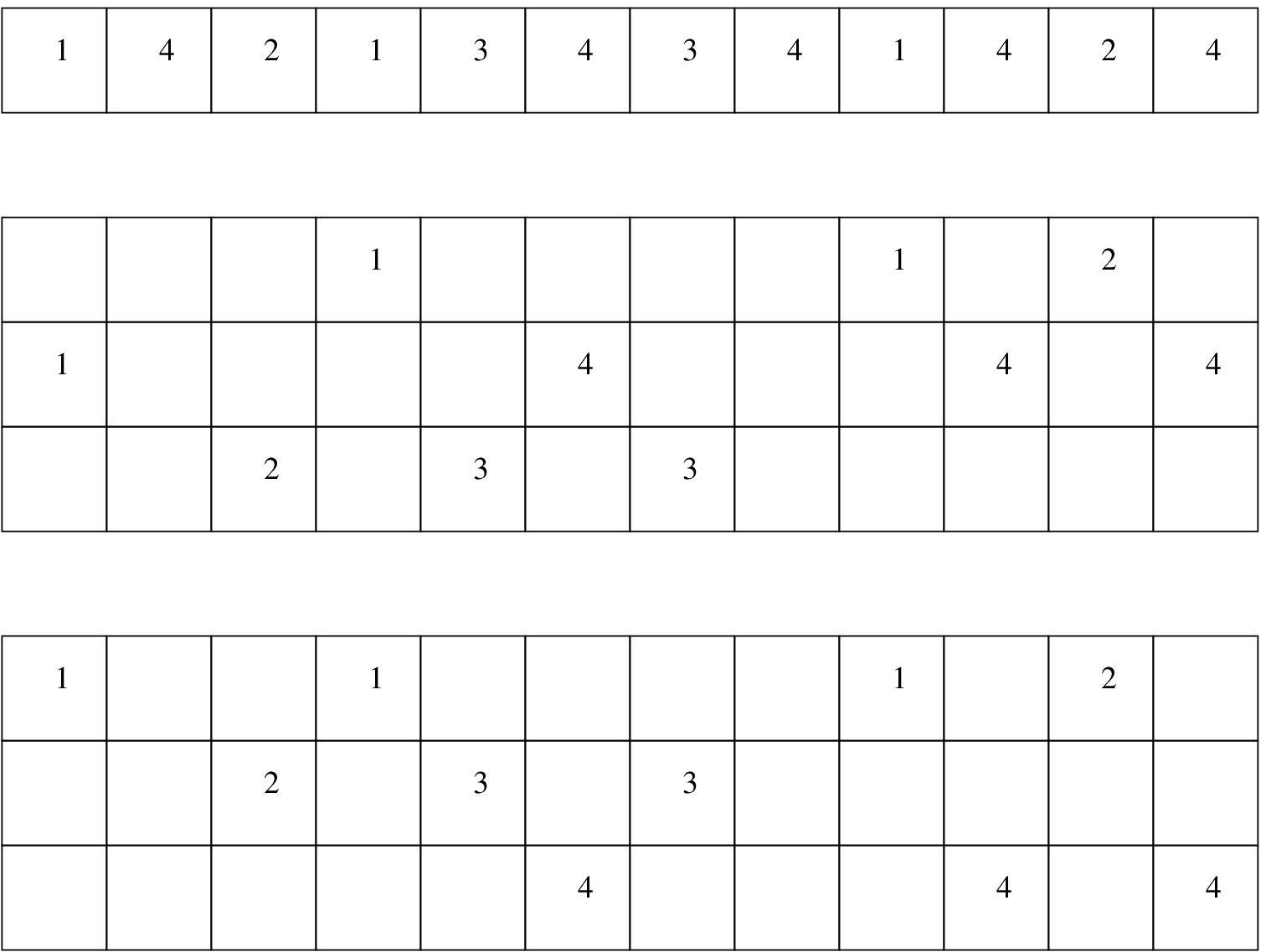}
  \end{center}
  \caption{Transformation of $r=3$ subsequences.}
  \label{figa}
\end{figure}

Hence, if we wish to find $r$ disjoint, increasing subsequences
whose length sum is maximal, it suffices to consider only
those disjoint, increasing subsequences for which
the final occurrence of the letter $\alpha_{\ell}$
in the subsequence $i$ happens after the final occurrence
in the subsequence $j$, whenever $i < j$.
Because such ranges do not overlap, if we wish to count
the number of $\alpha_{\ell}$s in the $j^{th}$ subsequence,
it suffices to simply count the number of $\alpha_s$s 
in $(X_k)_{1\le k \le n}$ over that range.

Indeed, returning to the fundamental combinatorial
objects of our development, the $a^j_k$, we see that
since $a^j_{\ell} - a^j_k$ counts the number
of $\alpha_j$s in the range $\ell + 1, \dots, k$,
we can describe the valid index ranges over which to
search for the maximal sum as
$J_{r,m} = \{(k_{j,\ell}, 1 \le j \le r, 0 \le \ell \le m): k_{j,j-1} = 0, k_{j,m-r+j} = n, 1 \le j \le r;
k_{j,\ell-1} \le k_{j,\ell}, 1 \le j \le r, 1 \le \ell \le m;
k_{j,\ell} \le k_{j-1,\ell}, 2 \le j \le r, 1 \le \ell \le m\}.$
The constraints on the $k_{j,\ell}$ follow simply from the
fact that each subsequence is increasing and that, moreover,
the intervals associated with a given letter do not overlap.
Figure~\ref{figb} indicates the relative positions of
each range, for $r=4$ and $m=7$.

\begin{figure}
  \begin{center}
    \includegraphics[width=0.8\textwidth]{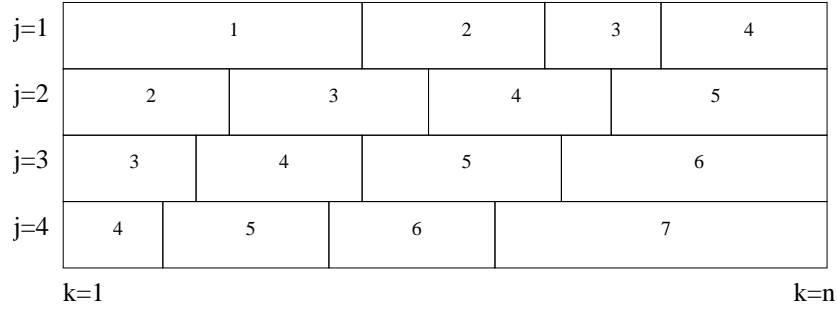}
  \end{center}
  \caption{Schematic diagram of $J_{r,m}$, for $r=4, m=7$.}
  \label{figb}
\end{figure}

Since the first possible letter of each subsequence grows
from $\alpha_1$ to $\alpha_r$, and the last possible letter grows
from $\alpha_{m+r-1}$ to $\alpha_m$, the result is proved.\CQFD
\end{Proof}

We are now ready to apply our asymptotic covariance results 
(Theorem \ref{thm2}), along 
with a Brownian sample-path approximation,
to the combinatorial expression \eqref{item7e}, 
and so obtain a Brownian functional
expression for the limiting shape of the Young tableau for
all irreducible, aperiodic, homogeneous Markov chains.

Indeed, for each $1 \le r \le m$, let

\begin{equation}\label{item7k}
V^r_n := \sum_{j=1}^r R^j_n =
  \max_{k_{j,\ell} \in J_{r,m}} \sum_{j=1}^r \sum_{\ell=j}^{m-r+j} \left( a^{\ell}_{k_{j,\ell}} - a^{\ell}_{k_{j,\ell-1}} \right),
\end{equation}

\noindent where the index set $J_{r,m}$ is defined 
as in Theorem \ref{thm3}.  Define as before 
$T_k^r = \sum^k_{i=1}(W^r_i - \pi_r) = a_k^r - \pi_r k$, 
and so rewrite \eqref{item7k} as

\begin{align}\label{item7l}
V^r_n 
&=  \max_{k_{j,\ell} \in J_{r,m}} \sum_{j=1}^r \sum_{\ell=j}^{m-r+j} \left( \left(T^{\ell}_{k_{j,\ell}}   + \pi_{\ell}k_{j,\ell}  \right)
                                                                       -    \left(T^{\ell}_{k_{j,\ell-1}} + \pi_{\ell}k_{j,\ell-1}\right) \right)\nonumber\\
&=  \max_{k_{j,\ell} \in J_{r,m}}  \sum_{j=1}^r \sum_{\ell=j}^{m-r+j} \left(\left( T^{\ell}_{k_{j,\ell}}  - T^{\ell}_{k_{j,\ell-1}}\right)
                                       +  \pi_{\ell} \left(k_{j,\ell}   -           k_{j,\ell-1}\right)\right).
\end{align}

Next, let $\tau$ be a permutation of the indices $1,2,\dots,m$
such that $\pi_{\tau(1)} \ge \pi_{\tau(2)} \ge \cdots \ge \pi_{\tau(m)} > 0$.
Moreover, we demand that if $\pi_{\tau(i)} = \pi_{\tau(j)}$ for $i < j$,
then $\tau(i) < \tau(j)$.  (The permutation so defined is thus unique.)
Let $\nu_r = \sum_{j=1}^r \pi_{\tau(j)}$ be
the sum of the $r$ largest values among $\pi_1,\pi_2,\dots,\pi_m$.
We obtain, below, the limiting distribution of
$(V^r_n - \nu_r n)/\sqrt{n}$ 
as a Brownian functional.

To introduce Brownian sample-path approximations,
and for each $1 \le r \le m$, we first
define the asymptotic variance of 
$T^r_n$ as in \eqref{item6k}, by

\begin{equation}\label{item7n}
\sigma_r^2 := \lim_{n \rightarrow \infty} \frac{\text{Var } T^r_n}{n} = e_{r} \Sigma e_{r}^T ,
\end{equation}

\noindent and, for $r_1 \ne r_2$,
the asymptotic covariance of
$T^{r_1}_n$ and $T^{r_2}_n$ by

\begin{align}\label{item7o}
\sigma_{r_1,r_2} 
&:= \lim_{n \rightarrow \infty} \frac{\mbox{Cov}(T^{r_1}_n, T^{r_2}_n)}{n} = e_{r_1} \Sigma e_{r_2}^T ,
\end{align}

\noindent where $\Sigma$ is the covariance matrix of 
Theorem \ref{thm2} associated with the transition matrix $P$.
For each $1 \le r \le m$,
we then let

\begin{equation}\label{item7p}
\hat B^r_n(t)=\frac{T^r_{[nt]} + (nt-[nt])(W^r_{[nt]+1} - \pi_r)}{\sigma_r\sqrt{n}},
\end{equation}

\noindent for $0 \le t \le 1$. This rescaling of
$[0,n]$ to $[0,1]$ calls for us to define a new
parameter set over which we will maximize a functional
arising from the expressions in \eqref{item7p}.
Indeed, for any positive integers $s$ and $d$, 
with $s \le d$, define
the set

\begin{align*}
I_{s,d}  = \Bigl\{(t_{j,\ell}, 1 \le j \le s, 0 \le \ell \le d): &t_{j,j-1} = 0, t_{j,d-s+j} = 1, 1 \le j \le s;\nonumber\\
							    &t_{j,\ell-1} \le t_{j,\ell}, 1 \le j \le s, 1 \le \ell \le d;\nonumber\\
                                                            &t_{j,\ell} \le t_{j-1,\ell}, 2 \le j \le s, 1 \le \ell \le d\Bigr\}.
\end{align*}

\begin{figure}
  \begin{center}
    \includegraphics[width=0.8\textwidth]{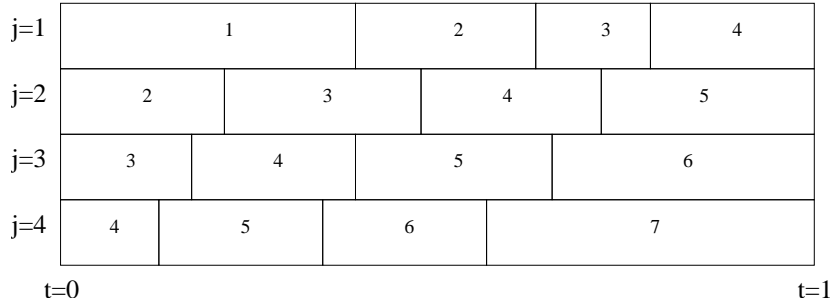}
  \end{center}
  \caption{Schematic diagram of $I_{s,d}$, for $s=4, d=7$.}
  \label{figc}
\end{figure}

\noindent Note that the constraints 
$t_{j,j-1} = 0$ and $t_{j,d-s+j} = 1$, 
for $1 \le j \le s$,
force many of the $t_{j,\ell}$ to be zero or one.
We will denote the $s \times (d+1)$-tuple 
elements of $I_{s,d}$, by $(t_{.,.})$.
Figure~\ref{figc} shows the structure of $I_{s,d}$,
for $s=4$ and $d=7$.  The locations of $t_{j,\ell}$
are indicated by the horizontal lines within the diagram.

With this notation, \eqref{item7l} becomes

\begin{align}\label{item7q}
\frac{V^r_n - \nu_r n}{\sqrt{n}}
&=  \max_{(t_{.,.}) \in I_{r,m}} 
    \Bigl\{      \sum_{j=1}^r \sum_{\ell=j}^{m-r+j} \sigma_{\ell}  \left( \hat B^{\ell}_n(t_{j,\ell})  - \hat B^{\ell}_n(t_{j,\ell-1}) \right) \nonumber\\
&\qquad \qquad \qquad + \sum_{j=1}^r \sum_{\ell=j}^{m-r+j}                \sqrt{n}(\pi_{\ell} - \pi_{\tau(j)})   \left( t_{j,\ell}  - t_{j,\ell-1} \right)      \Bigr\}.  
\end{align}

\noindent

Our analysis of \eqref{item7q} will yield the
following theorem, whose proof we defer to the
conclusion of the section.
This theorem gives, in particular, a full generalization 
of the limiting shape of the Young tableau
in the non-uniform iid case.\\

\begin{theorem}\label{thm4}
Let $(X_n)_{n \ge 0}$ be an irreducible, aperiodic,
homogeneous Markov chain with finite state space
${\cal A}_m = \{\alpha_1 < \cdots < \alpha_m\}$,
transition  matrix $P$, and stationary distribution
$(\pi_1,\pi_2,\dots,\pi_m)$.
Let $\Sigma = (\sigma_{r,s})_{1\le r,s \le m}$ 
be the associated asymptotic covariance matrix,
as given in \eqref{item6la},
and let the law of $X_0$ be given by the
stationary distribution. 
Let $\tau$ be the permutation of $\{1,2,\dots,m\}$
such that $\pi_{\tau(i)} \ge \pi_{\tau(i+1)}$,
and $\tau(i) < \tau(j)$ whenever
$\pi_{\tau(i)} = \pi_{\tau(j)}$ and $i < j$.
For each $1 \le r \le m$,
let $V^r_n$ be the sum of the lengths 
of the first $r$ rows of the associated
Young tableau,
and let $\nu_r = \sum_{j=1}^r \pi_{\tau(j)}$.  
Finally, let $d_r$ be
the multiplicity of $\pi_{\tau(r)}$, and let

\begin{equation*}
m_r = \begin{cases} 0,                                        &\text{if $\pi_{\tau(r)} = \pi_{\tau(1)}$,}\\
		    \max\{i: \pi_{\tau(i)} > \pi_{\tau(r)}\}, & \text{otherwise.}
      \end{cases}
\end{equation*}

\noindent Then, for each $1 \le r \le m$,

\begin{align}\label{item7r}
&\frac{V^r_n -\nu_r n}{\sqrt{n}} \Rightarrow  V^r_{\infty} := \sum_{i=1}^{m_r} \sigma_{\tau(i)} \tilde B^{\tau(i)}(1)  \nonumber\\
&\qquad  + \max_{I_{r-m_r,d_r}} 
   \sum_{j=1}^{r-m_r} \sum_{\ell=j}^{(d_r + m_r - r+j)} \sigma_{\tau(m_r + \ell)} \left(\tilde B^{\tau(m_r + \ell)}(t_{j,\ell})  - \tilde B^{\tau(m_r + \ell)}(t_{j,\ell-1})\right),
\end{align}

\noindent where the first sum on the right-hand side
of \eqref{item7r} is understood to be $0$, if $m_r = 0$.
Above, $\sigma_r^2 = \sigma_{r,r}$,
and $(\tilde B^1(t),\tilde B^2(t),\dots,\tilde B^m(t))$ is an $m$-dimensional Brownian
motion, with covariance matrix 
$\tilde{\Sigma} = ({\tilde\sigma}_{r,s})_{1\le r,s \le m}$ 
given by

\begin{equation}\label{item7s}
({\tilde\sigma}_{r,s})= t({\sigma}_{r,s})/\sigma_r \sigma_s,
\end{equation}

\noindent for $1\le r,s \le m$.
Moreover, for any $1 \le k \le m$,

\begin{equation}\label{item7saa}
\left(\frac{V^1_n -\nu_1 n}{\sqrt{n}},\frac{V^2_n -\nu_2 n}{\sqrt{n}},\dots,\frac{V^k_n -\nu_k n}{\sqrt{n}}\right) 
\Rightarrow \left(V^1_{\infty},V^2_{\infty},\dots,V^k_{\infty} \right).
\end{equation}



\end{theorem}

\begin{Rem}
The critical indices $d_r$ and $m_r$ in Theorem \ref{thm4} are chosen so
that 
$$\pi_{\tau(m_r)} > \pi_{\tau(m_r + 1)} = \pi_{\tau(r)} = \cdots =
	  \pi_{\tau(m_r + d_r)} > \pi_{\tau(m_r + d_r + 1)}.$$
	  
\noindent Thus, the functional in \eqref{item7r} consists of a sum of
$m_r$ Gaussian random variables and a maximal functional involving only $d_r$ of
the $m$ one-dimensional Brownian motions.

\end{Rem}

\begin{Rem}
Another, more natural, way of describing the covariance structure of the
$m$-dimensional Brownian motion in Theorem \ref{thm4} is to
note that $(\sigma_1 B^1(t),\sigma_2 B^2(t)$,
$\dots,\sigma_m B^m(t))$
has covariance matrix $t\Sigma$.
\end{Rem}

Let us now examine the case $r=1$. Here, as
previously noted, $V^1_n = LI_n$.  
Since $m_1 = 0$, 
\eqref{item7r} becomes

\begin{align}\label{item7ra}
\frac{LI_n - \pi_{max} n}{\sqrt{n}} \Rightarrow  
\max_{(t_{.,.}) \in I_{1,d_1}} 
    \sum_{\ell=1}^{d_1 } \sigma_{\tau(\ell)} \left(\tilde B^{\tau(\ell)}(t_{1,\ell})  - \tilde B^{\tau(\ell)}(t_{1,\ell-1})\right),
\end{align}

\noindent where we have written 
$\pi_{max}$ for $\pi_{\tau(1)}.$ 
The functional in \eqref{item7ra} is similar to the one
obtained in the iid case in \cite{HL},  
the essential difference being, 
not in the form of the Brownian functional,
but rather in the covariance structure of the
Brownian motions.

To see precisely where this difference comes into play, 
note that if the transition matrix $P$ is cyclic, 
then the covariance matrix of the Brownian motion 
is also cyclic. Consider then the $3$-letter aperiodic, homogeneous, 
doubly-stochastic Markov case. Since the
Brownian covariance matrix is symmetric, and, moreover,
degenerate, an additional cyclicity constraint forces 
it to have the permutation-symmetric structure seen in the 
iid uniform case.  In particular, $LI_n$ will have, up
to a scaling factor, the same limiting distribution as
in the iid uniform case:

\begin{align}\label{item7rab}
\frac{LI_n - n/3}{\sqrt{n}} \Rightarrow  
\sigma  \max_{(t_{.,.}) \in I_{1,3}} 
    \sum_{\ell=1}^{3}  \left( \tilde B^{\ell}(t_{1,\ell})  - \tilde B^{\ell}(t_{1,\ell-1})\right),
\end{align}

\noindent where $\sigma  = \sigma_{\ell},$
for all $1 \le \ell \le 3$, and with the
Brownian covariance matrix given by

\begin{equation*}
{\tilde \Sigma}= t\begin{pmatrix}
     1      &-1/2   &-1/2\\
    -1/2   &1       &-1/2\\
    -1/2   &-1/2    &1
  \end{pmatrix},
\end{equation*}

\noindent and where we have used the fact
that $\tau(\ell) = \ell$, for all $1 \le \ell \le 3$.

However, when $m\ge 4$, the cyclicity constraint
does not force the Brownian covariance matrix
to have the permutation-symmetric structure, 
as the following example shows for $m=4$.

\begin{Example}\label{Ex1}
Consider the following doubly-stochastic, aperiodic,
cyclic transition matrix:

\begin{equation}\label{item7t}
P =   \begin{pmatrix}
    0.4   &0.3    &0.2   &0.1\\
    0.1   &0.4    &0.3   &0.2\\
    0.2   &0.1    &0.4   &0.3\\
    0.3   &0.2    &0.1   &0.4
  \end{pmatrix}.
\end{equation}

\noindent While the doubly-stochastic nature of $P$
ensures that the stationary distribution is
uniform, the covariance matrix of the limiting Brownian motion,
at three-decimal accuracy, is computed to be

\begin{equation}\label{item7ta}
{\tilde \Sigma} = t\begin{pmatrix}
    1.000    &-0.357   &-0.287   &-0.357\\
    -0.357   &1.000    &-0.357   &-0.287\\
    -0.287   &-0.357   &1.000    &-0.357\\
    -0.357   &-0.287   &-0.357   &1.000
  \end{pmatrix},
\end{equation}

\noindent  and $\sigma_r^2 = \sigma^2 := 0.263$,
for each $1 \le r \le 4$.  Thus, the limiting
distribution of $LI_n$ is given by

\begin{align}\label{item7tb}
\frac{LI_n -n/4}{\sqrt{n}} \Rightarrow  \sigma \max_{(t_{.,.}) \in I_{1,4}} 
      \sum_{\ell=j}^{4}  \left(\tilde B^{\ell}(t_{1,\ell})  - \tilde B^{\ell}(t_{1,\ell-1})\right),
\end{align}

\noindent for $1 \le r \le 4$.  However, while the form 
of the functional is the same as in
the iid uniform case (up to the constant),
the covariance structure of the Brownian motion in \eqref{item7ta}
differs from that of the uniform iid case, {\it i.e.}, from

\begin{equation}\label{item7tc}
t\begin{pmatrix}
    1      &-1/3   &-1/3   &-1/3\\
    -1/3   &1      &-1/3   &-1/3\\
    -1/3   &-1/3   &1      &-1/3\\
    -1/3   &-1/3   &-1/3   &1
  \end{pmatrix},
\end{equation}

\noindent and so the limiting distribution in \eqref{item7tb}
is not that of the uniform iid case.      

\end{Example}

We thus see that Kuperberg's conjecture
regarding the shape of the Young tableau 
for random sequences generated by
aperiodic, homogeneous, and cyclic matrices
\cite{Ku} is not true for general $m$-alphabets. 
By simply extending the first-row analysis above to the
second and third rows, we see that it is true for $m=3$.
However, as could have been anticipated by \eqref{item7rab},
it fails for $m \ge 4$, as the previous example showed.
Furthermore, in the next section we shall see that
for the cyclic case the structure of $\Sigma$ can
be described in an elegant manner which more
clearly delineates when we obtain the uniform
iid limiting law.

In the more general doubly stochastic case, 
we have the following corollary:

\begin{corollary}\label{cor1a}
Let the transition matrix $P$ of Theorem \ref{thm4} be doubly stochastic.
Then, for every $1 \le r \le m$, $m_r = 0, d_r = m$, and 

\begin{align}\label{item7sa}
\frac{V^r_n -rn/m}{\sqrt{n}} \Rightarrow  \max_{(t_{.,.}) \in I_{r,m}} 
      \sum_{j=1}^{r} \sum_{\ell=j}^{m-r+j} \sigma_{\ell} \left(\tilde B^{\ell}(t_{j,\ell})  - \tilde B^{\ell}(t_{j,\ell-1})\right).
\end{align}

If, moreover, the matrix $P$ has all entries of $1/m$
({\it i.e.}, in the iid uniform alphabet case), then

\begin{align}\label{item7sb}
\frac{V^r_n -rn/m}{\sqrt{n}} \Rightarrow    \frac{\sqrt{m-1}}{m}  \max_{(t_{.,.}) \in I_{r,m}} 
    \sum_{j=1}^{r} \sum_{\ell=j}^{m-r+j}  \left(\tilde B^{\ell}(t_{j,\ell})  - \tilde B^{\ell}(t_{j,\ell-1})\right),
\end{align}

\noindent and the covariance matrix in \eqref{item7s} 
has all its off-diagonals equal to $-1/(m-1)$.

\end{corollary}

\noindent \begin{Proof}
For each $1 \le r \le m$, $\pi_r = 1/m$, and so $\nu_r = r/m$, $m_r = 0$, 
and the multiplicity $d_r = m$.  
Moreover, the 
permutation $\tau$ is simply the identity permutation.
This proves \eqref{item7sa}.  If, moreover, all the transition
probabilities are $1/m$, then
the multinomial nature of the underlying combinatorial
quantities $a_k^r$ tells us that
$\sigma_r^2 = (1/m)(1- 1/m)$, for each $1 \le r \le m$, and that
$\rho_{r_1,r_2} = -1/(m-1)$, for each $r_1 \ne r_2$,
thus  proving \eqref{item7sb}.\CQFD
\end{Proof}

To see that the functional in \eqref{item7sa}
is generally different from the uniform iid case,
even for $m=3$, consider the following 
non-cyclic example:

\begin{Example}\label{Ex2}
Let a doubly-stochastic (but non-cyclic), 
aperiodic Markov chain have 
transition matrix

\begin{equation}\label{Ex2a}
P =   \begin{pmatrix}
    0.4   &0.6    &0.0\\
    0.6   &0.0    &0.4  \\
    0.0   &0.4    &0.6  \\
  \end{pmatrix}.
\end{equation}

\noindent As in Example \ref{Ex1},
the doubly-stochastic nature of $P$
ensures that the stationary distribution is
uniform.  In the present example, 
the asymptotic covariance matrix, 
at three-decimal accuracy, is computed to be

\begin{equation}\label{Ex2b}
\begin{pmatrix}  
    0.459    &0.049   &-0.506\\
    0.049    &0.086   &-0.136\\
   -0.506   &-0.136    &0.642
  \end{pmatrix}.
\end{equation}

\noindent Note that, even though we have a uniform
stationary distribution, the asymptotic
variances ({\it i.e.}, the diagonals of \eqref{Ex2b}) 
have dramatically different values.
Moreover, according to
Remark \ref{Remiidnonunif},
in the uniform iid case, the only possibility for
the Brownian covariance matrix is that
the off-diagonals have value $-1/2$.
However, the Brownian motion 
covariance matrix obtained from 
\eqref{Ex2b} is

\begin{equation}\label{Ex2c}
t\begin{pmatrix}   
   1.000     &0.246   &-0.935\\
    0.246    &1.000   &-0.577\\
   -0.935   &-0.577    &1.000
  \end{pmatrix}.
\end{equation}

\noindent Not only are the off-diagonals
different from $-1/2$, but in some cases
are even positive. In
short, the functional in \eqref{item7sa}
has a distribution which differs from
{\it any} iid case (even non-uniform).

\end{Example}

\begin{Rem}
Generalizing a result of Baryshnikov \cite{Ba} and
of Gravner, Tracy, and Widom \cite{GTW}
on the representation of the maximal eigenvalue
of an $m \times m$ element of the GUE, Doumerc \cite{Do} 
found a Brownian functional expression for all the
eigenvalues of an $m \times m$ element of the GUE.
Our expression in \eqref{item7sb}
is similar, with the exception that
our $m$-dimensional Brownian motion is constrained by a
zero-sum condition, and, moreover, has a 
different covariance structure.
(We note, moreover, that the parameters over which his
Brownian functional is maximized in \cite{Do}  
might be intended to range over a slightly larger set
which corresponds to our $I_{r,m}$.)  
Using a path-transformation technique relating the
joint distribution of a certain transformation of
$n$ continuous processes to the joint distribution
of the processes conditioned never to leave the Weyl
chamber, O'Connell and Yor \cite{OY2} employed
queuing-theoretic arguments to obtain 
Brownian functional representations for the
entire spectrum of the $m \times m$ element of the GUE.
In a study of much more general transformations
of this type, Bougerol and Jeulin \cite{BJ} were able
to obtain this result as a special case.
\end{Rem}

If $d_r=1$, {\it i.e.,} if the
$r^{th}$ most probable state is unique, then
the following result can be
viewed as lying at the other extreme 
from Corollary \ref{cor1a}:

\begin{corollary}\label{cor1b}
Let $1 \le r \le m$, and let $d_r = 1 $ in Theorem \ref{thm4}.
Then

\begin{align}\label{item7sc}
\frac{V^r_n -\nu_r n}{\sqrt{n}} \Rightarrow  \sum_{i=1}^{r} \sigma_{\tau(i)} \tilde B^{\tau(i)}(1).
\end{align}

\end{corollary}

\noindent \begin{Proof}
If $d_r = 1$, then $m_r = r - 1$, and so the maximal term
of \eqref{item7r} contains only one summand, namely
$\sigma_{\tau(m_r + 1)} \tilde B^{\tau(m_r + 1)}(1) = \sigma_{\tau(r)} \tilde B^{\tau(r)}(1)$.
Including this term in the first summation term of \eqref{item7r}
proves \eqref{item7sc}. \CQFD
\end{Proof}

\begin{Rem}
The maximal term of the functional in \eqref{item7r}
is that of the doubly-stochastic,  $d_r$-letter case.
Indeed, the maximal term involves precisely $d_r$
Brownian motions over the $r-m_r$ rows.  Such a functional
would arise in a doubly-stochastic $d_r$-letter situation
with a covariance matrix consisting of the sub-matrix of
the original $\Sigma$ corresponding to the $d_r$ 
Brownian motions, as in Corollary \ref{cor1a}.
The Gaussian term corresponds to the functional 
of Corollary \ref{cor1b}. That is, in some sense,
the limiting law of \eqref{item7r} interpolates between
these two extreme cases.
\end{Rem}


\noindent \begin{Proof} {\bf (Theorem \ref{thm4})}
Since the $r=m$  case is trivial 
($V_n^m$ is then identically equal to $n$),
assume that $r < m$.
Recall the approximating functional \eqref{item7q}:

\begin{align}\label{item7u}
\frac{V^r_n - \nu_r n}{\sqrt{n}}
&=  \max_{I_{r,m}} 
    \biggl\{      \sum_{j=1}^r \sum_{\ell=j}^{m-r+j} \sigma_{\ell} \left( \hat B^{\ell}_n(t_{j,\ell})  - \hat B^{\ell}_n(t_{j,\ell-1})     \right)\nonumber\\
&\qquad \qquad + \sum_{j=1}^r \sum_{\ell=j}^{m-r+j}                \sqrt{n}(\pi_{\ell} - \pi_{\tau(j)})   \left( t_{j,\ell}  - t_{j,\ell-1} \right)      \biggr\}.
\end{align}

\noindent Introducing the notation $\Delta t_{j,\ell} := [t_{j,\ell-1}, t_{j,\ell-1}]$
and $M^{\ell}_n(\Delta t_{j,\ell}) := M^{\ell}_n(t_{j,\ell})  - M^{\ell}_n(t_{j,\ell-1})$,
for any $m$-dimensional process $M(t) = (M^1(t),M^2(t),\dots,M^m(t))$, $t \in [0,1]$,
we can rewrite \eqref{item7u} more compactly as

\begin{align}\label{item7ua}
\frac{V^r_n - \nu_r n}{\sqrt{n}}
=  \max_{I_{r,m}} 
    \biggl\{ \sum_{j=1}^r \sum_{\ell=j}^{m-r+j} \sigma_{\ell} \hat B^{\ell}_n(\Delta t_{j,\ell})
 -  \sqrt{n} \sum_{j=1}^r \sum_{\ell=j}^{m-r+j} (\pi_{\tau(j)} - \pi_{\ell})  |\Delta t_{j,\ell}| \biggr\}.  
\end{align}

The main idea of the proof to follow will be to show that the second summation of \eqref{item7ua} 
can, in effect, be eliminated by choosing the 
$(\Delta t_{j,\ell})$ in an appropriate
manner.  Now some of the coefficients $(\pi_{\tau(j) - \pi_{\ell}})$ are zero; 
such terms do not cause any problems.
Intuitively, however, the remaining terms should have
$|\Delta t_{j,\ell}| = 0$.  Defining the restricted set of parameters
$I^{*}_{r,m} = \{(t_{j,\ell}) \in I_{r,m}:\sum_{j=1}^{r} \sum_{\ell=j}^{m-r+j} (\pi_{\ell} - \pi_{\tau(j)})|\Delta t_{j,\ell}|=0, 1 \le \ell \le m  \}$,
we see that, provided $I_{r,m}^{*} \ne \emptyset$, 

\begin{align}\label{this1}
&\max_{I_{r,m}} \sum_{j=1}^{r} \sum_{\ell=j}^{m-r+j} \left(\sigma_{\ell} \hat{B}_n^{\ell}(\Delta t_{j,\ell}) 
- \sqrt{n} \left( \pi_{\tau(j)} - \pi_{\ell}\right) |\Delta t_{j,\ell}|\right) \nonumber\\
&\qquad \qquad \qquad \ge \max_{I_{r,m}^{*}} \sum_{j=1}^{r} \sum_{\ell=j}^{m-r+j} \sigma_{\ell} \hat{B}_n^{\ell}(\Delta t_{j,\ell}).
\end{align}

Moreover, by the
Invariance Principle and the Continuous Mapping Theorem,

\begin{align}\label{this2}
&\max_{I_{r,m}^{*}} \sum_{j=1}^{r} \sum_{\ell=j}^{m-r+j} \sigma_{\ell} \hat{B}_n^{\ell}(\Delta t_{j,\ell})
    \Rightarrow \max_{I_{r,m}^{*}} \sum_{j=1}^{r} \sum_{\ell=j}^{m-r+j} \sigma_{\ell} \tilde B^{\ell}(\Delta t_{j,\ell}).
\end{align}

We claim that, indeed, $I_{r,m}^{*} \ne \emptyset$, 
and that, moreover,

\begin{align}\label{this3}
&\max_{I_{r,m}} \sum_{j=1}^{r} \sum_{\ell=j}^{m-r+j} \left(\sigma_{\ell} \hat{B}_n^{\ell}(\Delta t_{j,\ell}) 
- \sqrt{n} \left( \pi_{\tau(j)} - \pi_{\ell}\right) |\Delta t_{j,\ell}| \right)\nonumber\\
&\qquad \qquad \qquad \Rightarrow \max_{I_{r,m}^{*}} \sum_{j=1}^{r} \sum_{\ell=j}^{m-r+j} \sigma_{\ell} \tilde B^{\ell}(\Delta t_{j,\ell}).
\end{align}

\noindent We will prove that $I_{r,m}^{*} \ne \emptyset$
by creating a bijection between 
$I_{r,m}^{*}$ and $I_{r-m_r,d_r}$.
To this end, for $1 \le i \le m_r$, let 
$\tilde{I}_{\tau(i),i} = [u_{\tau(i),i-1},u_{\tau(i),i}] = [0,1]$.
Next, choose any $(u_{.,.}) \in I_{r-m_r,d_r}$, 
and define further intervals 
$\tilde{I}_{\tau(m_r+j),\ell} = \Delta u_{j,\ell}$,
for $1 \le j \le r-m_r$ and $1 \le \ell \le d_r$.

We now create a partition of these intervals in a
manner which relies on the ideas used in the 
proof of Theorem \ref{thm3}.
Consider the {\it set} of points 
$\{ u_{j,\ell}\}_{(1 \le j \le r-m_r, 1 \le \ell \le d_r)}$,
and order them as
$s_0 := 0 < s_1 < \cdots < s_{\kappa-1} < s_{\kappa}:= 1$,
for some integer $\kappa$,
and let $\Delta s_q = [s_{q-1},s_q]$, 
for all $1 \le q \le \kappa$.

Trivially, for each $1 \le q \le \kappa$, 
and for each $1 \le i \le m_r$, 
$\Delta s_q \subset \tilde{I}_{\tau(i),i}$.
Moreover, for each $1 \le j \le r-m_r$, there exists a unique
$\ell(j,q)$ such that 
$\Delta s_q \subset \tilde{I}_{\tau(m_r+j),\ell(j,q)}$.
For each $q$, consider the set of indices
$A_q := \{\tau(1), \dots, \tau(m_r)\} \cup \{ \tau(m_r + \ell(1,q)), \dots, \tau(m_r + \ell(r-m_r,q)) \}$,
and order these $r$ elements of $A_q$ as
$1 \le \tilde{\ell}(1,q) < \cdots < \tilde{\ell}(r,q) \le m$.

Using these partitions, we examine,
with foresight, the following functional
of a general $m$-dimensional process $(M(t))_{t\ge 0}$:

{\allowdisplaybreaks
\begin{align}
& \sum_{i=1}^{m_r} M^{\tau(i)}(1) 
  + \sum_{j=1}^{(r-m_r)} \sum_{\ell=j}^{(r-m_r+d_r-1)} M^{\tau(m_r + \ell)}(\Delta u_{j,\ell})\label{this4a}\\
&\qquad= \sum_{i=1}^{m_r} \left( \sum_{q=1}^{\kappa}  M^{\tau(i)}(\Delta s_q)  \right)\nonumber\\
&\qquad\qquad + \sum_{j=1}^{(r-m_r)} \sum_{\ell=j}^{(r-m_r+d_r-1)} \left( \sum_{q:\Delta s_q \subset \tilde{I}_{\tau(m_r+j),\ell} } 
          M^{\tau(m_r + \ell)}(\Delta s_q)\right)\nonumber\\
&\qquad= \sum_{q=1}^{\kappa} \left( \sum_{i=1}^{m_r}  M^{\tau(i)}(\Delta s_q)  
             +  \sum_{j=1}^{(r-m_r)}  M^{\tau(m_r + \ell(j,q))}(\Delta s_q)   \right)\nonumber\\
&\qquad= \sum_{q=1}^{\kappa}   \sum_{j=1}^{r}        M^{\tilde{\ell}(j,q)}(\Delta s_q) 
 = \sum_{j=1}^{r}        \sum_{q=1}^{\kappa}   M^{\tilde{\ell}(j,q)}(\Delta s_q)\nonumber\\
&\qquad= \sum_{j=1}^{r}        \sum_{\ell=1}^{r}     M^{\tilde{\ell}(j,q)}(\Delta t_{j,\ell})\label{this4},
\end{align}
}

\noindent where, for each $1 \le j \le r$,
and for each $1 \le \ell \le m$,
$t_{j,\ell} := \max\{s_q: \ell \ge \tilde{\ell}(j,q) \}$.
(That is, for each $j$, we collapse together intervals
$\Delta s_q$ corresponding to the same component
$M^{\ell}$.)  Now, since our functional in 
\eqref{this4} has non-trivial summands only for
$\ell$ such that $\pi_{\tau(\ell)} \ge \pi_{\tau(r)}$, 
we have shown that
$(t_{.,.}) \in I_{r,m}^{*}$.

The following example illustrates this argument.  Suppose
we have an alphabet of size
$m=8$, with 
$$(\pi_1, \pi_2, \dots, \pi_8) = (0.07, 0.1, 0.2, 0.06, 0.2, 0.06, 0.1, 0.2).$$
Then,
$$\pi_{\tau(1)} = \pi_{\tau(2)} =\pi_{\tau(3)} = 0.2, \quad m_1 = m_2 = m_3 = 0, \quad d_1=d_2=d_3=3,$$
$$\pi_{\tau(4)} = \pi_{\tau(5)} = 0.1, \quad m_4 = m_5 = 3, \quad d_4=d_5=2,$$
$$\pi_{\tau(6)} = 0.07, \quad m_6 = 5, \quad d_6=1,$$
$$\pi_{\tau(7)} = \pi_{\tau(8)} = 0.06, \quad m_7 = m_8 = 6, \quad d_7=d_8=2.$$

In particular, note that the two largest, distinct
probability values are $0.2$ and $0.1$,
of multiplicities $3$ and $2$, respectively.
Next, consider the case $r=4$.  We now show how 
$I_{r-m_r,d_r} = I_{4-3,2} = I_{1,2}$ corresponds
to an element of $I_{r,m}^{*} = I_{4,8}^{*}$.
Figure~\ref{figd} shows a typical element of 
the unconstrained index set $I_{4,8}$.

\begin{figure}
  \begin{center}
    \includegraphics[width=0.8\textwidth]{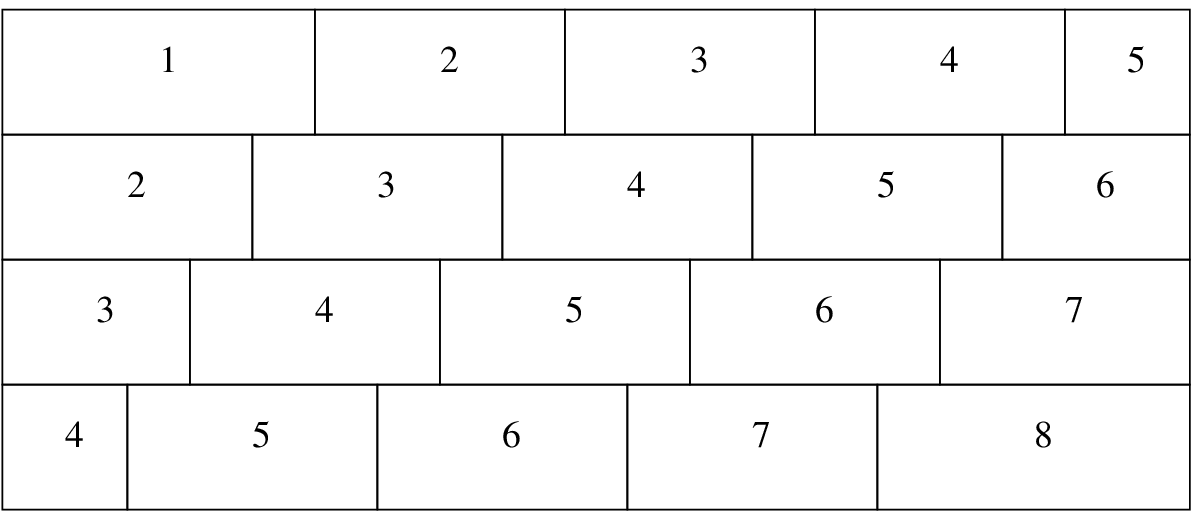}
  \end{center}
  \caption{A typical element of $I_{4,8}$.}
  \label{figd}
\end{figure}

Now $\tau(1) = 3, \tau(2)=5, \tau(3)=8, \tau(4) = 2,$
and $\tau(5) = 7$.  Our construction begins with
the amalgamation of $m_r=m_4=3$ rows, corresponding
to the three indices for which $\pi_i$ is strictly
less than $\pi_{\tau(r)} = \pi_{\tau(4)} = 0.1$, with
$I_{1,2}$.  This is shown in Figure~\ref{fige}.

\begin{figure}
  \begin{center}
    \includegraphics[width=0.7\textwidth]{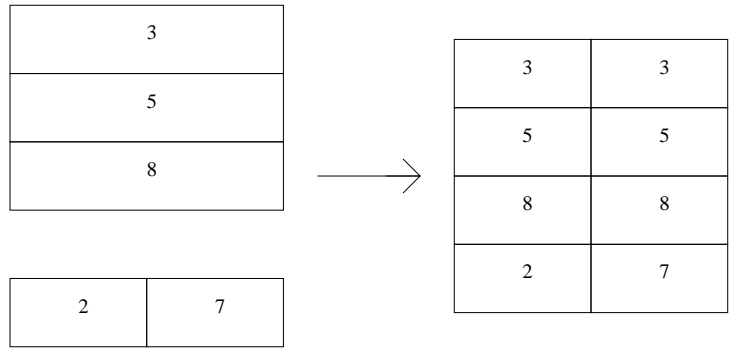}
  \end{center}
  \caption{Amalgamating $3$ rows with $I_{1,2}$.}
  \label{fige}
\end{figure}

Finally, we simply reorder each vertical column in the
original order of the indices,
as shown in Figure~\ref{figf}.
We see that, first of all,
we have constructed an element of $I_{4,8}$.  Moreover,
since we have three rows whose indices are associated
with the maximum value, and a remaining row whose indices
are associated with $\pi_{\tau(4)}$, we indeed have
an element of $I_{4,8}^{*}$.  Note that the $4 \times 4 = 16$
free indices in $I_{4,8}$ (corresponding to the locations of
the $16$ vertical bars in Figure~\ref{figd}) have been
reduced to a {\it single} one in $I_{4,8}^{*}$.

\begin{figure}
  \begin{center}
    \includegraphics[width=0.7\textwidth]{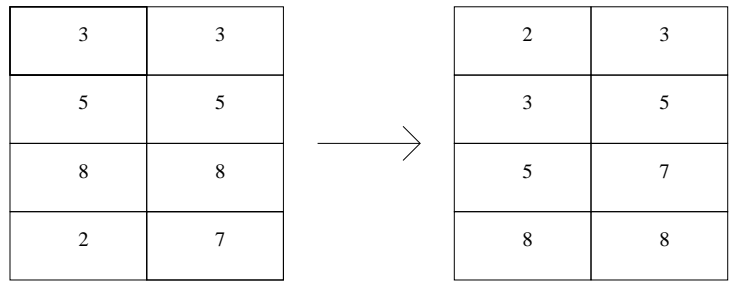}
  \end{center}
  \caption{Reordering vertically to obtain an element in $I_{4,8}^{*}.$}
  \label{figf}
\end{figure}

In addition, we may essentially reverse this construction,
starting with an element of $I_{r,m}^{*}$ ($\ne \emptyset$), 
and so obtain an
element of $I_{r-m_r,d_r}$.    Indeed,
from the definitions of
$I_{r,m}^{*}$ and $\nu_r$ we know that

$$\nu_r = \sum_{j=1}^{r} \pi_{\tau(j)} =  \sum_{j=1}^{r} \sum_{\ell=j}^{m-r+j} \pi_{\ell} |\Delta t_{j,\ell}|,$$

\noindent for any $(t_{.,.}) \in I_{r,m}^{*}$.
However, we also have

\begin{align*}
\sum_{j=1}^{r} \sum_{\ell=j}^{m-r+j} \pi_{\ell} |\Delta t_{j,\ell}| 
&= {\bf 1}_{\{m_r > 0\}} \Bigl( \sum_{j=1}^{r} \sum_{\ell=j}^{m-r+j} {\bf 1}_{\{\pi_{\tau(\ell)} \ge \pi_{\tau(m_r)}\}}  \pi_{\ell} |\Delta t_{j,\ell}|\nonumber\\
&  \qquad + \sum_{j=1}^{r} \sum_{\ell=j}^{m-r+j} {\bf 1}_{\{\pi_{\tau(\ell)} < \pi_{\tau(m_r)}\}}  \pi_{\ell} |\Delta t_{j,\ell}| \Bigr) \nonumber\\
& \qquad + {\bf 1}_{\{m_r = 0\}} \pi_{\tau(1)} \sum_{j=1}^{r} \sum_{\ell=j}^{m-r+j} |\Delta t_{j,\ell}|\nonumber\\
&\le {\bf 1}_{\{m_r > 0\}}( (\pi_{\tau(1)} + \cdots + \pi_{\tau(m_r)}) + (r-m_r)\pi_{\tau(r)})\nonumber\\
& \qquad +  {\bf 1}_{\{m_r = 0\}} r \pi_{\tau(1)} \nonumber\\
&= \nu_r,
\end{align*}

\noindent with equality holding throughout if and only if
$m_r = 0$ or $m_r > 0$ and 
$\sum_{j=1}^r |\Delta t_{j,\ell}| = 1,$
for all $\ell$ such that $\pi_{\tau(\ell)} \ge \pi_{\tau(m_r)}$,
and that, moreover,
$\sum_{j=1}^{r} \sum_{\ell=j}^{m-r+j} {\bf 1}_{\{\pi_{\tau(\ell)} = \pi_{\tau(r)}\}} |\Delta t_{j,\ell}| = r - m_r$.
If $m_r > 0$, then, for any $(t_{.,.}) \in I_{r,m}^{*}$,
we may start with \eqref{this4},
and use again the permutation of the indices employed there.
We thus obtain the first term of \eqref{this4a},
which corresponds to the condition 
$\sum_{j=1}^r |\Delta t_{j,\ell}| = 1,$
for all $\ell$ such that $\pi_{\tau(\ell)} \ge \pi_{\tau(m_r)}$,
and also the second term of \eqref{this4a},
which corresponds to the other condition
$\sum_{j=1}^{r} \sum_{\ell=j}^{m-r+j} {\bf 1}_{\{\pi_{\tau(\ell)} = \pi_{\tau(r)}\}} |\Delta t_{j,\ell}| = r - m_r$.
If $m_r = 0$ the same reasoning holds, except that the first
term in \eqref{this4a} is taken to be zero.

Having thus established
a bijection between
$I_{r,m}^{*}$ and $I_{r-m_r,d_r}$,
we may thus maximize over these two parameter sets, and so,
for any process $(M(t))_{t\ge 0}$, obtain the general result

\begin{align}\label{this5}
&\sum_{i=1}^{m_r} M^{\tau(i)}(1) 
  + \max_{I_{r-m_r,d_r}} \sum_{j=1}^{(r-m_r)} \sum_{\ell-j}^{(r-m_r+d_r-1)} M^{\tau(m_r + \ell)}(\Delta u_{j,\ell})\nonumber\\
&\qquad =  \max_{I_{r,m}^{*}} \sum_{j=1}^{r}  \sum_{\ell=1}^{m-r+j}         M^{\tilde{\ell}(j,q)}(\Delta t_{j,\ell}).
\end{align}

We now proceed to show that \eqref{this3} holds.
First, fix $c > 0$, and, for each $1 \le \ell \le m$, set

\begin{equation}\label{this6}
c_{\ell} =   \begin{cases}
  c, \qquad \text{if $\pi_{\ell} < \pi_{\tau(r)}$,}\\
  0, \qquad \text{otherwise.}
  \end{cases}
\end{equation}

Next, let 
$\widehat{M}_n^{\ell}(t) = \sigma_{\ell} \hat{B}_n^{\ell}(t) - c_{\ell}t$,
and let 
$M^{\ell}(t) = \sigma_{\ell} \tilde B^{\ell}(t) - c_{\ell}t$.
Then, for $n$ large enough, namely, 
for $n > c/(\pi_{\tau(r)} - \pi_{\tau(r+1)})$, 
we have that, almost surely,
for any $t_{.,.} \in I_{r,m}$,

\begin{align}\label{this7}
&\sum_{j=1}^{r}  \sum_{\ell=1}^{m-r+j} \widehat{M}_n^{\ell}(\Delta t_{j,\ell})\nonumber\\
&\qquad \qquad \ge \sum_{j=1}^{r} \sum_{\ell=j}^{m-r+j} 
  \left(\sigma_{\ell}  \hat{B}_n^{\ell}(\Delta t_{j,\ell}) 
- \sqrt{n} \left( \pi_{\tau(j)} - \pi_{\ell}\right) |\Delta t_{j,\ell}|\right).
\end{align}

\noindent Hence, almost surely, both

\begin{align}\label{this8}
&\max_{I_{r,m}} \sum_{j=1}^{r}  \sum_{\ell=1}^{m-r+j} \widehat{M}_n^{\ell}(\Delta t_{j,\ell})\nonumber\\
&\qquad \qquad \ge \max_{I_{r,m}}  \sum_{j=1}^{r} \sum_{\ell=j}^{m-r+j} \left( \sigma_{\ell} \hat{B}_n^{\ell}(\Delta t_{j,\ell}) 
- \sqrt{n} \left( \pi_{\tau(j)} - \pi_{\ell}\right) |\Delta t_{j,\ell}| \right),
\end{align}

\noindent and

\begin{align}\label{this9}
\max_{I_{r,m}^{*}} \sum_{j=1}^{r}  \sum_{\ell=1}^{m-r+j} \widehat{M}_n^{\ell}(\Delta t_{j,\ell})
  = \max_{I_{r,m}^{*}} \sum_{j=1}^{r} \sum_{\ell=j}^{m-r+j} \sigma_{\ell} \hat{B}_n^{\ell}(\Delta t_{j,\ell}).
\end{align}

Now choose any $z > 0$.  Then

\begin{align}\label{this10}
&\bbp \biggl(  \max_{I_{r,m}}  \sum_{j=1}^{r} \sum_{\ell=j}^{m-r+j} \left(\sigma_{\ell} \hat{B}_n^{\ell}(\Delta t_{j,\ell}) 
      - \sqrt{n} \left( \pi_{\tau(j)} - \pi_{\ell}\right) |\Delta t_{j,\ell}|\right)\nonumber\\
& \qquad   - \max_{I_{r,m}^{*}} \sum_{j=1}^{r} \sum_{\ell=j}^{m-r+j} \sigma_{\ell} \hat{B}_n^{\ell}(\Delta s_{j,\ell}) > z \biggr) \nonumber\\
& \le \bbp\left(\max_{I_{r,m}} \sum_{j=1}^{r}  \sum_{\ell=1}^{m-r+j} \widehat{M}_n^{\ell}(\Delta t_{j,\ell})
  - \max_{I_{r,m}^{*}} \sum_{j=1}^{r}  \sum_{\ell=1}^{m-r+j} \widehat{M}_n^{\ell}(\Delta t_{j,\ell}) > z \right),
\end{align}

\noindent so that

{\allowdisplaybreaks
\begin{align}\label{this11}
&\limsup_{n\rightarrow \infty} 
  \bbp \biggl( \max_{I_{r,m}}  \sum_{j=1}^{r} \sum_{\ell=j}^{m-r+j} \left(\sigma_{\ell} \hat{B}_n^{\ell}(\Delta t_{j,\ell}) 
       - \sqrt{n} \left( \pi_{\tau(j)} - \pi_{\ell}\right) |\Delta t_{j,\ell}|\right)\nonumber\\
&\qquad \qquad \qquad  - \max_{I_{r,m}^{*}} \sum_{j=1}^{r} \sum_{\ell=j}^{m-r+j} \sigma_{\ell} \hat{B}_n^{\ell}(\Delta s_{j,\ell}) > z \biggr)\nonumber\\
&\ \le \limsup_{n\rightarrow \infty} 
   \bbp\biggl( \max_{I_{r,m}} \sum_{j=1}^{r}  \sum_{\ell=1}^{m-r+j} \widehat{M}_n^{\ell}(\Delta t_{j,\ell})
    - \max_{I_{r,m}^{*}} \sum_{j=1}^{r}  \sum_{\ell=1}^{m-r+j} \widehat{M}_n^{\ell}(\Delta t_{j,\ell}) > z  \biggr)\nonumber\\
&=   \bbp\left(\max_{I_{r,m}} \sum_{j=1}^{r}  \sum_{\ell=1}^{m-r+j} M^{\ell}(\Delta t_{j,\ell})
     - \max_{I_{r,m}^{*}} \sum_{j=1}^{r}  \sum_{\ell=1}^{m-r+j} M^{\ell}(\Delta t_{j,\ell}) > z  \right),
\end{align}
}

\noindent by the Invariance Principle and the Continuous Mapping Theorem.  
Next, for any $0 \le \varepsilon \le 1$, let

$$I_{r,m}(\varepsilon) = \{(t_{j,\ell}) \in I_{r,m}: \sum_{j,\ell} |\Delta t_{j,\ell}|{\bf 1}_{\{\pi_{\ell} < \pi_{\tau(r)}\}} \le \varepsilon r\}.$$

\noindent Thus, 
$I_{r,m}^{*} = I_{r,m}(0) \subset I_{r,m}(\varepsilon) \subset I_{r,m}(1) = I_{r,m}$.
We bound \eqref{this11} using this family of subsets as follows:

{\allowdisplaybreaks
\begin{align}\label{this12}
& \bbp\left(\max_{I_{r,m}} \sum_{j=1}^{r}  \sum_{\ell=1}^{m-r+j} M^{\ell}(\Delta t_{j,\ell})
    - \max_{I_{r,m}^{*}} \sum_{j=1}^{r}  \sum_{\ell=1}^{m-r+j} M^{\ell}(\Delta t_{j,\ell}) > z  \right)\nonumber\\
&\le \bbp\left(\max_{I_{r,m}(\varepsilon)} \sum_{j=1}^{r}  \sum_{\ell=1}^{m-r+j} M^{\ell}(\Delta t_{j,\ell})
    - \max_{I_{r,m}^{*}} \sum_{j=1}^{r}  \sum_{\ell=1}^{m-r+j} M^{\ell}(\Delta t_{j,\ell}) > z  \right)\nonumber\\
&\quad +  \bbp\left(\max_{I_{r,m} \backslash I_{r,m}(\varepsilon)} \sum_{j=1}^{r}  \sum_{\ell=1}^{m-r+j} M^{\ell}(\Delta t_{j,\ell})
    - \max_{I_{r,m}^{*}} \sum_{j=1}^{r}  \sum_{\ell=1}^{m-r+j} M^{\ell}(\Delta t_{j,\ell}) > z  \right)\nonumber\\
&\le \bbp\left(\max_{I_{r,m}(\varepsilon)} \sum_{j=1}^{r}  \sum_{\ell=1}^{m-r+j} \tilde B^{\ell}(\Delta t_{j,\ell})
    - \max_{I_{r,m}^{*}} \sum_{j=1}^{r}  \sum_{\ell=1}^{m-r+j} \tilde B^{\ell}(\Delta s_{j,\ell}) > z  \right)\nonumber\\
&\quad +  \bbp\left(\max_{I_{r,m} \setminus I_{r,m}(\varepsilon)} \sum_{j=1}^{r}  \sum_{\ell=1}^{m-r+j} \tilde B^{\ell}(\Delta t_{j,\ell})
    - \max_{I_{r,m}^{*}} \sum_{j=1}^{r}  \sum_{\ell=1}^{m-r+j} \tilde B^{\ell}(\Delta s_{j,\ell}) > z + \varepsilon r c  \right)\nonumber\\
&\le \bbp\left(\max_{I_{r,m}(\varepsilon)} \sum_{j=1}^{r}  \sum_{\ell=1}^{m-r+j} \tilde B^{\ell}(\Delta t_{j,\ell})
    - \max_{I_{r,m}^{*}} \sum_{j=1}^{r}  \sum_{\ell=1}^{m-r+j} \tilde B^{\ell}(\Delta s_{j,\ell}) > z  \right)\nonumber\\
&\quad +  \bbp \biggl( \max_{I_{r,m}} \sum_{j=1}^{r}  \sum_{\ell=1}^{m-r+j} \tilde B^{\ell}(\Delta t_{j,\ell})
- \max_{I_{r,m}^{*}} \sum_{j=1}^{r}  \sum_{\ell=1}^{m-r+j} \tilde B^{\ell}(\Delta s_{j,\ell}) > z + \varepsilon r c  \biggr).
\end{align}
}

We can now take the limsup in \eqref{this12}, as
$c \rightarrow \infty$, and then, as
$\varepsilon \rightarrow 0$, and
so establish convergence to zero in probability.
Moreover, since

$$\bbp\left(\max_{I_{r,m}} \sum_{j=1}^{r}  \sum_{\ell=1}^{m-r+j} M^{\ell}(\Delta t_{j,\ell})
    - \max_{I_{r,m}^{*}} \sum_{j=1}^{r}  \sum_{\ell=1}^{m-r+j} M^{\ell}(\Delta t_{j,\ell}) \ge 0  \right) = 1,$$

\noindent we have in fact shown, with the help of \eqref{this11}, 
that with probability one,

\begin{equation*}
\max_{I_{r,m}} \sum_{j=1}^{r}  \sum_{\ell=1}^{m-r+j} M^{\ell}(\Delta t_{j,\ell})
=\max_{I_{r,m}^{*}} \sum_{j=1}^{r}  \sum_{\ell=1}^{m-r+j} M^{\ell}(\Delta t_{j,\ell}),
\end{equation*}

\noindent and thus

\begin{align}\label{this13}
&\max_{I_{r,m}}  \sum_{j=1}^{r} \sum_{\ell=j}^{m-r+j} \left(\sigma_{\ell} \hat{B}_n^{\ell}(\Delta t_{j,\ell}) 
       - \sqrt{n} \left( \pi_{\tau(j)} - \pi_{\ell}\right) |\Delta t_{j,\ell}|  \right)   \nonumber\\
&\qquad \qquad  - \max_{I_{r,m}^{*}} \sum_{j=1}^{r} \sum_{\ell=j}^{m-r+j} \sigma_{\ell} \hat{B}_n^{\ell}(\Delta s_{j,\ell}) 
\stackrel{\bbp}{\rightarrow} 0.
\end{align}

\noindent Since

\begin{equation}\label{this13b}
\max_{I_{r,m}^{*}} \sum_{j=1}^{r} \sum_{\ell=j}^{m-r+j} \sigma_{\ell} \hat{B}_n^{\ell}(\Delta s_{j,\ell})
\Rightarrow \max_{I_{r,m}^{*}} \sum_{j=1}^{r} \sum_{\ell=j}^{m-r+j} \sigma_{\ell} \tilde B^{\ell}(\Delta s_{j,\ell}),
\end{equation}

\noindent by the Converging Together Lemma, we have 
proved \eqref{this3}.  Equation \eqref{item7r}
of the theorem follows from the bijection between
$I_{r,m}^{*}$ and $I_{r-m_r,d_r}$ described in
the general result \eqref{this5}.

Finally, we can obtain the convergence of the joint distribution
in \eqref{item7saa} in the following
manner.  Given any $(\theta_1,\theta_2,\dots,\theta_r) \in \bbr^r$,
we have

{\allowdisplaybreaks
\begin{align}\label{this13c}
&\sum_{k=1}^{r} \theta_k \left( \frac{V^k_n -\nu_k n}{\sqrt{n}} \right)\nonumber\\
&\qquad=  \sum_{k=1}^{r} \theta_k \biggl( \max_{I_{k,m}}  \sum_{j=1}^{k} \sum_{\ell=j}^{m-k+j} \left(\sigma_{\ell} \hat{B}_n^{\ell}(\Delta t_{j,\ell}) 
       - \sqrt{n} \left( \pi_{\tau(j)} - \pi_{\ell}\right) |\Delta t_{j,\ell}|  \right)   \biggr) \nonumber\\
&\qquad=  \sum_{k=1}^{r} \theta_k \biggl( \max_{I_{k,m}}  \sum_{j=1}^{k} \sum_{\ell=j}^{m-k+j} \left(\sigma_{\ell} \hat{B}_n^{\ell}(\Delta t_{j,\ell}) 
       - \sqrt{n} \left( \pi_{\tau(j)} - \pi_{\ell}\right) |\Delta t_{j,\ell}|  \right)   \nonumber\\
&\qquad - \max_{I_{k,m}^{*}} \sum_{j=1}^{k} \sum_{\ell=j}^{m-k+j} \sigma_{\ell} \hat{B}_n^{\ell}(\Delta s_{j,\ell}) \biggr)
+ \sum_{k=1}^{r} \theta_k \biggl( \max_{I_{k,m}^{*}} \sum_{j=1}^{k} \sum_{\ell=j}^{m-k+j} \sigma_{\ell} \hat{B}_n^{\ell}(\Delta s_{j,\ell})   \biggr).
\end{align}
}

Now from \eqref{this13}, the first summation on the right-hand side of \eqref{this13c}
converges to zero in probability, as $n \rightarrow \infty$.  Moreover,
the second summation 
is a continuous functional of $(\hat{B}_n^{1},\hat{B}_n^{2},\dots,\hat{B}_n^{m})$,
and so, by the Invariance Principle and Continuous Mapping Theorem,
converges.  Then the Converging Together Lemma, along with
the bijection result \eqref{this5}, gives

\begin{align}\label{this13d}
&\sum_{k=1}^{r} \theta_k \left( \frac{V^k_n -\nu_k n}{\sqrt{n}} \right)\nonumber\\
&\qquad \Rightarrow \sum_{k=1}^{r} \theta_k \biggl( \max_{I_{k,m}^{*}} \sum_{j=1}^{k} \sum_{\ell=j}^{m-k+j} \sigma_{\ell} \tilde B^{\ell}(\Delta s_{j,\ell})  \biggr)
= \sum_{k=1}^{r} \theta_k V^k_{\infty}.
\end{align}

Since \eqref{this13d} holds for arbitrary $(\theta_1,\theta_2,\dots,\theta_r) \in \bbr^r$,
by the Cram\'er-Wold Theorem, we have the joint convergence result
\eqref{item7saa}.\CQFD
\end{Proof}

Since the shape of the Young tableau is
more naturally expressed in terms of the
$R_n^k$, rather than of the $V_n^k$,
we may restate the results of the
previous theorem as follows:

\begin{theorem}\label{thm5}
Let $(X_n)_{n \ge 0}$ be an irreducible, aperiodic,
homogeneous Markov chain with finite state space
${\cal A}_m = \{\alpha_1 < \cdots < \alpha_m\}$,
and with stationary distribution
$(\pi_1,\pi_2,\dots,\pi_m)$.
Then, in the notations of Theorem~\ref{thm4},

\begin{align}\label{this14}
\left( \frac{R_n^1 -\pi_{\tau(1)}n}{\sqrt{n}},\frac{R_n^2 -\pi_{\tau(2)}n}{\sqrt{n}}, \dots, \frac{R_n^m -\pi_{\tau(m)}n}{\sqrt{n}} \right) 
    \Rightarrow (R_{\infty}^1, R_{\infty}^2, \dots, R_{\infty}^m),
\end{align}

\noindent where

\begin{align}\label{this15a}
R_{\infty}^1
&=  \max_{I_{1,d_1}} 
    \sum_{\ell=1}^{d_1} \sigma_{\tau(\ell)} \left(\tilde B^{\tau(\ell)}(t_{1,\ell})  - \tilde B^{\tau(\ell)}(t_{1,\ell-1}) \right),
\end{align}

\noindent and, for each $2 \le k \le m$,

\begin{align}\label{this15b}
&R_{\infty}^k
=  \sum_{i=m_{k-1}+1}^{m_k} \sigma_{\tau(i)} \tilde B^{\tau(i)}(1)  \nonumber\\
&\quad  + \max_{I_{k-m_k,d_k}} 
   \sum_{j=1}^{k-m_k} 
     \sum_{\ell=j}^{(d_k + m_k- k+j)}
       \sigma_{\tau(m_k + \ell)} \tilde B^{\tau(m_k + \ell)}(\Delta t_{j,\ell})\nonumber\\
&\quad  - \max_{I_{k-1-m_{k-1},d_{k-1}}} 
   \sum_{j=1}^{k-1-m_{k-1}} 
     \sum_{\ell=j}^{(d_{k-1} + m_{k-1}- k+1+j)}
       \sigma_{\tau(m_{k-1} + \ell)} \tilde B^{\tau(m_{k-1} + \ell)}(\Delta t_{j,\ell}),
\end{align}

\noindent where we use the notation
$\tilde B^{s}(\Delta t_{j,\ell}) = \tilde B^{s}(t_{j,\ell}) - \tilde B^{s}(t_{j,\ell-1})$,
for any $1 \le s \le m$, $1 \le j \le k$, and $1 \le \ell \le m$,
and where the first sum on the right-hand side
of \eqref{this15b} is understood to be $0$, if $m_k = m_{k-1}$.

\end{theorem}

\noindent \begin{Proof}
First, $R_n^1 = V_n^1$, and,
for each $2 \le k \le m$, 
$R_n^k = V_n^k - V_n^{k-1}$.
Expressing these equalities
at the multivariate level,
we have

{\allowdisplaybreaks
\begin{align}\label{this15c}
&\left( \frac{R_n^1 -\pi_{\tau(1)}n}{\sqrt{n}},\frac{R_n^2 -\pi_{\tau(2)}n}{\sqrt{n}}, \dots, \frac{R_n^m -\pi_{\tau(m)}n}{\sqrt{n}} \right)  \nonumber\\
&\qquad\qquad= \left( \frac{V_n^1 -\pi_{\tau(1)}n}{\sqrt{n}},\frac{V_n^2 - V_n^1 -\pi_{\tau(2)}n}{\sqrt{n}}, \dots, \frac{V_n^m -V_n^{m-1} -\pi_{\tau(m)}n}{\sqrt{n}} \right) \nonumber\\
&\qquad\qquad= \left( \frac{V_n^1 -\nu_1 n}{\sqrt{n}},\frac{V_n^2 -\nu_2 n}{\sqrt{n}}, \dots, \frac{V_n^m - \nu_m n}{\sqrt{n}} \right) \nonumber\\
&\qquad\qquad\qquad - \left(0,\frac{V_n^1 -\nu_1 n}{\sqrt{n}}, \dots, \frac{V_n^m - \nu_{m-1} n}{\sqrt{n}} \right) \nonumber\\
&\qquad\qquad\Rightarrow (V_{\infty}^1, V_{\infty}^2,\dots, V_{\infty}^m) - (0, V_{\infty}^1, \dots, V_{\infty}^m) \nonumber\\
&\qquad\qquad:= (R_{\infty}^1, R_{\infty}^2, \dots, R_{\infty}^m),
\end{align}
}

\noindent where the weak convergence
follows immediately 
from the Continuous Mapping Theorem,
since the transformation is linear.

Equations \eqref{this15a} and \eqref{this15b}
follow simply from the Brownian expressions
for $(V_{\infty}^1, V_{\infty}^2,\dots, V_{\infty}^m)$
in Theorem \ref{thm4}.\CQFD
\end{Proof}

If all $m$ letters have unique stationary
probabilities, 
then we have the following corollary to
Theorem~\ref{thm5}:

\begin{corollary}\label{cor1c}
If the stationary distribution of Theorem~\ref{thm5}
is such that each $\pi_r$ is unique, then

\begin{align}\label{this15d}
\left( \frac{R_n^1 -\pi_{\tau(1)}n}{\sqrt{n}},\frac{R_n^2 -\pi_{\tau(2)}n}{\sqrt{n}}, \dots, \frac{R_n^m -\pi_{\tau(m)}n}{\sqrt{n}} \right) 
    \Rightarrow N((0,0,\dots,0),\Sigma).
\end{align}

In other words, the limiting distribution is identical in law to
the spectrum of the diagonal matrix
$D = diag\{Z_1,Z_2,\dots,Z_m\}$, where
$(Z_1,Z_2,\dots,Z_m)$ is a centered normal
random vector with covariance matrix $\Sigma$.
\end{corollary}

\noindent \begin{Proof}
Now, for all $1 \le k \le m$, 
$d_k = 1 $, and
$m_k = k-1$, so that

\begin{align*}
R_{\infty}^1
&=  \max_{I_{1,d_1}} 
    \sum_{\ell=1}^{d_1} \sigma_{\tau(\ell)} \left(\tilde B^{\tau(\ell)}(t_{1,\ell})  - \tilde B^{\tau(\ell)}(t_{1,\ell-1}) \right)\nonumber\\
&=  \sigma_{\tau(1)}\tilde B^{\tau(1)}(1),
\end{align*}

\noindent and, for each $2 \le k \le m$,

{\allowdisplaybreaks
\begin{align*}
&R_{\infty}^k
=  \sum_{i=m_{k-1}+1}^{m_k} \sigma_{\tau(i)} \tilde B^{\tau(i)}(1)  \nonumber\\
&\quad  + \max_{I_{k-m_k,d_k}} 
   \sum_{j=1}^{k-m_k} 
     \sum_{\ell=j}^{(d_k + m_k- k+j)}
       \sigma_{\tau(m_k + \ell)} \tilde B^{\tau(m_k + \ell)}(\Delta t_{j,\ell})\nonumber\\
&\quad  - \max_{I_{k-1-m_{k-1},d_{k-1}}} 
   \sum_{j=1}^{k-1-m_{k-1}} 
     \sum_{\ell=j}^{(d_{k-1} + m_{k-1}- k+1+j)}
       \sigma_{\tau(m_{k-1} + \ell)} \tilde B^{\tau(m_{k-1} + \ell)}(\Delta t_{j,\ell})\nonumber\\
&\quad=  \sum_{i=k-1}^{k-1} \sigma_{\tau(i)} \tilde B^{\tau(i)}(1)  \nonumber\\
&\quad  + \max_{I_{1,1}} 
   \sum_{j=1}^{1} 
     \sum_{\ell=j}^{j}
       \sigma_{\tau(k - 1 + \ell)} \tilde B^{\tau(k - 1 + \ell)}(\Delta t_{j,\ell})\nonumber\\
&\quad  - \max_{I_{1,1}} 
   \sum_{j=1}^{1} 
     \sum_{\ell=j}^{j}
       \sigma_{\tau(k-2 + \ell)} \tilde B^{\tau(k-2 + \ell)}(\Delta t_{j,\ell})\nonumber\\
&\quad=   \sigma_{\tau(k-1)} \tilde B^{\tau(k-1)}(1) 
         +\sigma_{\tau(k)} \tilde B^{\tau(k)}(1)
         -\sigma_{\tau(k-1)} \tilde B^{\tau(k-1)}(1)\nonumber\\
&\quad=   \sigma_{\tau(k)} \tilde B^{\tau(k)}(1).
\end{align*}
}

\noindent Moreover, the joint law result
 for $(R_{\infty}^1,R_{\infty}^2,\dots,R_{\infty}^m)$
holds as well, and this is clearly a multivariate normal distribution, with mean
$(0,0,\dots,0)$ and covariance matrix $\Sigma$.  
Since the spectrum of a diagonal matrix consists of its diagonal elements,
the final claim of the corollary holds.\CQFD
\end{Proof}

\begin{Rem}\label{RMTconnect}
We know that the joint law of 
$(R_{\infty}^1, R_{\infty}^2, \dots, R_{\infty}^m)$
in the iid uniform alphabet case is identical to the
joint law of the eigenvalues of an
$m \times m$ traceless GUE matrix.
Corollary \ref{cor1c} also gives a spectral characterization
for the unique probability case, in particular,
for a non-uniform iid alphabet with unique stationary
probabilities.  
This is consistent with
the characterization of the limiting law of
$LI_n$ in the non-uniform iid case, due to 
Its, Tracy, and Widom \cite{ITW1,ITW2}, as that of the
largest eigenvalue of the block associated with
the most probable letters among a
direct sum of independent GUE matrices
whose dimensions correspond
to the multiplicities $d_r$ of
Theorems \ref{thm4} and \ref{thm5}, 
subject to the condition
that $\sum_{r=1}^m \sqrt{\pi_{\tau(r)}}X_r = 0$,
where $X_1,X_2,\dots,X_m$ are the diagonal
elements of the random matrix.
\end{Rem}

\begin{Rem}\label{GenTrace}
The difference between the zero-trace condition
$\sum_{r=1}^m X_r = 0$ and
the generalized traceless condition
$\sum_{r=1}^m \sqrt{\pi_{\tau(r)}} X_r = 0$
amounts to nothing more than a difference
in the choice of scaling for each row
$R^{r}_n$.
We will find it more natural to express our
results using the normalization associated
with the zero-trace condition
$\sum_{r=1}^m X_r = 0$
\end{Rem}

\section{Fine Structure of the Brownian Functional}

So far, we have seen that the limiting shape of the random Young
tableau generated by an aperiodic, irreducible, homogeneous
Markov chain can be expressed as a Brownian functional.
The form of this functional
is similar to the iid case;
the only difference is in the covariance
structure of the Brownian motion.
We begin our study of the consequences of this difference.

In the iid uniform $m$-alphabet case, Johansson \cite{Jo} proved
that the limiting shape of the Young tableau had a
joint law which is that of the spectrum of an $m \times m$
traceless GUE matrix. An immediate consequence of this
result is that the limiting shape of the Young tableau 
contains simple symmetries, {\it e.g.,}  for each $1 \le r \le m$,
$R^r_{\infty} \stackrel{\cal{L}}{=} - R^{m-r}_{\infty}.$
Now, as was seen in Corollary \ref{cor1a} 
of Theorem \ref{thm4}, the form of the Brownian  functional 
in the doubly stochastic case involved only the maximal term.
We will see that that there
is also a pleasing symmetry to the limiting shape of 
Young tableaux in the doubly stochastic case 
by examining a natural bijection
between the parameter set $I_{r,m}$ and
$I_{m-r,m}$, for any $1 \le r \le m-1.$
Indeed, this result will follow as a corollary
to the following, more general, theorem:

\begin{theorem}\label{thm6}
The limiting functionals of Theorem \ref{thm4} enjoy
the following symmetry property:
for every $1 \le r \le m-1$,

{\allowdisplaybreaks
\begin{align}\label{this16}
&V^r_{\infty} := \sum_{i=1}^{m_r} \sigma_{\tau(i)} \tilde B^{\tau(i)}(1)  \nonumber\\
&\qquad  + \max_{t(\cdot,\cdot) \in I_{r-m_r,d_r}} 
   \sum_{j=1}^{r-m_r} 
   \sum_{\ell=j}^{(m_r+d_r-r+j)} \sigma_{\tau(m_r + \ell)} \tilde B^{\tau(m_r + \ell)}(\Delta t_{j,\ell})\nonumber\\
&\qquad \stackrel{\cal{L}}{=}
\sum_{i=m_r+d_r+1}^{m} \sigma_{\tau(i)} \tilde B^{\tau(i)}(1)  \nonumber\\
&\qquad + \max_{u(\cdot,\cdot) \in I_{m_r+d_r-r,d_r}} 
     \sum_{j=1}^{m_r+d_r-r} 
     \quad \sum_{\ell=j}^{r-m_r+j} \sigma_{\tau(m_r + \ell)} \tilde B^{\tau(m_r + \ell)}(\Delta u_{j,\ell}),
\end{align}
}

\noindent where $\tilde B^{\ell}(\Delta) :=
\tilde B^{\ell}(t) - \tilde B^{\ell}(s)$, for $\Delta = [s,t]$,
and where the non-maximal terms on the left and right-hand
sides of \eqref{this16} are identically zero if $m_r = 0$,
or $m_r + d_r = m$, respectively. 

\end{theorem}

\begin{Rem}\label{Remsym}
Recall that, from the definitions of $m_r$ and $d_r$,
the non-maximal summation terms on the left and 
right-hand sides of \eqref{this16}
reflect the letters which have, respectively, greater and smaller
stationary probabilities than $\pi_{\tau(r)}$.
Recall, moreover, that the maximal terms are associated with the indices
having the same stationary probability as $\pi_{\tau(r)}$.
The maximal term on the left-hand side of \eqref{this16}
involves a summation over $r-m_r$ rows,
while the one on the right-hand side involves
$m_{r+1}-r$ rows.  Thus, in a sense, the two maximal
terms in \eqref{this16} split $d_r = m_{r+1}-m_r$ rows between
themselves.  In summary, the functional on the
right-hand side of \eqref{this16} corresponds to
the sum of the $m-r$ {\it bottom} rows of the 
Young tableau.
\end{Rem}

\noindent \begin{Proof} 
Without loss of generality, we may assume that
$\tau(j) = j$, for all $1 \le j \le m$.
Fix $1 \le r \le m-1$, and
for any point $t$ in the index set $I_{r-m_r,d_r}$,
define $\Delta t_{j+m_r,\ell} = [t_{j,\ell-1}, t_{j,\ell}]$,
for $1 \le j \le r-m_r$ and $1 \le \ell \le d_r$.
Furthermore, for each 
$1 \le j \le m_r$ or $m_{r+1} < j \le m$, 
set $\Delta t_{j,\ell} = [0,1]$, for $j = \ell$,
$\Delta t_{j,\ell} = \{0\}$, for $0 \le \ell < j$, and
$\Delta t_{j,\ell} = \{1\}$, for $j < \ell \le m.$
Next, as in the proof of Theorem \ref{thm4},
consider the {\it set} of points 
$\{ t_{j,\ell}\}_{(1 \le j \le r-m_r, 1 \le \ell \le d_r)}$,
and order them as
$s_0 := 0 < s_1 < \cdots < s_{\kappa-1} < s_{\kappa}:= 1$,
for some integer $\kappa$,
and let $\Delta s_q = [s_{q-1},s_q]$, 
for each $1 \le q \le \kappa$.

Now, for each $1 \le q \le \kappa$, 
let $A_q$ consist of the indices $\ell$ for which
$\Delta s_q \cap \Delta t_{j,\ell} \ne \emptyset$.
Then, almost surely, 

{\allowdisplaybreaks
\begin{align}\label{this18}
\sum_{i=1}^{m_r} \sigma_{i} \tilde B^{i}(1)  
&+ \sum_{j=1}^{r-m_r}    \sum_{\ell=j}^{(m_r+d_r-r+j)} \sigma_{m_r + \ell} \tilde B^{m_r + \ell}(\Delta t_{j,\ell})\nonumber\\
&\qquad = \sum_{j=1}^{r} \sum_{\ell=1}^{m} \sigma_{\ell} \tilde B^{\ell}(\Delta t_{j,\ell})\nonumber\\
&\qquad = \sum_{j=1}^{r} \sum_{q=1}^{\kappa}  \sum_{\ell=1}^{m} \sigma_{\ell} \tilde B^{\ell}(\Delta t_{j,\ell} \cap \Delta s_q)\nonumber\\
&\qquad = \sum_{j=1}^{r} \sum_{q=1}^{\kappa}  \sum_{\ell \in A_q} \sigma_{\ell} \tilde B^{\ell}(\Delta s_q).
\end{align}
}

Now by the ``stairstep'' properties of $I_{r,m}$
there are precisely
$r$ elements in each $A_q$.  
Letting  $\tilde{A}_q = \{1,\dots,m\} \setminus A_q$,
for each $1 \le q \le \kappa$, we thus see that each 
$\tilde{A}_q$ contains exactly $m-r$ elements.
Let  $\tilde{\ell}_{j,q}$ be the $j^{th}$ smallest
element of $\tilde{A}_q$.  We claim that for each $1 \le j \le m-r$,
the sequence 
$\tilde{\ell}_{j,1}, \tilde{\ell}_{j,2}, \dots, \tilde{\ell}_{j,\kappa}$.
is weakly decreasing. 

Indeed, fix $1 \le j \le m-r$ and $1 \le q \le \kappa-1$,
and suppose that $\tilde{\ell}_{j,q}$ is less than all
the elements of $A_q$.  Then, by the properties of $I_{r,m}$,
the least element of $A_{q+1}$ is no smaller, so that
the $j^{th}$ smallest element of $\tilde{A}_q$,
$\tilde{\ell}_{j,q+1}$ is also $\tilde{\ell}_{j,q}$.
Next, suppose that $\tilde{\ell}_{j,q}$ is
greater than $k \ge 1$ elements of $A_q$.
Thus, $\tilde{\ell}_{j,q} = j + k$.  Then 
there are at most $k$ elements of $A_{q+1}$ which
are less than or equal to 
$\tilde{\ell}_{j,q}$, by the properties
of $I_{r,m}$.  But this implies that there
are at least $j$ elements of $\tilde{A}_{q+1}$
which are less than or equal to 
$\tilde{\ell}_{j,q}$.  Thus,
$\tilde{\ell}_{j,q+1} \le \tilde{\ell}_{j,q}$,
and the claim is proved.

Moreover, since each $A_q$ contains 
$\{1,2,\dots,m_r\}$, we see that necessarily
each $\tilde{A}_q$ contains 
$\{m_r+d_r+1, m_r+d_r+2,\dots, m\}.$

For each $1 \le j \le m-r$,
we may now amalgamate the intervals $\Delta s_q$
to obtain a partition of the unit interval.
Specifically, for each $1 \le j \le m-r$,
and each $1 \le \ell \le m$, let $\tilde{u}_{j,\ell}$ be the
smallest $s_q$ such that
$\tilde{\ell}_{j,q+1} \le \ell$. 
(We define $\tilde{u}_{j,0} = 1$, for all $1 \le j \le m-r$.)

Finally, and most crucially, recall that
$\sum_{\ell=1}^m \sigma_{\ell} \tilde B^{\ell}(t) = 0$, for all $t$.  
Then since 
$(\tilde B^1,\tilde B^2,\dots,\tilde B^m)  
\stackrel{{\cal L}}{=} (-\tilde B^1,-\tilde B^2,\dots,-\tilde B^m)$,

{\allowdisplaybreaks
\begin{align}\label{this19}
&\sum_{j=1}^{r}   \sum_{q=1}^{\kappa}  \sum_{\ell \in A_q}         \sigma_{\ell} \tilde B^{\ell}(\Delta s_q)\nonumber\\
&\qquad = \sum_{j=1}^{m-r} \sum_{q=1}^{\kappa}  \sum_{\ell \in \tilde{A}_q} \left(-\sigma_{\ell} \tilde B^{\ell}(\Delta s_q)\right)\nonumber\\
&\qquad = -\sum_{i=m_r+d_r+1}^{m} \sigma_{i} \tilde B^{i}(1) 
           - \sum_{j=1}^{m_r+d_r-r} 
               \quad \sum_{\ell=1}^{m} \sigma_{m_r+\ell} \tilde B^{m_r+\ell}(\Delta u_{j,\ell})\nonumber\\
&\qquad \stackrel{\cal{L}}{=} 
\sum_{i=m_r+d_r+1}^{m} \sigma_{i} \tilde B^{i}(1) 
+ \sum_{j=1}^{m_r+d_r-r} 
     \quad \sum_{\ell=1}^{m}  \sigma_{m_r+\ell} \tilde B^{m_r+\ell}(\Delta u_{j,\ell}),
\end{align}
}

\noindent where 
$\Delta u_{j,\ell} = [u_{j,\ell-1},u_{j,\ell}]$. 
But, by the way we ordered each $A_q$, we must have 
$\Delta u_{j_1,\ell} \cap \Delta u_{j_2,\ell} = \emptyset$, for any $j_1 \ne j_2$.  Thus,
$u \in I_{m_r+d_r-r,d_r}$, and so we may restrict the summation over 
$\ell$ in \eqref{this19} to $\ell=j,\dots, r-m_r+j$, since the remaining terms are zero.
Equation \eqref{this16} follows immediately
by taking the maxima over $I_{r-m_r,d_r}$ 
and $I_{m_r+d_r-r,d_r}$ over the left-hand
and right-hand sides, respectively, of \eqref{this19}.\CQFD
\end{Proof}

For doubly stochastic transition matrices, the
symmetry is even more apparent:

\begin{corollary}\label{cor2}
Let the transition matrix $P$ of Theorem \ref{thm4} be doubly stochastic.
Then, for every $1 \le r \le m-1$,

\begin{align}\label{this19a}
&V_{\infty}^r : = \max_{t(\cdot,\cdot) \in I_{r,m}} \sum_{j=1}^{r} \sum_{\ell=j}^{m-r+j} 
\sigma_{\ell} \left(\tilde B^{\ell}(t_{j,\ell})  - \tilde B^{\ell}(t_{j,\ell-1})\right)\nonumber\\
&\qquad \stackrel{\cal{L}}{=}
\max_{u(\cdot,\cdot) \in I_{m-r,m}} \sum_{j=1}^{m-r} \sum_{\ell=j}^{r+j} 
\sigma_{\ell} \left(\tilde B^{\ell}(u_{j,\ell})  - \tilde B^{\ell}(u_{j,\ell-1})\right) := V_{\infty}^{m-r},
\end{align}

\noindent and so

\begin{align}\label{this19b}
\lim_{n \rightarrow \infty}   \frac{\sum_{j=1}^r      R_n^j - rn/m}{\sqrt{n}} 
\stackrel{\cal{L}}{=} \lim_{n \rightarrow \infty} \frac{rn/m - \sum_{j=m-r+1}^m  R_n^j}{\sqrt{n}}.
\end{align}

\noindent Moreover,

\begin{equation}\label{this19c}
(V_{\infty}^1,\dots,V_{\infty}^r) \stackrel{\cal{L}}{=}
(V_{\infty}^{m-1},\dots,V_{\infty}^{m-r}) .
\end{equation}

\end{corollary}

\noindent \begin{Proof} 
Since $m_r = 0$ and $d_r = m$ for all $1 \le r \le m$, 
the non-maximal terms on both sides of \eqref{this16}
disappear, and we have \eqref{this19a}.

To prove \eqref{this19b}, recall that $V_n^m = \sum_{j=1}^m R_n^j = n$.
Then, from the result just proved,

{\allowdisplaybreaks
\begin{align}\label{this20}
\frac{V_n^{m-r} - (m-r)n/m}{\sqrt{n}} 
&= \frac{\sum_{j=1}^{m-r} R_n^j- (m-r)n/m}{\sqrt{n}} \nonumber\\
&= \frac{\left(n - \sum_{j=m-r+1}^m R_n^j\right)- (m-r)n/m}{\sqrt{n}} \nonumber\\
&= \frac{rn/m - \sum_{j=m-r+1}^m R_n^j}{\sqrt{n}}\nonumber\\
&\Rightarrow V_{\infty}^{m-r} \stackrel{\cal{L}}{=} V_{\infty}^{r},
\end{align}
}

\noindent and we have established the claimed symmetry.

Finally, the extension of \eqref{this19a} to \eqref{this19c} 
follows from a standard Cram\'er-Wold argument.\CQFD
\end{Proof}

Turning again to the cyclic case, recall that,
for $m \ge 4$, the limiting shape of the Young tableau
in general differs from that of the iid uniform case.  The
following theorem characterizes the asymptotic covariance
matrices of such Markov chains.

\begin{theorem}\label{thm7}
Let $P$ be the $m \times m$ transition matrix
of an aperiodic, irreducible, cyclic Markov chain
on an $m$-letter, ordered alphabet, 
${\cal A}_m = \{\alpha_1 < \alpha_2 < \cdots < \alpha_m\}$, 
with 

\begin{equation}\label{this21}
P =
  \begin{pmatrix} 
  a_1 & a_m &\cdots & a_3 & a_2\\
  a_2 & a_1 &\ddots &&     a_3\\
  \vdots & \ddots & \ddots & \ddots &\vdots\\
  a_{m-1} && \ddots  & a_1  &a_m\\
  a_m &a_{m-1} &\cdots &a_2 &a_1
  \end{pmatrix}.
\end{equation}

\noindent Then, for $1 \le j \le m$,
$\lambda_j = \sum_{k=1}^m a_k \omega^{(k-1)(j-1)}$
is an eigenvalue of $P$,
where $\omega = exp (2\pi i/m)$ is the $m^{th}$
principal root of unity.  Moreover,
letting 
$\gamma_j =  \lambda_j/(1-\lambda_j)$,  
for $2 \le j \le m$, and
$\beta_j = \cos (2\pi j/m)$,
for $0 \le j \le m$, 
the asymptotic covariance matrix $\Sigma$
is given by:\\

\noindent For $m = 2m_0 + 1$,

\begin{equation}\label{this22}
\Sigma = \frac{m-1}{m^2}M^{(1)} + \frac{4}{m^2} \sum_{j=2}^{m_0+1} Re(\gamma_j) M^{(j)},
\end{equation}

\noindent and for $m = 2m_0$,

\begin{equation}\label{this23}
\Sigma = \frac{m-1}{m^2}M^{(1)} + \frac{4}{m^2} \sum_{j=2}^{m_0} Re(\gamma_j) M^{(j)} + \frac{2}{m^2} \gamma_{m_0+1} M^{(m_0+1)},
\end{equation}

\noindent where $M^{(j)}$ 
is an $m \times m$  Toeplitz matrix with entries
$(M^{(j)})_{k,\ell} = \beta_{(j-1)|k-\ell|}$, for
$2 \le j \le m$, 
and $(M^{(1)})_{k,\ell} = \delta_{k,\ell} - (1-\delta_{k,\ell})/(m-1)$,
for $j=1$.
\end{theorem}

\noindent \begin{Proof} It is straightforward, and classical,
to verify that, for each $1 \le j \le m$,
$(1,\omega^{j-1}, \omega^{2(j-1)},\dots, \omega^{(m-1)(j-1)})$
is a left eigenvector of $P$, with eigenvalue 
$\lambda_j = \sum_{k=1}^m a_k \omega^{(k-1)(j-1)}$.  We can thus
write our standard diagonalization of $P$ as $P = S^{-1} \Lambda S$,
where $\Lambda = diag(1,\lambda_2,\dots,\lambda_m)$,

\begin{equation}\label{this24}
S = 
  \begin{pmatrix} 
  1 & 1 &\cdots & 1 & 1\\
  1 & \omega   & \omega^2 &\cdots & \omega^{m-1}\\
  1 & \omega^2 & \omega^4 &\cdots & \omega^{2(m-1)}\\
  \vdots & \vdots & \ddots & \ddots &\vdots\\
  1 & \omega^{m-1} & \omega^{2(m-1)} &\cdots & \omega^{(m-1)^2}
  \end{pmatrix},
\end{equation}

\noindent and 

{\allowdisplaybreaks
\begin{equation}\label{this25}
S^{-1} = 
  \frac{1}{m}
  \begin{pmatrix} 
  1 & 1 &\cdots & 1 & 1\\
  1 & \omega^{-1}   & \omega^{-2} &\cdots & \omega^{-(m-1)}\\
  1 & \omega^{-2} & \omega^{-4} &\cdots & \omega^{-2(m-1)}\\
  \vdots & \vdots & \ddots & \ddots &\vdots\\
  1 & \omega^{-(m-1)} & \omega^{-2(m-1)} &\cdots & \omega^{-(m-1)^2}
  \end{pmatrix}.
\end{equation}
}

In the present cyclic, and hence, doubly stochastic case, we know
that $\Sigma = (1/m)(I + S^{-1}DS + (S^{-1}DS)^T)$,
where, as usual, 
$D =diag(\gamma_1,\gamma_2,\dots, \gamma_m)$
$= diag(-1/2,\lambda_2/(1-\lambda_2),\dots, \lambda_m/(1-\lambda_m))$.
We can then compute the entries of
$S^{-1}DS$ as follows:

{\allowdisplaybreaks
\begin{align}\label{this26}
(S^{-1}DS )_{j_1,j_2}
  &=  \sum_{k,\ell} (S^{-1})_{j_1,k} (D)_{k,\ell} (S)_{\ell,j_2}\nonumber\\
  &=  \sum_{k,\ell}  \frac{1}{m} (\omega^{-(j_1-1)(k-1)}) (\delta_{k,\ell} \gamma_k) (\omega^{(j_2-1)(\ell-1)})\nonumber\\
  &=  \sum_{k=1}^m   \frac{\gamma_k}{m} \omega^{(j_2-j_1)(k-1)}\nonumber\\
  &=  \frac{1}{m} \left(-\frac{1}{2} + \sum_{k=2}^m   \gamma_k \omega^{(j_2-j_1)(k-1)}\right),
\end{align}
}

\noindent for all $1 \le j_1, j_2, \le m$.
The entries of the asymptotic covariance matrix 
can thus be written as

\begin{align}\label{this27}
\sigma_{j_1,j_2}
  &= \frac{1}{m}\left(\delta_{j_1,j_2} + (S^{-1}DS )_{j_1,j_2} + (S^{-1}DS )_{j_2,j_1}\right)\nonumber\\
  &= \frac{1}{m}\left(\delta_{j_1,j_2} + \frac{1}{m} \left(-1 + \sum_{k=2}^m  \gamma_k ( \omega^{(j_2-j_1)(k-1)} + \omega^{(j_1-j_2)(k-1)} ) \right) \right)\nonumber\\
  &= \frac{m-1}{m^2} M^{(1)}_{j_1,j_2} + \frac{2}{m^2} \sum_{k=2}^m  \gamma_k  \beta_{|j_2-j_1|(k-1)} ,
\end{align}

\noindent for all $1 \le j_1, j_2, \le m$.

Next, note that since $\lambda_{m+2-k} = \bar \lambda_k$,
we have $\gamma_{m+2-k} = \bar \gamma_k$,
for all $2 \le k \le m$.
Moreover, since 
$\beta_{|j_2-j_1|(k-1)} = \beta_{|j_2-j_1|((m+2-k)-1)}$,
we can write \eqref{this27} more symmetrically as
\eqref{this22} or \eqref{this23}, depending on whether
$m$ is odd or even, respectively, and in the latter case,
we also use that $\gamma_{m_0+1}$ is real, since 
$\omega^{m_0} = -1$.\CQFD
\end{Proof}

Let us again examine the cases $m=3$ and $m=4$.
In the former case, we have  

\begin{equation*}
M^{(1)} = 
\begin{pmatrix} 
  1 & -1/2 & -1/2\\
  -1/2 & 1 & -1/2\\
  -1/2 & -1/2 &1
  \end{pmatrix}.
\end{equation*}

\noindent But for $m=3$, $\beta_1 = -1/2 = \beta_2$, 
and so $M^{(1)} = M^{(2)}$.  Hence

\begin{equation}\label{Sig3}
\Sigma = \frac{2}{9} M^{(1)} + \frac{4}{9}Re(\gamma_2)M^{(2)} = \frac{2}{9}(1 + 2Re(\gamma_2))M^{(1)}.
\end{equation}

Hence, for $m=3$, cyclicity {\it always} produces a rescaled
version of the uniform iid case, with the rescaling
factor given by $1 + 2Re(\gamma_2)$.

For $m=4$, however,

\begin{equation*}
M^{(1)} = 
\begin{pmatrix} 
  1 & -1/3 & -1/3 &-1/3\\
  -1/3 & 1 & -1/3 &-1/3\\
  -1/3 & -1/3 &1  &-1/3\\
  -1/3 & -1/3 & -1/3 &1
  \end{pmatrix},
\end{equation*}

\noindent and $\beta_1 = 0$,
$\beta_2 = -1$, and $\beta_3 = 0$.  Thus,

\begin{equation*}
M^{(2)} = 
  \begin{pmatrix} 
    1 & 0 & -1 & 0\\
    0 & 1 & 0 & -1\\
    -1 & 0 & 1 & 0\\
    0 & -1 & 0 & 1
  \end{pmatrix},
\end{equation*}

\noindent and 
\begin{equation*}
M^{(3)} = 
  \begin{pmatrix} 
    1 & -1 & 1 & -1\\
    -1 & 1 & -1 & 1\\
    1 & -1 & 1 & -1\\
    -1 & 1 & -1 & 1
 \end{pmatrix}.
\end{equation*}

\noindent In this case, we have

\begin{align*}
\Sigma = \frac{3}{16} M^{(1)} + \frac{4}{16}Re(\gamma_2)M^{(2)} + \frac{2}{16}\gamma_3M^{(3)}.
\end{align*}

\noindent Next, note that
$2M^{(2)} + M^{(3)} = 3M^{(1)}$.
Then, if $Re(\gamma_2) = \gamma_3$,

\begin{align}\label{Sig4}
\Sigma &= \frac{3}{16} M^{(1)} + \frac{4}{16}Re(\gamma_2)M^{(2)} + \frac{2}{16}\gamma_3M^{(3)}\nonumber\\
&= \frac{3}{16} M^{(1)} + \frac{2}{16}(2Re(\gamma_2)M^{(1)})\nonumber\\
&= \frac{3}{16}(1 + 2Re(\gamma_2))M^{(1)},
\end{align}

\noindent so that there is still a
rescaled version of the iid case in a non-iid 
cyclic setting.  Indeed, since we know that
$\lambda_2 = a_1 + ia_2 - a_3 - ia_4 = (a_1-a_3) + i(a_2-a_4)$
and $\lambda_3 = a_1 - a_2 + a_3 - a_4$, 
we find that 

$$Re(\gamma_2) = \frac{1-a_2-2a_3-a_4}{(a_2+2a_3+a_4)^2 + (a_2-a_4)^2} - 1,$$

\noindent and $\gamma_3 = 1/(2(a_2+a_4)) - 1$.
A short calculation then shows that
$Re(\gamma_2) = \gamma_3$ if and only if
$a_3^2 = a_2 a_4$.  We thus have a
complete characterization of all $4$-letter,
cyclic Markov chains whose Young tableaux have
the same limiting shape as the uniform iid case.
In particular, choosing $a_2=a_4=a$,
for some $0 < a < 1/3$, leads to $a_3 = a$
and $a_1 = 1 - 3a$.  If, moreover,
$a = 1/4$, we have again the iid uniform case.
For $a \ne 1/4$, however, we may view the
Markov chain as a ``lazy'' version
of the uniform iid case.

Note that the scaling factor in both \eqref{Sig3}
and \eqref{Sig4} is $1 + 2Re(\gamma_2)$.
The following theorem shows that, in fact,
such a scaling factor occurs for general $m$,
and gives a spectral characterization of all
transition matrices which lead to an iid
limiting shape.

\begin{theorem}\label{thm7b}
Let $P$ be the $m \times m$ transition matrix
of an aperiodic, irreducible, cyclic Markov chain
on an $m$-letter, ordered alphabet given in
Theorem \ref{thm7}.  Then the asymptotic
covariance matrix $\Sigma$ is a rescaled version
of the iid uniform covariance matrix $\Sigma_{iidu}:= ((m-1)/m^2)M^{(1)}$ if
and only if the constants 
$\gamma_j =  \lambda_j/(1-\lambda_j)$,  
for $2 \le j \le m$, satisfy the condition

\begin{align}\label{this21b}
Re(\gamma_j) = \gamma,  \qquad \text{for all $2 \le j \le m$},
\end{align}

\noindent for some real constant $\gamma$.
Moreover, the scaling is then given by

\begin{equation}\label{this21ba}
\Sigma = (1 + 2\gamma)\Sigma_{iidu}.
\end{equation}

\end{theorem}

\noindent \begin{Proof} We first claim
that the system of matrix equations

\begin{equation}\label{this23c}
\sum_{j=2}^{m} b_j M^{(j)} = M^{(1)}
\end{equation}

\noindent has a unique solution
$b_j = 1/(m-1)$, for all $2 \le j \le m$.
Indeed, revisiting \eqref{this27}, we can
express each $M^{(j)}$ as

\begin{align}\label{this23d}
M^{(j)} 
&= \tilde M^{(j)} + \tilde M^{(-j)}\nonumber\\
&= \tilde M^{(j)} + \tilde M^{(m-j+1)},
\end{align}

\noindent where 
$(\tilde M^{(j)})_{k,\ell} = \omega^{(j-1)(\ell-k)}/2$,
for all $1 \le k,\ell \le m$, so that
\eqref{this23c} becomes

\begin{align}\label{this23e}
M^{(1)}
&= \sum_{j=2}^{m} b_j \left(\tilde M^{(j)} + \tilde M^{(m-j+1)}\right)\nonumber\\
&= \sum_{j=2}^{m} (b_j+b_{m-j+1}) \tilde M^{(j)}\nonumber\\
&= \sum_{j=2}^{m}  \tilde b_j \tilde M^{(j)},
\end{align}

\noindent where $\tilde b_j := (b_j+b_{m-j+1})/2$, 
for $2 \le j \le m$.

Now, clearly,
each $\tilde M^{(j)}$ is cyclic, so that in solving
\eqref{this23e} we need only examine
the $m$ entries in the first rows of the 
matrices.  We can thus reduce \eqref{this23e}
to the  
$m \times (m-1)$ system of equations

\begin{equation}\label{this23f}
  \begin{pmatrix}
  1 & 1 &1 &\cdots & 1\\
  \omega &\omega^2 &\omega^3 &\cdots &\omega^{m-1}\\
  \omega^2 &\omega^4 &\omega^6 &\cdots &\omega^{2(m-1)}\\
  \vdots & \vdots  &\vdots  &\ddots &\vdots\\
  \omega^{m-1} &\omega^{2(m-1)} &\omega^{3(m-1)} &\cdots &\omega^{(m-1)^2}
  \end{pmatrix}
  \begin{pmatrix}
  \tilde b_2\\
  \tilde b_3\\
  \vdots\\
  \tilde b_{m}
  \end{pmatrix}
  =
  \begin{pmatrix}
  1\\
  \frac{-1}{m-1}\\
  \frac{-1}{m-1}\\
  \vdots\\
  \frac{-1}{m-1}
  \end{pmatrix}.
\end{equation}

\noindent Since each of the last $m-1$ rows of the matrix
in \eqref{this23f} sums to $-1$, it is clear that
$\tilde b_j = 1/(m-1)$ is a solution to
the system.  To see that this solution is, 
in fact, unique, consider the 
$(m-1) \times (m-1)$ sub-matrix consisting of
the last $m-1$ rows of the matrix in \eqref{this23f},
namely,

\begin{equation}\label{this23g}
  \begin{pmatrix}
  \omega &\omega^2 &\omega^3 &\cdots &\omega^{m-1}\\
  \omega^2 &\omega^4 &\omega^6 &\cdots &\omega^{2(m-1)}\\
  \vdots & \vdots  &\vdots  &\ddots &\vdots\\
  \omega^{m-1} &\omega^{2(m-1)} &\omega^{3(m-1)} &\cdots &\omega^{(m-1)^2}
  \end{pmatrix}.
\end{equation}

\noindent Now this matrix, which is very closely
related to the Fourier matrix which arises in discrete
Fourier transform problems, is in fact invertible,
and can be shown to have one eigenvalue of $-1$,
and $m-2$ eigenvalues of the form
$\pm \sqrt{m}$ and $\pm i\sqrt{m}$, so that the
modulus of the determinant is $m^{(m-2)/2} \ne 0$.
Thus, the solution $\tilde b_j = 1/(m-1)$ is
unique, and since $b_j = (b_j+b_{m-j+1})/2 = b_{m-j+1}$,
for all $2 \le j \le m$, we conclude that
$b_j = 1/(m-1)$ as well, for all $2 \le j \le m$,
and the claim is proved.

We can now use Theorem \ref{thm7} to simplify
the asymptotic covariance matrix decomposition
as follows:

{\allowdisplaybreaks
\begin{align}\label{this23h}
\Sigma 
&= \frac{m-1}{m^2} M^{(1)} + \frac{2}{m^2}\sum_{k=2}^{m}\gamma_k M^{(k)}\nonumber\\
&= \frac{m-1}{m^2} M^{(1)} + 2\gamma \frac{1}{m^2}\sum_{k=2}^{m}M^{(k)}\nonumber\\
&= \frac{m-1}{m^2} M^{(1)} + 2\gamma \frac{m-1}{m^2}M^{(2)}\nonumber\\
&= (1+2\gamma)\frac{m-1}{m^2} M^{(1)}\nonumber\\
&= (1+2\gamma)\Sigma_{iidu},
\end{align}
}

\noindent where $\gamma = Re(\gamma_j)$, 
for all $2 \le j \le m$.  If the real parts of
$\gamma_j$ are not all identical, then the
uniqueness of the solution of \eqref{this23c}
implies that no such simplification is possible,
and the theorem is proved.\CQFD
\end{Proof}

\begin{Rem}\label{nonvac}
To see that the condition in \eqref{this21b}
is not vacuous for any $m$, recall that
for $m=4$, the ``lazy'' chain
has the iid limiting shape.  This is true
for general $m$: if $a_2 = a_3 = \cdots = a_m = a$,
for some $0 < a < 1/(m-1)$, then 
$\lambda_j = 1-(m-1)a$, for all $2 \le j \le m$.
Trivially, then, $\gamma_j = 1/((m-1)a) - 1 := \gamma$,
for all $2 \le j \le m$, so that the conditions
of Theorem \ref{thm7b} are satisfied, and the
scaling factor is given by
$1 + 2\gamma = (2 -(m-1)a)/((m-1)a)$.
Even in the $m=4$ case, however, we saw that
there were other, more general, cyclic
transition matrices which gave rise to the
iid limiting distribution.
\end{Rem}

The previous theorem indicates precisely when we
may expect the limiting shape of a cyclic Markov
chain to be identical to that of the iid uniform case.
Now the first-order behavior of all rows
of the Young tableau is $n/m + O(\sqrt{n})$ for cyclic 
Markov chains.  Although this differs from the first-order
behavior in the non-uniform iid case, one may still 
ask whether the limiting shape for a cyclic Markov chain might
still be that of some non-uniform iid case.  In fact,
this can never occur: cyclicity ensures that the 
asymptotic covariance matrix is also cyclic, and thus
cannot be equal to the asymptotic covariance matrix
of any non-uniform iid case.

Still, we may ask how to relate the iid non-uniform
limiting shape to that of a general Markov chain
having the same stationary distribution.
The following interpolation result describes the
asymptotic covariance matrix for a Markov chain whose
transition matrix is a convex combination of
an iid (uniform or non-uniform) 
transition matrix and another arbitrary
transition matrix having the same stationary distribution:

\begin{theorem}\label{thm9}
For any $m\ge 3$, let $P_0$ be the $m \times m$ transition matrix
of an irreducible, aperiodic, homogeneous Markov chain,
and let its associated asymptotic covariance matrix be
given by 

\begin{equation}\label{this33b}
\Sigma_0 = \Pi_0 + \Pi_0(S_0^{-1} D_0 S_0) + (S_0^{-1} D_0 S_0)^T \Pi_0,
\end{equation}

\noindent in the standard notations of Theorem \ref{thm2}.
Then, for $0 < \delta \le 1$, 
the transition matrix 
$P = (1-\delta)I_m + \delta P_0$
has an asymptotic covariance matrix given by

\begin{equation}\label{this34}
\Sigma = \frac{1}{\delta}\left(\Sigma_0 + (1 -\delta)\Sigma_{\Pi_0}\right),
\end{equation}

\noindent where $\Sigma_{\Pi_0}$ is the covariance matrix
associated with the iid Markov chain having the same
stationary distribution as $P_0$.

\end{theorem}

\noindent \begin{Proof} Using the standard notations of
Theorem \ref{thm2}, we will write 

$$\Sigma = \Pi + \Pi(S^{-1} D S) + (S^{-1} D S)^T \Pi$$

\noindent in terms of the decomposition $\Sigma_0$ in
\eqref{this33b}.  Now, clearly, the stationary distribution
under $P$ is that of $P_0$, so that $\Pi = \Pi_0$.  
We will thus write the stationary distribution simply as
$(\pi_1,\pi_2,\dots,\pi_m)$.
Moreover, the eigenvectors are also unchanged, 
so that $S = S_0$.  However,
for each eigenvalue $\lambda_{k,0}$ of $P_0$, we have that
$\lambda_k = (1-\delta) + \delta \lambda_{k,0}$ 
is an eigenvalue of $P$, for $1 \le k \le m$.  
Thus, for each $2 \le k \le m$, 
the diagonal entries of $D$ are given by

\begin{align*}
\gamma_k 
&:= \frac{\lambda_k}{1-\lambda_k}\nonumber\\
&= \frac{(1-\delta) + \delta \lambda_{k,0}} {\delta(1-\lambda_{k,0})}\nonumber\\
&= \frac{1-\delta}{\delta} + \gamma_{k,0},
\end{align*}

\noindent where $\gamma_{k,0}$ are the diagonal
entries of $D_0$.  We can thus decompose $D$
as follows:

{\allowdisplaybreaks
\begin{align}\label{this35}
D 
&= \text{\text{diag}}(-1/2,\gamma_2,\dots,\gamma_m)\nonumber\\
&= \text{diag}(-1/2,0,\dots,0) + \left(\frac{1-\delta}{\delta}\right) \text{diag}(0,1,\dots,1)\nonumber\\
&\qquad + \left(\frac{1}{\delta}\right) \text{diag}(0,\gamma_{2,0},\dots, \gamma_{m,0})\nonumber\\
&= \text{diag}\left(-\left(\frac{1-\delta}{2\delta}\right),0,\dots,0\right) + \left(\frac{1-\delta}{\delta}\right) I_m 
   +  \left(\frac{1}{\delta}\right)D_0.
\end{align}
}

Next, recall from Proposition \ref{prop1} that the 
first column of $S^{-1}$ is $(1,1,\dots,1)^T$.  Hence,

{\allowdisplaybreaks\begin{align}\label{this36}
S^{-1}DS 
&= S_0^{-1}DS_0\nonumber\\
&= \begin{pmatrix}
      1 &* &\cdots &*\\
      \vdots &\vdots &\cdots &\vdots\\
      \vdots &\vdots &\cdots &\vdots\\
      1 &* &\cdots &*
   \end{pmatrix}
   \begin{pmatrix}
      -\frac{1-\delta}{2\delta} &0 &\cdots &0\\
      0 &\ddots &\cdots &\vdots\\
      \vdots &\vdots &\ddots &\vdots\\
      0 &\cdots &\cdots &0
   \end{pmatrix}
   \begin{pmatrix}
      \pi_1 &\pi_2 &\cdots &\pi_m\\
      * &* &\cdots &*\\
      \vdots &\vdots &\cdots &\vdots\\
      * &* &\cdots &*
   \end{pmatrix}\nonumber\\
&\qquad \qquad + \left(\frac{1-\delta}{\delta}\right) S_0^{-1}I_m S_0
   +  \left(\frac{1}{\delta}\right) S_0^{-1}D_0 S_0\nonumber\\
&=-\left(\frac{1-\delta}{2\delta}\right)
   \begin{pmatrix}
      \pi_1 &\pi_2 &\cdots &\pi_m\\
      \pi_1 &\pi_2 &\cdots &\pi_m\\
      \vdots &\vdots &\cdots &\vdots\\
      \pi_1 &\pi_2 &\cdots &\pi_m
   \end{pmatrix}
   + \left(\frac{1-\delta}{\delta}\right) I_m
   +  \left(\frac{1}{\delta}\right) S_0^{-1}D_0 S_0,
\end{align}
}

\noindent  which gives us

{\allowdisplaybreaks
\begin{align}\label{this37}
\Pi S^{-1}DS 
&= \Pi_0 S^{-1}DS\nonumber\\
&=-\left(\frac{1-\delta}{2\delta}\right)
   \begin{pmatrix}
      \pi_1 &0 &\cdots &0\\
      0 &\pi_2 &\cdots &\vdots\\
      0 &\vdots &\ddots &\vdots\\
      0 &\cdots &\cdots &\pi_m
   \end{pmatrix}
   \begin{pmatrix}
      \pi_1 &\pi_2 &\cdots &\pi_m\\
      \pi_1 &\pi_2 &\cdots &\pi_m\\
      \vdots &\vdots &\cdots &\vdots\\
      \pi_1 &\pi_2 &\cdots &\pi_m
   \end{pmatrix}\nonumber\\
&\qquad \qquad   + \left(\frac{1-\delta}{\delta}\right) \Pi_0
   +  \left(\frac{1}{\delta}\right) \Pi_0 S_0^{-1}D_0 S_0\nonumber\\
&=-\left(\frac{1-\delta}{2\delta}\right) 
   \begin{pmatrix}
      \pi_1^2 &\pi_1 \pi_2 &\cdots &\pi_1 \pi_m\\
      \pi_2 \pi_1 &\pi_2^2 &\cdots &\pi_2 \pi_m\\
      \vdots &\vdots &\ddots &\vdots\\
      \pi_m \pi_1 &\pi_m \pi_2 &\cdots &\pi_m^2
   \end{pmatrix}\nonumber\\
&\qquad \qquad  + \left(\frac{1-\delta}{\delta}\right) \Pi_0
   +  \left(\frac{1}{\delta}\right) \Pi_0 S_0^{-1}D_0 S_0.
\end{align}
}

\noindent  Finally, we can express $\Sigma$ as

{\allowdisplaybreaks
\begin{align}\label{this38}
\Sigma 
&= \Pi + \Pi(S^{-1} D S) + (S^{-1} D S)^T \Pi\nonumber\\
&= \left(\frac{1}{\delta}\right) \Pi_0 
   + \left(1-\frac{1}{\delta}\right) \Pi_0 
+ \Pi_0(S^{-1} D S) + (\Pi_0(S^{-1} D S))^T \nonumber\\
&=\left(\frac{1}{\delta}\right) \Sigma_0
  +\left(1 - \frac{1}{\delta} \right)(\Pi_0 - 2\Pi_0)\nonumber\\
&\qquad \qquad   +\left(1 - \frac{1}{\delta} \right)
   \begin{pmatrix}
      \pi_1^2 &\pi_1 \pi_2 &\cdots &\pi_1 \pi_m\\
      \pi_2 \pi_1 &\pi_2^2 &\cdots &\pi_2 \pi_m\\
      \vdots &\vdots &\ddots &\vdots\\
      \pi_m \pi_1 &\pi_m \pi_2 &\cdots &\pi_m^2
   \end{pmatrix}\nonumber\\
&=\left(\frac{1}{\delta}\right) \Sigma_0
  +\left(1 - \frac{1}{\delta} \right) (-\Sigma_{\Pi_0})\nonumber\\
&= \frac{1}{\delta}\left(\Sigma_0 + (1 -\delta)\Sigma_{\Pi_0}\right),
\end{align}
}

\noindent and we are done.\CQFD
\end{Proof}

Thus far we have expressed our limiting laws in terms
of Brownian functionals whose Brownian motions have a
non-trivial covariance structure arising directly from
the specific nature of the transition matrix.  It is of
interest to instead express the limiting laws in
terms of {\it standard} Brownian motions.

Since the asymptotic covariance matrix $\Sigma$ 
is non-negative definite, we can find an $m \times m$
matrix $C$ such that $\Sigma = CC^T$. 
(The matrix $C$ is not unique, since
$(CQ)(CQ)^T = CC^T = \Sigma$ for any
orthogonal matrix $Q$.)  Clearly, we then have

\begin{equation}\label{this40}
(\sigma_1 \tilde B^1(t),\sigma_2 \tilde B^2(t),\dots,\sigma_m \tilde B^m(t))^T
= C(B^1(t),B^2(t),\dots,B^m(t))^T,
\end{equation}

\noindent where $(B^1(t),B^2(t),\dots,B^m(t))^T$ is a standard,
$m$-dimensional Brownian motion, since

\begin{align*}
&\bbe \bigl[(\sigma_1 \tilde B^1(t),\sigma_2 \tilde B^2(t),\dots,\sigma_m \tilde B^m(t))^T
      (\sigma_1 \tilde B^1(t),\sigma_2 \tilde B^2(t),\dots,\sigma_m \tilde B^m(t)) \bigr]\nonumber\\
&\qquad =\bbe  \bigl[C(B^1(t),B^2(t),\dots,B^m(t))^T\bigr]\bigl[(C(B^1(t),B^2(t),\dots,B^m(t))^T\bigr]^T\nonumber\\
&\qquad =C \bigl[ \bbe(B^1(t),B^2(t),\dots,B^m(t))^T)(B^1(t),B^2(t),\dots,B^m(t)) \bigr]C^T\nonumber\\
&\qquad =C (t I_m) C^T\nonumber\\
&\qquad = t\Sigma.
\end{align*}

Next, we can, without loss of generality, assume that $\tau(\ell) = \ell$,
for all $\ell$, and so write our main result \eqref{item7r} 
in Theorem \ref{thm4} as

\begin{align}\label{this41}
\frac{V^r_n -\nu_r n}{\sqrt{n}} &\Rightarrow  \sum_{k=1}^{m_r} \sigma_{k} \tilde B^{k}(1) 
+ \max_{I_{r-m_r,d_r}} 
    \sum_{j=1}^{r-m_r} \sum_{\ell=j}^{(d_r + m_r - r+j)} \sigma_{m_r + \ell} \tilde B^{m_r + \ell}(\Delta t_{j,\ell})\nonumber\\
&:= V^r_{\infty}.
\end{align}

\noindent Simply substituting \eqref{this40} into
\eqref{this41} immediately yields

\begin{align}\label{this42}
V^r_{\infty}
&=\sum_{k=1}^{m_r} \left(\sum_{i=1}^m C_{k,i} B^{i}(1) \right)\nonumber\\
&\qquad  + \max_{I_{r-m_r,d_r}} 
             \sum_{j=1}^{r-m_r} \sum_{\ell=j}^{(d_r + m_r - r+j)} 
               \left(\sum_{i=1}^m C_{m_r + \ell,i} B^{i}(\Delta t_{j,\ell})\right)\nonumber\\
&=\sum_{i=1}^m \sum_{k=1}^{m_r} C_{k,i} B^{i}(1) \nonumber\\
&\qquad  + \max_{I_{r-m_r,d_r}} 
             \sum_{i=1}^m  \sum_{j=1}^{r-m_r} \sum_{\ell=j}^{(d_r + m_r - r+j)} 
                C_{m_r + \ell,i} B^{i}(\Delta t_{j,\ell}).
\end{align}

Now the first term in \eqref{this42} is simply a Gaussian term
whose variance can be computed explicitly.  Unfortunately,
the maximal term does not in general succumb to any
significant simplifications.  However, in the 
iid case, we can further simplify \eqref{this42} in a
very satisfying way.

  Indeed, since, in the iid case, we have
$\sigma_k^2 = \pi_k(1-\pi_k)$ and, for $k \ne \ell$,
$\sigma_{k,\ell} = - \pi_k \pi_{\ell}$, one can quickly
check that $C$ can be chosen so that
$C_{k,k} = \sqrt{\pi_k} -\sqrt{\pi_k}\pi_k$, and, for $k \ne \ell$,
$C_{k,\ell} = -\sqrt{\pi_{\ell}}\pi_k$.  
Moreover, for all $m_r + 1 \le k \le m_r + d_r$,
$\pi_k$ = $\pi_{m_r + 1} = \pi_r$.
Then, within the
maximal term, 
$C_{m_r + \ell,i} = \sqrt{\pi_r} - \pi_r\sqrt{\pi_r}$,
for $i = m_r + \ell$, and
$C_{m_r + \ell,i} = - \pi_{r}\sqrt{\pi_i}$,
for $i \ne m_r + \ell$.
With the convention that $\nu_0 = 0$, 
we can then express \eqref{this42} as 

{\allowdisplaybreaks
\begin{align}\label{this43}
V^r_{\infty}
&=\sum_{i=1}^{m_r}  \sqrt{\pi_i} B^{i}(1) 
  + \sum_{i=1}^m \sum_{k=1}^{m_r} (-\sqrt{\pi_{i}}\pi_k) B^{i}(1) \nonumber\\
&\qquad  +  \max_{I_{r-m_r,d_r}} \biggl\{
                 \sum_{j=1}^{r-m_r} \sum_{\ell=j}^{(d_r + m_r - r+j)} 
                        \sqrt{\pi_r} B^{m_r+\ell}(\Delta t_{j,\ell})\nonumber\\
&\qquad \qquad +  \sum_{i=1}^m  \sum_{j=1}^{r-m_r} \sum_{\ell=j}^{(d_r + m_r - r+j)} 
                (- \pi_{r}\sqrt{\pi_i}) B^{i}(\Delta t_{j,\ell})\biggr\}\nonumber\\
&=\sum_{i=1}^{m_r}  \sqrt{\pi_i} B^{i}(1) 
  - \sum_{i=1}^m \sqrt{\pi_{i}} B^{i}(1)  \sum_{k=1}^{m_r}  \pi_k \nonumber\\
&\qquad  +  \sqrt{\pi_r} \max_{I_{r-m_r,d_r}} \biggl\{
                 \sum_{j=1}^{r-m_r} \sum_{\ell=j}^{(d_r + m_r - r+j)} 
                        B^{m_r+\ell}(\Delta t_{j,\ell})\nonumber\\
&\qquad \qquad -  \sqrt{\pi_{r}}\sum_{i=1}^m \sqrt{\pi_{i}} \sum_{j=1}^{r-m_r} \sum_{\ell=j}^{(d_r + m_r - r+j)} 
                 B^{i}(\Delta t_{j,\ell})\biggr\}\nonumber\\
&=\biggl\{\sum_{i=1}^{m_r}  \sqrt{\pi_i} B^{i}(1) 
  - \nu_{m_r} \sum_{i=1}^m \sqrt{\pi_{i}} B^{i}(1) 
  - \pi_r r\sum_{i=1}^m \sqrt{\pi_{i}} B^{i}(1)\biggr\}\nonumber\\
&\qquad\qquad   + \sqrt{\pi_r} \max_{I_{r-m_r,d_r}} 
                 \sum_{j=1}^{r-m_r} \sum_{\ell=j}^{(d_r + m_r - r+j)} 
                        B^{m_r+\ell}(\Delta t_{j,\ell})\nonumber\\
&=\biggl\{\sum_{i=1}^{m_r}  \sqrt{\pi_i} B^{i}(1) 
  - (\nu_{m_r} + \pi_r r) \sum_{i=1}^m \sqrt{\pi_{i}}  B^{i}(1)\biggr\}\nonumber\\
&\qquad  + \sqrt{\pi_r} \max_{I_{r-m_r,d_r}} 
                 \sum_{j=1}^{r-m_r} \sum_{\ell=j}^{(d_r + m_r - r+j)} 
                        B^{m_r+\ell}(\Delta t_{j,\ell})\nonumber\\
&=\Bigl\{(1-\nu_{m_r} - \pi_{r}r)   \sum_{i=1}^{m_r}  \sqrt{\pi_i}B^{i}(1)\nonumber\\
&\qquad\qquad  - (\nu_{m_r} + \pi_r r) \sum_{i=m_r+d_r+1}^m \sqrt{\pi_{i}} B^{i}(1)\Bigr\}\nonumber\\
&\qquad  + \sqrt{\pi_r}\Bigl\{- (\nu_{m_r} + \pi_r r)\sum_{i=m_r+1}^{m_r+d_r} B^{i}(1)\nonumber\\
&\qquad\qquad    + \max_{I_{r-m_r,d_r}} 
                 \sum_{j=1}^{r-m_r} \sum_{\ell=j}^{(d_r + m_r - r+j)} 
                        B^{m_r+\ell}(\Delta t_{j,\ell})\Bigr\}.
\end{align}
}

Note that the first two Gaussian term of \eqref{this43}
are independent of the remaining two Gaussian-maximal expression
terms.

Following Glynn and Whitt\cite{GW} and Barishnykov\cite{Ba},
who studied the Brownian functional

$$D_m = \max_{I_{1,m}} \sum_{\ell=1}^{m} B^{\ell}(\Delta t_{\ell}),$$

\noindent we define the following, more general, Brownian functional:

\begin{equation}\label{this44}
D_{r,m} := \max_{I_{r,m}} 
                 \sum_{j=1}^{r} \sum_{\ell=j}^{(m-r+j)} 
                        B^{\ell}(\Delta t_{j,\ell}),
\end{equation}

\noindent where $1 \le r \le m$.  Clearly, the maximal term in
\eqref{this43} has just such a form.  We also remark that 
$D_{r,m}$ corresponds to the sum of the $r$ largest eigenvalues
of an $m \times m$ GUE matrix.

To better understand \eqref{this43}, we may, without much
loss in generality, focus on the first block, that is,
values of $r$ such that $m_r = 0$.  The first  Gaussian
term of \eqref{this43} thus vanishes, and, writing
$\pi_{max}$ for $\pi_r$, we have

{\allowdisplaybreaks
\begin{align}\label{this45}
V^r_{\infty}
&= -r\pi_{max}  \sum_{i=d_1+1}^m \sqrt{\pi_{i}} B^{i}(1) \nonumber\\
& \qquad + \sqrt{\pi_{max}} 
     \left(-r \pi_{max} \sum_{i=1}^{d_1} B^{i}(1) + D_{r,d_r}\right).
\end{align}
}

\noindent  In the uniform iid case, the first Gaussian term of \eqref{this45}
itself vanishes, since $d_r = d_1 = m$, and we have

\begin{align}\label{this46}
V^r_{\infty}
&= \frac{1}{\sqrt{m}}
     \left(-\frac{r}{m}\sum_{i=1}^{m}  B^{i}(1) + D_{r,m}\right)\nonumber\\
&:= \frac{H_{r,m}}{\sqrt{m}}.
\end{align}

\noindent For $r=1$, this result corresponds to Theorem $4.1$ of
the authors' previous paper \cite{HL}.  Furthermore, specializing
\eqref{this45} to $r=1$, 

{\allowdisplaybreaks
\begin{align}\label{this47}
\frac{LI_n - \pi_{max}n}{\sqrt{n}} 
&\Rightarrow -\pi_{max}  \sum_{i=d_1+1}^m \sqrt{\pi_{i}} B^{i}(1) \nonumber\\
& \qquad + \sqrt{\pi_{max}} 
     \left(-\pi_{max} \sum_{i=1}^{d_1} B^{i}(1) + D_{1,d_1}\right)\nonumber\\
&= -\pi_{max}  \sum_{i=d_1+1}^m \sqrt{\pi_{i}} B^{i}(1) \nonumber\\
& \qquad + \sqrt{\pi_{max}} 
       \left(\frac{1}{d_1}-\pi_{max}\right) \sum_{i=1}^{d_1}B^{i}(1)\nonumber\\
& \qquad + \sqrt{\pi_{max}}H_{1,d_1}.
\end{align}
}

One can easily compute the variance of the Gaussian terms
in \eqref{this47} to be $\pi_{max}(1-d_1\pi_{max})/d_1$,
which is consistent with Proposition $4.1$ of
the authors' previous paper \cite{HL}.

The iid development above suggests that we
can find additional cases which yield simple functionals
of standard Brownian motions.  Indeed, the first property of
the matrix $C$ in the iid case that allowed
the functionals to be simplified was that $C_{k,\ell} = c_{\ell}$,
for all $k \ne \ell, m_r+1 \le k \le m_r+d_r$, and $ 1 \le \ell \le m$, where
$c_1, c_2, \dots, c_m$ were real numbers.  
Then, writing the diagonal terms of $C$ as 
$C_{k,k} = b_k + c_k$, for $m_r+1 \le k \le m_r+d_r$,
we may revisit \eqref{this42}, and write

{\allowdisplaybreaks
\begin{align}\label{this49}
V^r_{\infty}
&=\sum_{i=1}^m \sum_{k=1}^{m_r} C_{k,i} B^{i}(1) 
  + \max_{I_{r-m_r,d_r}} 
             \sum_{i=1}^m  \sum_{j=1}^{r-m_r} \sum_{\ell=j}^{(d_r + m_r - r+j)} 
                C_{m_r + \ell,i} B^{i}(\Delta t_{j,\ell})\nonumber\\
&=\sum_{i=1}^m \sum_{k=1}^{m_r} C_{k,i} B^{i}(1) 
 +  \max_{I_{r-m_r,d_r}} \biggl\{
                 \sum_{j=1}^{r-m_r} \sum_{\ell=j}^{(d_r + m_r - r+j)} 
                        b_{m_r+\ell} B^{m_r+\ell}(\Delta t_{j,\ell})\nonumber\\
&\qquad \qquad +  \sum_{i=1}^m  \sum_{j=1}^{r-m_r} \sum_{\ell=j}^{(d_r + m_r - r+j)} 
                c_i B^{i}(\Delta t_{j,\ell})\biggr\}\nonumber\\
&=\sum_{i=1}^m \sum_{k=1}^{m_r} C_{k,i} B^{i}(1)  
  +  r  \sum_{i=1}^m  c_i B^{i}(1)\nonumber\\
&\qquad  +  \max_{I_{r-m_r,d_r}} 
                 \sum_{j=1}^{r-m_r} \sum_{\ell=j}^{(d_r + m_r - r+j)} 
                        b_{m_r+\ell} B^{m_r+\ell}(\Delta t_{j,\ell})
\end{align}
}

Except for the fact that we have written the functional in terms of
standard Brownian motions, the maximal term in \eqref{this49} 
is no simpler than that of our original functional.  
However, the second property of the iid case 
that yielded further simplifications
was that $b_k = b$, for all $m_r+1 \le k \le m_r + d_r$.  In this
case, \eqref{this49} becomes

{\allowdisplaybreaks
\begin{align}\label{this50}
V^r_{\infty}
&=\sum_{i=1}^m \sum_{k=1}^{m_r} C_{k,i} B^{i}(1)  
  +  r  \sum_{i=1}^m  c_i B^{i}(1)\nonumber\\
&\qquad  +  b\max_{I_{r-m_r,d_r}} 
                 \sum_{j=1}^{r-m_r} \sum_{\ell=j}^{(d_r + m_r - r+j)} 
                        B^{m_r+\ell}(\Delta t_{j,\ell})
\end{align}
}

Again, by focusing on the first block, we no longer have the initial
Gaussian term, and \eqref{this50} becomes

{\allowdisplaybreaks
\begin{align}\label{this51}
V^r_{\infty}
&= r\sum_{i=d_1+1}^{m}  c_i B^{i}(1) \nonumber\\
&\qquad + r\sum_{i=1}^{d_1}  c_i B^{i}(1) 
     +  b\max_{I_{r,d_1}} 
            \sum_{j=1}^{r} \sum_{\ell=j}^{(d_1- r+j)} 
                      B^{\ell}(\Delta t_{j,\ell})\nonumber\\
&= r\sum_{i=d_1+1}^{m}  c_i B^{i}(1) 
  + \left(r\sum_{i=1}^{d_1}  c_i B^{i}(1) 
     +  b D_{r,d_1}\right) \nonumber\\
&= r\sum_{i=d_1+1}^{m}  c_i B^{i}(1) + r\sum_{i=1}^{d_1}  \left(c_i + \frac{b}{d_1}\right) B^{i}(1) \nonumber\\
&\qquad +  b \left(-\frac{r}{d_1}\sum_{i=1}^{d_1} B^{i}(1) + D_{r,d_1}\right) \nonumber\\
&= r\sum_{i=d_1+1}^{m}  c_i B^{i}(1) + r\sum_{i=1}^{d_1}  \left(c_i + \frac{b}{d_1}\right) B^{i}(1) 
 +  b H_{r,d_1}.
\end{align}
}

We restate these results in the following theorem:

\begin{theorem}\label{thmindbm} Assume, without
loss of generality, that $\tau(\ell) = \ell$,
for all $1 \le \ell \le m$, in the
notations of Theorem \ref{thm4}.
Moreover, let the asymptotic covariance matrix 
be given by $\Sigma = CC^T$, where $C$ is an
$m \times m$ matrix whose first $d_1$ rows
are given by

\begin{equation}\label{this51a}
	\begin{cases}
        C_{k,\ell} = c_{\ell}, &k \ne \ell, 1 \le k \le d_1, 1 \le \ell \le m\\
        C_{k,k}    = b + c_k,  &1 \le k \le d_1,
	\end{cases}
\end{equation}

\noindent for some real constants $c_1, c_2, \dots, c_m$ and $b$.
Then, for $1 \le r \le d_1$,

\begin{equation}\label{this51b}
V^r_{\infty}
= r\sum_{i=d_1+1}^{m}  c_i B^{i}(1) + r\sum_{i=1}^{d_1}  \left(c_i + \frac{b}{d_1}\right) B^{i}(1) 
 +  b H_{r,d_1},
\end{equation}

\noindent where $H_{r,d_1}$ is the maximal functional

\begin{equation*}
H_{r,d_1} := \frac{1}{\sqrt{d_1}}
     \left(-\frac{r}{d_1}\sum_{i=1}^{d_1}  B^{i}(1) 
    + \max_{I_{r,d_1}} 
                 \sum_{j=1}^{r} \sum_{\ell=j}^{(d_1-r+j)} 
                        B^{\ell}(\Delta t_{j,\ell})\right).
\end{equation*}

\end{theorem}

\begin{Rem}\label{remindbm}
One can generalize Theorem \ref{thmindbm} to non-initial
blocks ({\it i.e.}, to $r > d_1$) by 
extending the conditions in \eqref{this51a}
to non-initial blocks and then
applying the theorem to
$V^r_{\infty} - V^{m_r}_{\infty}$.
\end{Rem}

To better understand which asymptotic covariance matrices
$\Sigma$ can be decomposed in this manner, 
the conditions $C_{k,\ell} = c_{\ell}$,
for all $k \ne \ell, 1 \le k \le d_1,1 \le \ell \le m$,
and $b_k = b$, for all $1 \le k \le d_1$,
imply that

\begin{equation}\label{this52}
\sigma_k^2 = b^2 + 2bc_k + \sum_{i=1}^m c_i^2,
\end{equation}

\noindent for $1 \le k \le d_1$, and

\begin{equation}\label{this53}
\sigma_{k,\ell} = bc_k + bc_{\ell} + \sum_{i=1}^m c_i^2,
\end{equation}

\noindent for $1 \le k < \ell \le d_1$.  

If we let $(Z_1,Z_2,\dots,Z_m)$ be a centered
Gaussian random vector with covariance matrix 
$\Sigma$, then \eqref{this52} and \eqref{this53}
give us

\begin{align}\label{this54}
\bbe (Z_k - Z_{\ell})^2 
&= \sigma_k^2 - 2\sigma_{k,\ell} + \sigma_{\ell}^2\nonumber\\
&= 2b^2,
\end{align}

\noindent for all $1 \le k < \ell \le d_1$.  That is,
the $L^2$-distance between any pair
$(Z_k,Z_{\ell})$ is the same,
for $1 \le k < \ell \le d_1$.  

Notice that if
$\sigma_k^2 = \sigma^2$, for all $1 \le k \le d_1$,
then in fact \eqref{this54} implies that
$\rho_{k,\ell} = \sigma_{k,\ell}/\sigma_k \sigma_{\ell} = 1-b^2/\sigma^2 := \rho$,
for all $1 \le k < \ell \le d_1$.  That is, 
the $d_1 \times d_1$ submatrix of $\Sigma$
must be permutation-symmetric.

Next, we note that, for
$1 \le k <\ell \le d_1$,

\begin{equation}\label{this55}
\sigma_k^2 - \sigma_{\ell}^2  = 2b(c_k - c_{\ell}),
\end{equation}

\noindent so that $c_k = \sigma_k^2/(2b) + c_0$,
for some constant $c_0$.  Substituting this
expression into \eqref{this52} and,
writing $\Gamma = \sum_{i=d_1+1}^mc_i^2$,
we obtain

{\allowdisplaybreaks
\begin{align}\label{this56}
\sigma_k^2
&= b^2 + 2b\left(\frac{\sigma_k^2}{2b} + c_0\right) 
  + \sum_{i=1}^{d_1}\left(\frac{\sigma_i^2}{2b} + c_0\right)^2 + \Gamma\nonumber\\
&= b^2 + \sigma_k^2  + 2bc_0
  + \sum_{i=1}^{d_1}\left(\frac{\sigma_i^2}{2b} + c_0\right)^2 + \Gamma.
\end{align}
}

\noindent Writing $\overline{\sigma^r} = (\sum_{i=1}^{d_1}\sigma_i^r)/d_1$,
for any $r>0$, \eqref{this56} gives us

{\allowdisplaybreaks
\begin{align}\label{this57}
&b^2 + 2bc_0
  + \sum_{i=1}^{d_1}\left(\frac{\sigma_i^2}{2b} + c_0\right)^2 + \Gamma\nonumber\\
&\qquad = d_1 c_0^2 
        + \left(2b+\frac{d_1}{b}\overline{\sigma^2}\right)c_0 
        + \left(b^2 + \frac{d_1\overline{\sigma^4}}{4b^2} + \Gamma\right)\nonumber\\
&\qquad = 0.
\end{align}
}

In order for $c_0$ to be a real number, 
the discriminant of the quadratic equation in \eqref{this57}
must satisfy

\begin{align}\label{this58}
\left(2b+\frac{d_1}{b}\overline{\sigma^2}\right)^2
  - 4d_1\left(b^2 + \frac{d_1\overline{\sigma^4}}{4b^2} + \Gamma\right)
  \ge 0,
\end{align}

\noindent which leads to the inequality

\begin{align}\label{this58c}
(d_1-1)b^4 - d_1(\overline{\sigma^2} - \Gamma)
                +\frac{d_1^2}{4}\left(\overline{\sigma^4} - \left(\overline{\sigma^2}\right)^2\right) \le 0.
\end{align}

This inequality, in turn, gives 
us constraints on $b^2$.  Indeed, the necessary 
and sufficient condition needed 
for such a $b^2$ to exist is given 
by examining the quadratic in $b^2$ in \eqref{this58c}
at its extremal point, namely,
at $b^2 = d_1(\overline{\sigma^2} - \Gamma))/(2(d_1-1))$.  
Doing so leads to the condition

\begin{align}\label{this58d}
-\left(\frac{d_1^2(\overline{\sigma^2} - \Gamma)^2}{4(d_1-1)}\right) 
  + \frac{d_1^2}{4}\left(\overline{\sigma^4} - \left(\overline{\sigma^2}\right)^2\right) \le 0,
\end{align}

\noindent or simply,

\begin{align}\label{this58e}
 \overline{\sigma^4}  - \left(\overline{\sigma^2}\right)^2
  \le \left(\frac{(\overline{\sigma^2} - \Gamma)^2}{d_1-1}\right),
\end{align}

\noindent since 
$\overline{\sigma^4} \ge  \left(\overline{\sigma^2}\right)^2$.  
The closer that $\Gamma$ is
to $\overline{\sigma^2}$, the more similar that the
$d_1$ variances must be.  Thus, \eqref{this58e} functions
as a bound on the variability among these $d_1$ variances.
Provided that the variances satisfy \eqref{this58e},
the condition on $b^2$ is given by

\begin{align}\label{this58ea}
b^2 \in \Biggl( \frac{d_1}{2(d_1-1)}\Biggl\{ (\overline{\sigma^2} - \Gamma)
    -\sqrt{(\overline{\sigma^2} - \Gamma)^2 
               - (d_1-1) \left(\overline{\sigma^4} - \left(\overline{\sigma^2}\right)^2\right) }  \Biggr\},\nonumber\\
\frac{d_1}{2(d_1-1)}\Biggl\{ (\overline{\sigma^2} - \Gamma)
    + \sqrt{(\overline{\sigma^2} - \Gamma)^2 
               - (d_1-1) \left(\overline{\sigma^4} - \left(\overline{\sigma^2}\right)^2\right) }  \Biggr\} \Biggr) .
\end{align}

Now consider the doubly stochastic case, where $d_1 = m$.
Applying the general fact that each row of $\Sigma$ must
necessarily sum to zero, we use \eqref{this52} and \eqref{this53}
to find that, for each $1 \le k \le m$,

\begin{align}\label{this58eb}
\sum_{\ell=1}^{m}\sigma_{k,\ell} 
&=\sigma_k^2 + \sum_{\ell \ne k}\sigma_{k,\ell} \nonumber\\
&= \left(b^2 + 2bc_k + \sum_{i=1}^m c_i^2\right)
  + \sum_{\ell \ne k} \left(bc_k + bc_{\ell} + \sum_{i=1}^m c_i^2\right)\nonumber\\
&= b^2 + b\left(m c_k + \sum_{\ell=1}^{m}c_{\ell}\right) + m \sum_{i=1}^m c_i^2\nonumber\\
&=0,
\end{align}

\noindent  so that
$c_k = c \in \bbr$, for all $1 \le k \le m$. 
Substituting $c$ back into \eqref{this58eb}
gives us

\begin{align}\label{this58ec}
\sum_{\ell=1}^{m}\sigma_{k,\ell} 
&= b^2 + b(mc + mc) + m(mc^2)\nonumber\\
&= (b+mc)^2 = 0,
\end{align}

\noindent so that $b = -mc$.  This then implies
that $\sigma_k^2 = m(m-1)c^2$ and 
$\sigma_{k,\ell} = -mc^2$, for all $1 \le k \le m$,
$\ell \ne k$.  But this is precisely a permutation-symmetric
covariance matrix, 
which in the iid case
corresponds to the class of Markov chains having
a uniform stationary distribution.

We summarize these results in the following:

\begin{theorem}\label{thmcondbm}
In order that the asymptotic covariance matrix
$\Sigma$ have a decomposition $\Sigma = CC^T$,
where 

\begin{equation}\label{this58g}
	\begin{cases}
        C_{k,\ell} = c_{\ell}, &k \ne \ell, \quad 1 \le k \le d_1, 
                               \quad 1 \le \ell \le m,\\
        C_{k,k}    = b + c_k,  &1 \le k \le d_1,
	\end{cases}
\end{equation}

\noindent for some real constants $c_1, c_2, \dots, c_m$ and $b$,
it is necessary and sufficient that

\begin{align}\label{this58h}
 \overline{\sigma^4}  - \left(\overline{\sigma^2}\right)^2
  \le \frac{|\overline{\sigma^2} - \Gamma|}{\sqrt{d_1-1}},
\end{align}

\noindent where $\Gamma = \sum_{i=d_1+1}^mc_i^2$, and
$\overline{\sigma^r} = (\sum_{i=1}^{d_1}\sigma_i^r)/d_1$,
for any $r>0$.  In this case,

\begin{align}\label{this58i}
b^2 \in \Biggl( \frac{d_1}{2(d_1-1)}\Biggl\{ (\overline{\sigma^2} - \Gamma)
    -\sqrt{(\overline{\sigma^2} - \Gamma)^2 
               - (d_1-1) \left(\overline{\sigma^4} - \left(\overline{\sigma^2}\right)^2\right) }  \Biggr\},\nonumber\\
\frac{d_1}{2(d_1-1)}\Biggl\{ (\overline{\sigma^2} - \Gamma)
    + \sqrt{(\overline{\sigma^2} - \Gamma)^2 
               - (d_1-1) \left(\overline{\sigma^4} - \left(\overline{\sigma^2}\right)^2\right) }  \Biggr\} \Biggr).
\end{align}

\noindent In particular, if $d_1 = m$, the asymptotic
covariance matrix must be permutation-symmetric, with
$c_k = c$, for all $k$, and $b = -mc$, so that
the common variance is $m(m-1)c^2$ and the common
covariances are all $-mc^2$.
\end{theorem}



\section{Connections to Random Matrix Theory}

For iid uniform $m$-letter alphabets, 
the limiting law of the Young tableau 
corresponds to the joint distribution 
of the eigenvalues of an $m \times m$ 
matrix from the traceless GUE \cite{Jo}. 
In the non-uniform iid case, we further noted that
Its, Tracy, and Widom \cite{ITW1,ITW2} have essentially
described the limiting shape as that of the 
joint distribution of the eigenvalues 
of a random matrix consisting of independent diagonal
blocks, each of which is a matrix
from the GUE.  The size
of each block depends
upon the multiplicity
of the corresponding stationary probability.
In addition, there is a zero-trace
condition involving the stationary probabilities
on the composite matrix.

As a first step in extending these
connections between Brownian functionals and
spectra of random matrices, recall the general case when 
the stationary probabilities are all distinct
(see Remark \ref{RMTconnect}).
Our Brownian functionals then have
no true maximal terms, so that the limiting
shape, $(R_{\infty}^1,R_{\infty}^2,\dots,R_{\infty}^m)$
is simply multivariate normal, with
covariance matrix $\Sigma$ (or, more precisely,
the matrix obtained by permuting the rows and columns 
of $\Sigma$ using $\tau$, the 
permutation of $\{1,2,\dots,m\}$
previously defined).  Trivially,
this limiting law corresponds to the spectrum of a diagonal
matrix whose elements are multivariate normal
with the same covariance matrix $\Sigma$.

We can see that this general result is consistent
with the non-uniform iid case having distinct
probabilities.  Indeed, each block is of size $1$,
and is rescaled so that the variance is 
$\pi_{\tau(i)}(1-\pi_{\tau(i)})$, for
$1 \le i \le m$.  
Because of this rescaling, instead of having
a generalized zero-trace condition, as in the
non-rescaled matrices used in \cite{ITW1,ITW2},
our condition is rather a true zero-trace
condition.  This zero-trace condition is clear,
since the covariance matrix for {\it any} iid
case (uniform and non-uniform alike) is that of a 
multinomial distribution with parameters
$(n=1;\pi_{\tau(1)},\pi_{\tau(2)},\dots,\pi_{\tau(m)})$,
and any $(Y_1,Y_2,\dots,Y_m)$ having such a 
distribution of course satisfies
$\sum_{i=1}^m Y_i = 1$, so that
$Var(\sum_{i=1}^m Y_i) = 0$, which implies 
the zero-trace condition for
 $(R_{\infty}^1,R_{\infty}^2,\dots,R_{\infty}^m)$.

Next, consider the case when each stationary 
probability has multiplicity no greater than $2$.  
We conjecture that the limiting shape
 $(R_{\infty}^1,R_{\infty}^2,\dots,R_{\infty}^m)$
is that of the spectrum of
a direct sum of certain $1 \times 1$
and/or $2 \times 2$ random matrices.
 Specifically,
let $\kappa \le m$ be the number of distinct probabilities
among the stationary distributions. Then the
composite matrix consists of a direct sum
of $\kappa$ GUE matrices which are as follows.
First, the overall diagonal 
$(X_1, X_2, \dots, X_m)$ of the matrix has
a $N(0,\Sigma)$ distribution.
Next, if $d_r = 1$, then the GUE matrix is
simply the $1 \times 1$ matrix $(X_r)$.
Finally, if $d_r=2$, then the GUE matrix
is the $2 \times 2$ matrix

\begin{equation*}
  \begin{pmatrix}
     X_{m_r+1} &Y_{m_r+1} + iZ_{m_r+1}\\
     Y_{m_r+1} - iZ_{m_r+1}   &X_{m_r+2}
  \end{pmatrix},
\end{equation*}

\noindent whose off-diagonal random variables 
$Y_{m_r+1}$ and $Z_{m_r+1}$
are iid, centered, normal random variables,
independent of all other random variables in the
overall matrix, 
with variance

\begin{equation*}
(\sigma^2_{m_r+1} 
   - 2\rho_{m_r+1,m_r+2}\sigma_{m_r+1}\sigma_{m_r+2}
   +\sigma^2_{m_r+2})/4.
\end{equation*}

If such a conjecture were true, it would
imply the following, more modest marginal
result regarding a single block of such a matrix,
which without loss of generality we take to
be the first block.  Specifically, if
$d_1 = 2$ and
$\tau(r) = r$, for all $1 \le r \le m$,
we claim that
$(R_{\infty}^1,R_{\infty}^2) = (V_{\infty}^1,V_{\infty}^2-V_{\infty}^1)$
is distributed as the spectrum 
$(\lambda_1,\lambda_2)$ of the
$2 \times 2$ GUE matrix

\begin{equation}\label{rmt1}
A_1 :=
  \begin{pmatrix}
     X_1 &Y_1 + iZ_1 \\
     Y_1 - iZ_1   &X_2
  \end{pmatrix},
\end{equation}

\noindent where $\lambda_1 \ge \lambda_2$.
Equivalently, we will show that
$(V_{\infty}^1,V_{\infty}^2)$
is distributed as
$(\lambda_1,\lambda_1+\lambda_2)$.

Let the $2 \times 2$ submatrix $\Sigma_2$
of $\Sigma$ be written as

\begin{equation}\label{rmt2}
\Sigma_2 =
  \begin{pmatrix}
     \tilde{\sigma}^2_1              &\tilde{\rho}\tilde{\sigma}_1\tilde{\sigma}_2\\
     \tilde{\rho}\tilde{\sigma}_1\tilde{\sigma}_2      &\tilde{\sigma}^2_2
  \end{pmatrix}.
\end{equation}

\noindent Then

\begin{align}\label{rmt3}
(V_{\infty}^1,V_{\infty}^2)
&= \Bigl( \max_{0\le t \le 1} 
     \bigl(\tilde{\sigma}_1\tilde{B}^1(t) 
       + \tilde{\sigma}_2\tilde{B}^2(1) - \tilde{\sigma}_2\tilde{B}^2(t) \bigr),\nonumber\\
&\qquad \qquad   \tilde{\sigma}_1\tilde{B}^1(1) + \tilde{\sigma}_1\tilde{B}^2(1)  \Bigr)\nonumber\\
&= \Bigl(\tilde{\sigma}_2\tilde{B}^2(1)
     + \max_{0\le t \le 1} 
     \bigl(\tilde{\sigma}_1\tilde{B}^1(t) 
       - \tilde{\sigma}_2\tilde{B}^2(t) \bigr),\nonumber\\
&\qquad \qquad    \tilde{\sigma}_1\tilde{B}^1(1) + \tilde{\sigma}_1\tilde{B}^2(1)  \Bigr).
\end{align}

We simplify \eqref{rmt3}, by
introducing new Brownian motions and then
decomposing the resulting expression into
two independent parts.  To do so, begin by
defining the new variances and correlation coefficients
$\sigma_1^2 :=  \tilde{\sigma}^2_2$,
$\sigma_2^2 :=  \tilde{\sigma}^2_1 - 
2\tilde{\rho}\tilde{\sigma}_1\tilde{\sigma}_2 + \tilde{\sigma}^2_2$,
and
$\rho := (\tilde{\rho}\tilde{\sigma}_1 - \tilde{\sigma}_2)/
\sqrt{\tilde{\sigma}^2_1 - 
2\tilde{\rho}\tilde{\sigma}_1\tilde{\sigma}_2 + \tilde{\sigma}^2_2}$.
Then it is easily verified that 
$B^1(t) := \tilde B^2(t)$, and
$B^2(t) := (\tilde \sigma_1 \tilde B^1(t) - \tilde \sigma_2 \tilde B^2(t))/\sigma_2 $
are (dependent) standard Brownian motions, and
\eqref{rmt3} becomes

\begin{align}\label{rmt4}
(V_{\infty}^1,V_{\infty}^2)
&= \bigl(\sigma_1 B^1(1) + \sigma_2 \max_{0\le t \le 1} B^2(t),
          2\sigma_1 B^1(1) + \sigma_2 B^2(1)\bigr)\nonumber\\
&= \Bigl( (\sigma_1 B^1(1) - \rho \sigma_1 B^2(1)) 
         + \sigma_2 \Bigl(\rho \frac{\sigma_1}{\sigma_2} +  \max_{0\le t \le 1} B^2(t) \Bigr),\nonumber\\
&\qquad \qquad 2(\sigma_1 B^1(1) - \rho \sigma_1 B^2(1))  
              + (\sigma_2 + 2\rho \sigma_1) B^2(1)\Bigr).
\end{align}

\noindent Note that 
$B^1(t) - \rho B^2(t)$ is independent of
$B^2(t)$ and has variance 
$\sigma_1^2(1-\rho^2)$.
Introducing the Brownian functional

\begin{equation}\label{rmt4b}
U(\beta) = \left(\beta - \frac12\right)B^2(1) + \max_{0\le t \le 1} B^2(t),
\end{equation}

\noindent $\beta \in \bbr$,
and using 
$\sigma_1^2, \sigma_2^2$, and $\rho$ above,
\eqref{rmt4} becomes

\begin{align}\label{rmt5}
&(V_{\infty}^1,V_{\infty}^2)
\stackrel{{\cal L}}{=} \sigma_1 \sqrt{1-\rho^2}Z (1,2)
 + \biggl(\sigma_2 U\left(\frac{1}{2} - \rho \frac{\sigma_1}{\sigma_2} \right),(\sigma_2 + 2\rho \sigma_1) B^2(1) \biggr)\nonumber\\
&\quad= \frac{\tilde{\sigma}_1\tilde{\sigma}_2\sqrt{1-\tilde \rho^2}}
        {\sqrt{\tilde{\sigma}^2_1 - 2\tilde{\rho}\tilde{\sigma}_1\tilde{\sigma}_2 + \tilde{\sigma}^2_2}}
        Z(1,2)\nonumber\\
&\qquad + \biggl(\sqrt{\tilde{\sigma}^2_1 - 2\tilde{\rho}\tilde{\sigma}_1\tilde{\sigma}_2 + \tilde{\sigma}^2_2}\quad
             U\biggl(\frac{\tilde{\sigma}^2_1 - \tilde{\sigma}^2_2}
                         {2\sqrt{\tilde{\sigma}^2_1 - 2\tilde{\rho}\tilde{\sigma}_1\tilde{\sigma}_2 
                          + \tilde{\sigma}^2_2} }  \biggr),\nonumber\\    
&\qquad \qquad \qquad \qquad \qquad  2(\tilde{\sigma}^2_1 - \tilde{\sigma}^2_2) B^2(1) \biggr),
\end{align}

\noindent where $Z$ is a standard normal
random variable independent of the sigma-field
generated by $B^2$.

Turning now to the eigenvalues' distributions,
we first consider the centered, multivariate normal
random variables $(W_1,W_2)$, having covariance matrix

\begin{equation*}
  \begin{pmatrix}
     \sigma^2_1              &\rho\sigma_1\sigma_2 \\
     \rho\sigma_1\sigma_2      &\sigma^2_2
  \end{pmatrix},
\end{equation*}

\noindent and let $W_3$ and $W_4$ be two iid, 
centered, normal random variables, independent of
$(W_1,W_2)$, with variance $\sigma_2^2$.
Then it is classical that

\begin{equation*}
\Bigl(W_2, \sqrt{W_2^2 + W_3^2 + W_4^2}\Bigr) \stackrel{{\cal L}}{=}
\sigma_2\bigl(B(1),2\max_{0\le t \le 1}B(t)-B(1)\bigr),
\end{equation*}

\noindent or, equivalently,

\begin{equation}\label{rmt6}
\biggl(W_2, \beta W_2 + \frac12 \sqrt{W_2^2 + W_3^2 + W_4^2}\biggr) \stackrel{{\cal L}}{=}
\sigma_2(B(1),U(\beta)),
\end{equation}

\noindent where $B$ is a standard Brownian
motion, and $U(\beta)$, $\beta \in \bbr$,
is defined in terms of $B$, rather than
in terms of $B^2$, as in \eqref{rmt4b}.
Then consider the random variable

\begin{align}\label{rmt7}
\tilde \lambda 
&:= W_1 + \sqrt{W_2^2 + W_3^2 + W_4^2}\nonumber\\
&= \biggl(W_1 - \rho\frac{\sigma_1}{\sigma_2}\biggr) 
       + \biggl(\rho\frac{\sigma_1}{\sigma_2} + \sqrt{W_2^2 + W_3^2 + W_4^2}\biggr).
\end{align}

\noindent Using \eqref{rmt6},
and noting that the variance of the first term
in \eqref{rmt7} is $\sigma_1^2(1-\rho^2)$,
it is easy to see that

\begin{equation}\label{rmt8}
\tilde \lambda 
\stackrel{{\cal L}}{=} \sigma_1\sqrt{1-\rho^2}Z 
                         + 2\sigma_2 U\Bigl(\frac{\rho\sigma_1}{2\sigma_2}\Bigr),
\end{equation}

\noindent where $Z$ is a standard normal random variable
independent of $B$.

We now apply this result to the eigenvalues of the 
matrix $A_1$ in \eqref{rmt1}, namely, to

\begin{equation}\label{rmt9}
\lambda_1 
= \biggl(\frac{X_1 + X_2}{2}\biggr)
  + \sqrt{\biggl(\frac{X_1 - X_2}{2}\biggr) + Y_1^2 + Z_1^2}, 
\end{equation}

\noindent and

\begin{equation}\label{rmt10}
\lambda_2 
= \biggl(\frac{X_1 + X_2}{2}\biggr)
  - \sqrt{\biggl(\frac{X_1 - X_2}{2}\biggr) + Y_1^2 + Z_1^2}.
\end{equation}

\noindent Letting $W_1 = (X_1 + X_2)/2$, $W_2 = (X_1 - X_2)/2$,
$W_3 = Y_1$, and $W_4 = Z_1$, we have

\begin{align}\label{rmt11}
(\lambda_1, \lambda_1+\lambda_2) 
&= \biggl( W_1 + \sqrt{W_2^2 + W_3^2 + W_4^2}, 2W_1 \biggr)\nonumber\\
&=  \biggl( \biggl(W_1 -\hat \rho \frac{\hat \sigma_1}{\hat \sigma_2}W_2\biggr)
       + 2\biggl( \hat \rho \frac{\hat \sigma_1}{2\hat \sigma_2}W_2
                 + \frac12 \sqrt{W_2^2 + W_3^2 + W_4^2}\biggr),\nonumber\\ 
&\qquad \qquad 2\biggl(W_1 -\hat \rho \frac{\hat \sigma_1}{\hat \sigma_2}W_2\biggr)
       + 2\hat \rho \frac{\hat \sigma_1}{\hat \sigma_2}W_2\biggr)\nonumber\\
&= \biggl(W_1 -\hat \rho \frac{\hat \sigma_1}{\hat \sigma_2}W_2\biggr)(1,2)\nonumber\\
&\qquad \qquad + \biggl(\hat \rho \frac{\hat \sigma_1}{2\hat \sigma_2}W_2
                           + \frac12 \sqrt{W_2^2 + W_3^2 + W_4^2},
               2\hat \rho \frac{\hat \sigma_1}{\hat \sigma_2}W_2\biggr),
\end{align}

\noindent where
$\hat \sigma_1^2 = (\tilde{\sigma}^2_1 
+2\tilde{\rho}\tilde{\sigma}_1\tilde{\sigma}_2 + \tilde{\sigma}^2_2)/4$,
$\hat \sigma_2^2 = (\tilde{\sigma}^2_1 
-2\tilde{\rho}\tilde{\sigma}_1\tilde{\sigma}_2 + \tilde{\sigma}^2_2)/4$,
and
$\hat \rho \hat \sigma_1^2 \hat \sigma_2^2 = (\tilde{\sigma}^2_1 - \tilde{\sigma}^2_2)/4$.
Noting that the variance of 
$W_1 - (\hat \rho \hat \sigma_1/ \hat \sigma_2) W_2$
is $\hat \sigma_1^2(1-\hat \rho^2) = \sigma_1^2(1-\rho^2)$,
and that, moreover,
$\beta := \hat \rho \hat \sigma_1 / 2\hat \sigma_2 = 
(\tilde{\sigma}^2_1 - \tilde{\sigma}^2_2) /
(2\sqrt{\tilde{\sigma}^2_1 - 2\tilde{\rho}\tilde{\sigma}_1\tilde{\sigma}_2 
                          + \tilde{\sigma}^2_2})$,
we find that

\begin{align}\label{rmt12}
(\lambda_1, \lambda_1+\lambda_2) 
&= \hat \sigma_1 \sqrt{1 - \hat \rho^2}Z(1,2)
   + \Bigl(2\hat \sigma_2 U\Bigl(\frac{\hat \rho \hat \sigma_1}{2 \sigma_2}\Bigr),
      2 \hat \rho \frac{\hat \sigma_1}{\hat \sigma_2}B^2(1) \Bigr)\nonumber\\
&= \sigma_1\sqrt{1- \rho^2}Z(1,2) + \sigma_2\bigl(U(\beta), 4\beta B^2(1) \bigr)\nonumber\\
&\stackrel{{\cal L}}{=} (V_{\infty}^1,V_{\infty}^2),
\end{align}

\noindent and we have our identity in law.\\

To illustrate the ways in which random matrix interpretations
might potentially illuminate other, apparently unrelated, Brownian
functionals, consider the following example.  Let
$(\varepsilon_k)_{k\ge 1}$ be a sequence of positive numbers decreasing to
zero.  Then it is possible to find an increasing 
sequence of integers $(m_k)_{k \ge 1}$
so that, for each $k$,
there is a Markov chain on $m_k$ letters such that:

\textbullet \quad the maximal stationary probability $\pi_{max}(k)$
is of multiplicity $3$, and

\textbullet \quad the $3 \times 3$ covariance submatrix
$\Sigma_3(k)$ governing the associated Brownian functional
$V^1_{\infty}(k)$ is of the form

\begin{equation}\label{locscore1}
  \Sigma_3(k) =
  \sigma(k)^2
  \begin{pmatrix}
  \varepsilon_k^2  &0 &0\\
  0 &1 &0\\
  0 &0 &\varepsilon_k^2 
  \end{pmatrix}.
\end{equation}

\noindent That is, the variance of $B^{\tau(2)}$ becomes arbitrarily
large in comparison to that of
$B^{\tau(1)}$ and $B^{\tau(3)}$.

Then, since $LI_n(k) = V^1_n(k)$, we have, as $n \rightarrow \infty$,

{\allowdisplaybreaks
\begin{align}\label{locscore2}
&\frac{LI_n(k) - \pi_{max}(k)}{\sqrt{n}} \Rightarrow
  \max_{I_{1,3}}\sum_{\ell=1}^3 \sigma_{\tau(\ell)} B^{\tau(\ell)}(\Delta t_{\ell})\nonumber\\
&\qquad = \sigma(k) \max_{I_{1,3}} 
           \bigl(\varepsilon_k (B^{\tau(1)}(t_1) - B^{\tau(1)}(0)  )
          +               (B^{\tau(2)}(t_2) - B^{\tau(2)}(t_1))\nonumber\\
&\qquad \qquad  \qquad +(\varepsilon_k (B^{\tau(3)}(1)   - B^{\tau(3)}(t_2))\bigr)\nonumber\\
&\qquad := V^1_{\infty}(k),
\end{align}
}

\noindent so that, as $k \rightarrow \infty$,

{\allowdisplaybreaks
\begin{align}\label{locscore3}
\frac{V^1_{\infty}(k)}{\sigma(k)}
&\Rightarrow \max_{0 \le t_1 \le t_2 \le 1} (B(t_2) - B(t_1)),
\end{align}
}

\noindent where $B(t)$ is a standard Brownian motion.  
The right-hand side of \eqref{locscore3} is known as the 
{\it local score}, and describes the largest positive increase
that $B$ makes within the unit interval.
Such functionals are of great importance in sequence comparison,
particularly in bioinformatics ({\it e.g.},
see Daudin, Ettienne, and Vallois \cite{DEV}.)
Moreover, 
\begin{equation}\label{locscore4}
\max_{0 \le t_1 \le t_2 \le 1} (B(t_2) - B(t_1)) \stackrel{\cal L}{=}
\max_{0 \le t \le 1} |B(t)|,
\end{equation}

\noindent which follows immediately 
from the classical equality in law,
due to L\'evy, 
$(|B(t)|)_{t\ge 0} \stackrel{\cal L}{=} (\max_{0 \le s \le t} B(s) - B(t))_{t\ge 0}$.
Thus, if we have a random matrix connection to
$V^1_{\infty}(k)$, we can extend it to 
$\max_{0 \le t \le 1} |B(t)|$, at least in some limiting sense.
This is also interesting from the following point of
view.  
Classically, the Brownian functional 
$\max_{0 \le t \le 1} B(t) \stackrel{\cal L}{=} |B(1)|$,
and a trivial random
matrix connection can be seen
by examining the eigenvalues of the random matrix

\begin{equation}\label{locscore5}
   \begin{pmatrix}
  Z &0\\
  0 &-Z
  \end{pmatrix},
\end{equation}

\noindent where $Z$ is a standard normal
random variable.  Then, clearly,
$\lambda_{max}$ has
the half-normal law, since
$\lambda_{max} = \max(Z,-Z) = |Z|$.
Thus, the functional 
$\max_{0 \le t \le 1} B(t)$ has
a random matrix interpretation, one which is
considerably simpler than any
potential random matrix interpretation for
$\max_{0 \le t \le 1} |B(t)|$.

\section{Concluding Remarks}

In this paper, we have obtained the limiting shape of
Young tableaux generated by an aperiodic, irreducible,
homogeneous Markov chain on a finite state alphabet.
The following remarks indicate natural directions
in which our results in some cases can, 
and in other cases, may hope to, be extended.\\

\noindent \textbullet \quad  Our limiting 
theorems have all been proved assuming 
that the initial distribution is 
the stationary one. However, 
such results as Theorem 2 of Derriennic and Lin
\cite{DL} allow to extend our framework to initial
distributions started at a specified state.
Indeed, in this case, {\it i.e.,} if for
some $k=1,\dots,m$, $\bbp(X_0 = \alpha_k) = 1$,
the asymptotic covariance matrix is still given
by \eqref{item6la}, and, for example, Theorem \ref{thm4}
remains valid.  For an arbitrary initial distribution,
what is needed in this non-stationary context is
an invariance principle.  More generally, our
results continue to hold for $k^{th}$-order Markov chains,
and in fact, they extend to {\it any} sequence for which
both an asymptotic covariance matrix and an 
invariance principle exist.\\

\noindent \textbullet \quad Our limiting theorems 
have only been proved for finite alphabets.
However, from the authors' previous work \cite{HL},
it is known that for countably infinite iid alphabets, 
$LI_n$ has a limiting law corresponding to that 
of a non-uniform, finite-alphabet.  Hence,
for a countably infinite-alphabet Markov chain (subject
to additional constraints such as Harris recurrence?), 
we might still be able to obtain limiting laws 
of the form developed in this paper.\\

\noindent \textbullet \quad By using
appropriate existing concentration inequalities,
one can expect to establish
the convergence of the moments of
the rows of the tableaux.\\

\noindent \textbullet \quad 
One field in which the connection between Brownian
functionals and random matrix theory has been
exploited is in Queuing Theory.  The development
below, following O'Connell and Yor \cite{OY1},
shows how Brownian functionals of the sort we have
studied arise as generalizations of standard queuing
models.

Let $A(s,t]$ and $S(s,t]$,
$-\infty < s < t < \infty$,
be two independent Poisson
point process on $\bbr$, with intensity measures
$\lambda$ and $\mu$, respectively, with $0 < \lambda < \mu$.
Here $A$ represents the {\it arrivals} process, and $S$ 
the {\it service time} process, at a queue consisting
of a single server. The condition $\lambda < \mu$ 
ensures that the {\it queue length}

\begin{equation}\label{qitem1}
Q(t) = \sup_{-\infty < s \le t} \left\{ A(s,t] - S(s,t] \right\},
\end{equation}

\noindent is a.s.~finite, for any $t \in \bbr$.  
Then, defining the {\it departure} process

\begin{equation}\label{qitem2}
D(s,t] = A(s,t] - (Q(t)-Q(s)),
\end{equation}

\noindent which is simply the number of arrivals
during $(s,t]$ less the change in the queue length
during $(s,t]$,
the classical problem is to determine the distribution
of $D(s,t]$.  The answer to this problem is given
by {\it Burke's Theorem} \cite{Burke} (see Theorem $1$ of \cite{OY1}):

\begin{theorem}\label{Bkthm}
$D$ is a Poisson process with intensity $\lambda$, and
$\{D(s,t],s\le t\}$ is independent of $\{Q(s),s\ge t\}$.
\end{theorem}

That is, $D$ has the same law as the arrivals process $A$.
Moreover, since,the queue length after time $t$ is independent
of the process $D$ up to time $t$, 
one may take the departures from the first queue
and use them as inputs to a second queue, and observe that
the departure process from the second queue also has the law
of $A$.  Proceeding in this way, one generalizes to a {\it tandem queue} 
of $n$ servers, each taking the departures from the previous queue
as its arrivals process.

One can further generalize this model to a {\it Brownian queue
in tandem} in the following manner.  Let
$B,B^1,B^2,\dots, B^n$ be independent, standard Brownian motions
on $\bbr$,
and write $B^k(s,t) = B^k(t) - B^k(s)$, for each $k$ and $s < t$, 
and similarly for $B$.  Let $m > 0$ be a constant, and define,
in complete analogy to \eqref{qitem1} and \eqref{qitem2},

\begin{equation}\label{qitem3}
q_1(t) = \sup_{-\infty < s \le t} \left\{ B(s,t) + B^1(s,t) - m(t-s) \right\},
\end{equation}

\noindent and, for $s < t$,

\begin{equation}\label{qitem4}
d_1(s,t) = B(s,t) - (q_1(t)-q_1(s)).
\end{equation}

\noindent  For $k = 2,3,\dots,n$, let

\begin{equation}\label{qitem5}
q_k(t) = \sup_{-\infty < s \le t} \left\{ d_{k-1}(s,t) + B^k(s,t) - m(t-s) \right\},
\end{equation}

\noindent and, for $s < t$,

\begin{equation}\label{qitem6}
d_k(s,t) = d_{k-1}(s,t) - (q_k(t)-q_k(s)).
\end{equation}

Here $B$ is the arrivals process for the first queue, 
$d_{k-1}$ is the arrivals process for the $k^{th}$ queue ($k \ge 2)$,
and  $mt - B^k(t)$ is the service process for the $k^{th}$ queue,
for all $k$.
Using the ideas employed in Burke's Theorem, 
it can be shown that the generalized queue lengths
$q_1(0), q_2(0), \dots, q_n(0)$
are iid random variables.  Moreover, they are 
exponentially distributed with mean $1/m$.

Using the definitions in \eqref{qitem3}-\eqref{qitem6},
and a simple inductive argument, one finds that

\begin{equation}\label{qitem7}
\sum_{k=0}^n q_k(0) 
=\sup_{t > 0} 
        \Bigl\{ B(-t,0) - mt + L_n(t)\Bigr\},
\end{equation}

\noindent where

\begin{equation}\label{qitem8}
L_n(t) = \sup_{\stackrel{\scriptstyle 0\le s_1\le\cdots}{\le s_{m-1}\le t} }
                 \{B^1(-t,-s_{n-1}) + \cdots + B^n(-s_1,0)\}.
\end{equation}

By Brownian rescaling, we observe that

\begin{align}\label{qitem9}
L_n(t) 
&\stackrel{ {\cal L} }{=} \sqrt{t}\sup_{\stackrel{\scriptstyle 0\le s_1\le\cdots}{\le s_{m-1}\le 1} }
                 \{B^1(-1,-s_{n-1}) + \cdots + B^n(-s_1,0)\}\nonumber\\
&\stackrel{ {\cal L} }{=} \sqrt{t}V^1_{\infty},
\end{align}

\noindent where the functional $V^1_{\infty}$ is
as in Theorem \ref{thm4}, with associated $n \times n$
covariance matrix $\Sigma = tI_n$ and parameter set
$I_{1,n}$.  Thus, $L_n(t)$ may
be thought of as a process version of this $V^1_{\infty}$.

The generalized Brownian queues 
in \eqref{qitem3}-\eqref{qitem6}
involved
{\it independent} Brownian motions.  These can 
be extended with Brownian motions 
$B^1$,
$\dots $,
$B^n$ for which
$(\sigma_1 B^1(t),\dots \sigma_n B^n(t))$
has (nontrivial) covariance matrix $t\Sigma$.  
Whether or not we keep the initial 
arrival process $B(t)$
independent of $(B^1$,$\dots$,$B^n)$, we now
no longer have that
$q_1(0), q_2(0), \dots, q_n(0)$
are iid random variables, 
due to the dependence among the service times
$mt-B^k(t)$, but we do still have the identity
\eqref{qitem7} and \eqref{qitem9} 
relating the total occupancy
of the queue at time zero to 
$V_{\infty}^1$.
More importantly, our generalizations of the Brownian
functionals $L_n(t)$ above can be used to describe
the joint law of the input/output of each queue.\\

\noindent \textbullet \quad An important topic connecting
much of random matrix theory to other problems, such
as the shape of random Young tableaux, is the field of
orthogonal polynomials. (See, {\it e.g.}, \cite{Jo}.)
It would be of great interest to see what, if any, classes of
orthogonal polynomials are associated with the present paper.\\


\begin{thebibliography}{9} 

\bibitem{BDJ2} Baik, J. Deift, P., and Johansson, K.,  
``On the distribution of the length of the second row of a Young diagram under Plancherel measure,''
{\it Geom. Funct. Anal.}, vol. 10, no. 4, pp. 702-731, 2000.  


\bibitem{BDJ2Add} Baik, J. Deift, P., and Johansson, K.,  
''Addendum to: ``On the distribution of the length of the second row of a Young diagram under Plancherel measure'','' 
{\it Geom. Funct. Anal.}, vol. 10, no. 6, pp. 1606-1607, 2000.  

\bibitem{BDJ1} Baik, J., Deift, P., and Johansson, K.,
``On the distribution of the length of the longest increasing subsequence of random permutations,''
{\it J. Amer. Math. Soc.}, vol. 12, no. 4, pp. 1119-1178, 1999.

\bibitem{Ba} Baryshnikov, Y.,
``GUEs and queues,'' 
{\it Probab. Theory Related Fields}, vol. 119, no. 2, pp. 256-274, 2001.

\bibitem{Bill} Billingsley, P.,
{\it Convergence of probability measures}.
Wiley Series in Probability and Statistics: Probability and Statistics, New York:
John Wiley and Sons Inc., second ed., 1999.
A Wiley-Interscience Publication.

\bibitem{BOO} Borodin, A., Okounkov, A., and Olshanski, G,
``Asymptotics of Plancherel measures for symmetric groups,''
{\it J. Amer. Math. Soc.}, vol. 13, no. 3, pp. 481-515 (electronic), 2000.

\bibitem{BJ} Bougerol, P. and Jeulin, T.,
``Paths in Weyl chambers and random matrices,'' 
{\it Probab. Theory Related Fields}, vol. 124, no. 4, pp. 517-543, 2002.

\bibitem{Burke} Burke, P.J.,
``The output of a queueing system,''
{\it Operations. Res.}, vol. 4, pp. 699-704 (1957), 1956.

\bibitem{ChG} Chistyakov, G.P. and G\"otze, F.,
``Distribution of the shape of Markovian random words,''
{\it Probab. Theory Related Fields}, vol. 129, no. 1, pp. 18-36, 2004.

\bibitem{DEV} Daudin, J.-J., Etienne, M.P., and Vallois, P.,
``Asymptotic behavior of the local score of independent and identically distributed random sequences,''
{\it Stochastic Process. Appl.},  vol. 107, no. 1, pp. 1-28, 2003.

\bibitem{DL} Derriennic, Y. and Lin, M.,
``The central limit theorem for Markov chains with normal transition operators, started at a point,''
{\it Probab. Theory Related Fields}, vol. 119, no. 4, pp. 508-528, 2001.

\bibitem{Do} Doumerc, Y.,
``A note on representations of eigenvalues of classical Gaussian matrices,'' 
in {\it S\'eminaire de Probabilit\'es XXXVII}, vol. 1832 of
{\it Lecture Notes in Math.}, pp. 370-384, Berlin: Springer, 2003.

\bibitem{Fu} Fulton, W.,
{\it Young tableaux}, vol. 35 of {\it London Mathematical Society Student Texts.}
Cambridge: Cambridge University Press, 1997.
With applications to representation theory and geometry.

\bibitem{GW}  Glynn, P.W. and Whitt, W.,
``Departure from many queues in series,''
{\it Ann. Appl. Probab.}, vol. 1, no. 4, pp. 546-572, 1991.

\bibitem{Go} Gordin, M.I.,  
``The central limit theorem for stationary processes,'' 
{\it Dokl. Akad. Nauk SSSR}, vol. 188, pp. 739-741, 1969.

\bibitem{GTW} Gravner, J., Tracy, C., and Widom, H.,
``Limit theorems for height fluctuations in a class of discrete space and time growth models,''
{\it J. Statist. Phys.}, vol. 102, no. 5-6, pp. 1085-1132, 2001.


\bibitem{Greene} Greene, C.,
``Some partitions associated with a partially ordered set,''
{\it J. Combin. Theory Ser. A}, vol. 20, no. 1, pp. 69-79, 1976.


\bibitem{HL} Houdr\'e C. and Litherland, T.,
``On the longest increasing subsequence for finite and countable alphabets.'' 
(Preprint: arXiv:math/0612364)

\bibitem{HLM} Houdr\'e C., Lember, J., and Matzinger, H.,
``On the longest common increasing binary subsequence,''
{\it C. R. Acad. Sci. Paris}, vol. 343, no. 9, pp. 589-594, 2006.

\bibitem{ITW1} Its, A.R., Tracy, C. A., and Widom, H.,  
``Random words, Toeplitz determinants, and integrable systems. I,''
in {\it Random matrix models and their applications}, vol. 40 of
{\it Math. Sci. Res. Inst. Publ.}, pp. 245-258, Cambridge: Cambridge Univ. Press, 2001.

\bibitem{ITW2} Its, A.R., Tracy, C. A., and Widom, H.,  
``Random words, Toeplitz determinants, and integrable systems. II,''
{\it Phys. D}, vol 152/153, pp. 199-224, 2001.
Advances in Nonlinear Mathematics and Science.

\bibitem{Jo} Johansson, K.,
``Discrete orthogonal polynomial ensembles and the Plancherel measure,''
{\it Ann. of Math. (2)}, vol. 153, no.1, pp. 259-296, 2001.

\bibitem{Jo2} Johansson, K.,  
``Shape fluctuations and random matrices,'' 
{\it Comm. Math. Phys}, vol. 209, no. 2, pp. 437-476, 2000.


\bibitem{Ku} Kuperberg, G.,
``Random words, quantum statistics, central limits, random matrices,''
{\it Methods. Appl. Anal.}, vol. 9, no. 1, pp. 99-118, 2002.


\bibitem{OY1} O'Connell, N. and Yor, M., 
``Brownian analogues of Burke's theorm,'' 
{\it Stochastic Process. Appl.}, vol. 96, no. 2, pp. 285-304, 2001.

\bibitem{OY2} O'Connell, N. and Yor, M., 
``A representation for Non-colliding random walks,''
{\it Electron. Comm. Probab.}, vol. 7, pp. 1-12 (electronic), 2002.

\bibitem{Ok} Okounkov, A.,
``Random matrices and random permutations,'' 
{\it Internat. Math. Res. Notices.}, , no. 20, pp. 1043-1095, 2000.

\bibitem{TW}  Tracy, C. A. and Widom, H.,  
``On the distributions of the lengths of the longest monotone subsequences in random words,''  
{\it Probab. Theory Related Fields}, vol. 119, no. 3, pp. 350-380, 2001.


\end{thebibliography}
\end{document}